\documentclass{article}
\usepackage[a4paper]{geometry}

\usepackage{amsmath}
\usepackage{amssymb}
\usepackage[utf8]{inputenc} 
\usepackage{xcolor}
\usepackage{subcaption}
\usepackage{gnuplot-lua-tikz}
\usepackage{booktabs}
\usepackage{cite}

\usepackage{algorithm}
\usepackage{algorithmic}

\usepackage{hyperref}

\title{Performance of the low-rank tensor-train SVD (TT-SVD) for large dense tensors on modern multi-core CPUs}
\author{Melven Röhrig-Zöllner \and Jonas Thies \and Achim Basermann}

\begin{document}
\maketitle

\begin{abstract}
There are several factorizations of multi-dimensional tensors into lower-dimensional components, known as `tensor networks'.
We consider the popular `tensor-train' (TT) format and ask: How efficiently can we compute a low-rank approximation from a full tensor on current multi-core CPUs?
Compared to sparse and dense linear algebra, kernel libraries for multi-linear algebra are rare and typically not as well optimized.
Linear algebra libraries like BLAS and LAPACK may provide the required operations in principle, but often at the cost of additional data movements for rearranging memory layouts.
Furthermore, these libraries are typically optimized for the compute-bound case (e.g.\ square matrix operations) whereas low-rank tensor decompositions lead to memory bandwidth limited operations.
We propose a `tensor-train singular value decomposition' (TT-SVD) algorithm based on two building blocks:
a `Q-less tall-skinny QR' factorization, and a fused tall-skinny matrix-matrix multiplication and reshape operation.
We analyze the performance of the resulting TT-SVD algorithm using the Roofline performance model.
In addition, we present performance results for different algorithmic variants for shared-memory as well as distributed-memory architectures.
Our experiments show that commonly used TT-SVD implementations suffer severe performance penalties.
We conclude that a dedicated library for tensor factorization kernels would benefit the community:
Computing a low-rank approximation can be as cheap as reading the data twice from main memory.
As a consequence, an implementation that achieves realistic performance will move the limit at which one has to resort to randomized methods that only process part of the data.
\end{abstract}

\clearpage

\section{Introduction}
The tensor-train (TT) decomposition is a particular form of a tensor network representation of a high-dimensional tensor
in which the 3D `core tensors' are aligned in a 1D format and connected by a contraction with their direct neighbors only to represent (or approximate)
a $d$-dimensional tensor. It was introduced as such by Oseledets~\cite{oseledets2009a,oseledets2009b}), but in fact has been known to (and used by) computational physicists under the name of Matrix Product States (MPS) since the 1980s~\cite{affleck1985,affleck1987}; see also~\cite{verstraete2006mps} for a more recent reference.
Closely related is the Density Matrix Renormalization Group (DMRG) algorithm~\cite{white1992dmrg}, an optimization method that operates on the space of MPS.
An overview on numerical algorithms based on low-rank tensor approximations can be found in~\cite{Grasedyck2013}.
Recent research also focuses on applications of tensor-trains in data science, see e.g.\ \cite{Klus2019, Chen2019, Chen2022, Kour2020} for a few examples.
The performance of common arithmetic operations in tensor-train format (such as additions and scalar products) are discussed in~\cite{Daas2022}.

One can construct an approximate TT-decomposition of high-dimensional data $X\in\mathbb{R}^{n_1 \times n_2 \times \dots \times n_d}$ using a high order singular value decomposition.
An algorithm for this, called TT-SVD, is presented in~\cite{Oseledets2011}.
Given $X$ and a maximum `bond dimension' $r_{max}$, it successively determines the core tensors $T^{(j)} \in\mathbb{R}^{r_{j-1} \times n_j \times r_{j}}, j=1\dots d$,
such that       $r_0=r_d=1$, 
                $r_j\le r_{max}$,
                the rows or columns of some matricization of all but one $T^{(j)}$ are orthonormal, 
                and the approximation error (difference between the tensor-train formed by the $T^{(j)}$ and the original data tensor $X$)
                is minimized (up to a constant factor) in the Frobenius norm on the manifold of rank-$r_{max}$ tensor-trains~\cite{Oseledets2011}.
Definitions of some of these concepts are obviously needed, and will be given in \autoref{sec:background}.

The aim of this paper is to develop an efficient formulation and implementation of this algorithm for modern multi-core CPUs.
We focus on situations where the data is large and dense, but it is feasible to process the complete data set for which a low-rank representation is sought (i.e., to read the data $\mathcal{O}(1)$ times).
In contrast, randomized (sampling) algorithms only access part of the data and can be used if the data set is too large~\cite{Larsen2020,Kolda2020}.
For the deterministic case, error bounds and asymptotic complexity estimates (for the size of the result) exist but differ slightly depending on the desired tensor format, see~\cite{Grasedyck2013} and the references therein.
One usually seeks an approximation with a specific accuracy (in terms of maximal size of the resulting approximation or a tolerance, or both).
However, common implementations often provide sub-optimal performance for this case as they do not take into account that the computation is limited by data transfers on current computers (see \autoref{sec:numerical_experiments}).
We investigate the TT-SVD because this is a simple and popular choice, but the ideas can be transferred to other tree tensor networks (see e.g.\ \cite{Grasedyck2011}) as the algorithmic building blocks are similar.
An important ingredient in our implementation is a Q-less `tall and skinny QR' (TSQR, see~\cite{Demmel2012}) variant that is described in detail in \autoref{sec:TSQR_algorithm}.
The idea to avoid computing and storing the large matrix $Q$ of a QR decomposition was already exploited for e.g.\ sparse matrix decompositions and tensor calculus in~\cite{Berry2005, Fan2016}.

Our contribution is twofold.
First, based on the example of the TT-SVD algorithm we show that low-rank tensor approximation is a memory-bound problem in high dimensions (in contrast to the SVD in two dimensions for square matrices).
Second, we discuss how the TT-SVD algorithm can be implemented efficiently on current hardware.
In order to underline our findings, we present performance results for the required building blocks and for different TT-SVD variants and implementations on a small CPU cluster.

The remainder of this paper is organized as follows.
In \autoref{sec:background}, we introduce the basic concepts and notation for tensor networks and performance engineering that we will use to describe our algorithms and implementation.
In \autoref{sec:algorithm} we describe the TT-SVD algorithm with focus on our tailored Q-less TSQR variant.
In \autoref{sec:performance_analysis} we present a performance model for the two key components of TT-SVD (Q-less TSQR and a `tall and skinny' matrix-matrix multiplication), as well as the overall algorithm.
Numerical experiments comparing actual implementations of TT-SVD (including our own optimized version) against the performance model can be found in \autoref{sec:numerical_experiments},
and the paper closes with a summary of our findings in \autoref{sec:conclusion}.

\section{Background and notation}\label{sec:background}
\subsection{Tensor notation and operations}
Classical linear algebra considers matrices and vectors ($n\times1$ matrices) and provides a notation for operations between them based on matrix-matrix products and matrix transpositions.
We make use of this common notation where possible. 
In this paper, a dense $d$-dimensional array or tensor is denoted by $X \in \mathbf R^{n_1 \times \dots \times n_d}$.
We can combine and split dimensions through $\operatorname{reshape}$ operations, e.g.:
\begin{align*}
Y &= \operatorname{reshape}\left(X, \begin{pmatrix} n_1 & \frac{\bar n}{n_1 n_d} & n_d\end{pmatrix}\right) \; \in \mathbf R^{n_1 \,\times \,\bar n/(n_1n_d)\,\times\, n_d}, \quad\text{with } \bar n := \prod_{i=1}^d n_i,\nonumber\\
X &= \operatorname{reshape}\left(Y, \begin{pmatrix} n_1 & \dots & n_d\end{pmatrix}\right).
\end{align*}
This assumes that the dimensions of a tensor are ordered and provides a notation for unfolding a d-dimensional tensor into a lower-dimensional tensor, respectively into a matrix (matricization), and folding it back into a d-dimensional tensor.
It only allows us to combine neighboring dimensions, which is sufficient for all cases in this paper.
In practice, many tensor algorithms can be written as series of matrix-operations of different matricizations of tensors,
but more general $\operatorname{reshape}$ operations can often be implemented without overhead by just reinterpreting the data in memory.

\subsubsection{Matrix decompositions}
In two dimensions, the singular value decomposition defines the (unique) decomposition of a rectangular matrix $M\in\mathbf R^{n_1\times n_2}$,
\begin{align}
M &= U \Sigma V^T &\Leftrightarrow&& M_{i_1,i_2} &= \sum_{j=1}^r U_{i_1,j}\, \sigma_j\, V_{i_2,j}
\end{align}
into the orthonormal matrices of left and right singular vectors $U \in \mathbf R^{n_1 \times r}$, $U^TU=I$ and $V\in\mathbf R^{n_2\times r}$, $V^TV=I$ and a diagonal $\Sigma = \operatorname{diag}(\sigma_1, \dots, \sigma_r)$ with singular values $\sigma_1 \ge \sigma_2 \ge \dots \ge \sigma_r > 0$.
The decomposistion is unique if $\sigma_1 > \sigma_2 > \dots > \sigma_r > 0$.
The rank of the matrix is defined as $r=\mathrm{card}(\{\sigma_j>0\}) \le \min(n_1, n_2)$.

In the steps of the TT-SVD algorithm, we also use the QR decomposition
\begin{align}
M &= QR,
\end{align}
with an orthonormal matrix $Q \in \mathbf R^{n_1 \times n_2}$, $Q^TQ=I$ and an upper triangular matrix $R\in\mathbf R^{n_2\times n_2}$ and $n_1 \ge n_2$.

\subsubsection{Tensor-train decomposition}
The tensor-train (TT) decomposition introduced in \cite{Oseledets2011} generalizes the idea of the SVD to $d$ dimensions:
\begin{align}
X_{i_1,i_2,\dots,i_d} &= \sum_{j_1=1}^{r_1} \sum_{j_2=1}^{r_2} \cdots \sum_{j_{d-1}=1}^{r_{d-1}}
T^{(1)}_{1,i_1,j_1} \, T^{(2)}_{j_1,i_2,j_2} \, \cdots \, T^{(d)}_{j_{d-1},i_d,1}.
\end{align}
Here, the three-dimensional tensors $T^{(j)}$ are called `core tensors' of the decomposition and $r_1,\dots,r_{d-1}$ the ranks.
In contrast to the SVD, the TT decomposition is not unique but a best-approximation with given maximal rank $r_\text{max} \ge r_j$ exists and the TT-SVD algorithm in \autoref{sec:ttsvd} calculates a quasi-optimal solution.
For a detailed discussion, we refer to~\cite{Oseledets2011}.

\subsection{Performance characteristics on current hardware}
Supercomputers consist of a set of compute nodes that are connected by a network (see e.g.\ \cite{Hager2010}).
For the performance modeling, we concentrate on the node-level performance of the required algorithmic building blocks.
However, we also show results with a distributed memory variant of the TT-SVD algorithm that allows scaling beyond a single node.
Our algorithmic choices and performance optimizations are motivated by hardware characteristics of multi-core processors, which we therefore briefly introduce.

Each compute node has one or several multi-core CPU sockets with dedicated memory.
The CPU cores can access the memory of the complete node but accesses to the dedicated memory of the socket are faster (ccNUMA architecture).
To reduce the complexity of the shared memory parallelization, we use OpenMP for parallelizing over the cores of one socket, and MPI for communicating between sockets and nodes.

An important aspect of multi-core optimization is the increasing gap between the memory bandwidth and the floating point performance.
To alleviate this problem, multiple levels of caches are used, where the larger and slower levels are shared between multiple cores.
Efficient algorithms need to exploit spatial and temporal locality (accessing memory addresses close to each other and accessing the same memory address multiple times).
In addition, the floating point performance increased due to specialized wider SIMD units as well as optimized out-of-order execution of pipelined instructions.
So algorithms can only achieve high performance if they contain many independent instructions for contiguous chunks of data (e.g.\ no data dependencies/conditional branches). 

The actual run-time of a program on a specific hardware may be determined by many factors.
Therefore, it is helpful to model the performance based on a few simple characteristics that are anticipated to be potential bottlenecks.
For our numerical application, we use the Roofline performance model~\cite{Williams2009}, which considers two limiting factors.
The algorithm is either compute-bound (limited by the floating point rate) or bandwidth-bound (limited by data transfers). The upper bound for the performance is thus given by
\begin{align}
P_\text{Roofline} = \min\left(P_\text{max}, I_c b_s\right). \label{eq:Roofline}
\end{align}
Here $P_\text{max}$ and $b_s$ characterize the hardware:
$P_\text{max}$ denotes the applicable peak performance. That is, the maximal performance possible when executing the required floating point operations.
$b_s$ is the obtainable bandwidth on the slowest data path (e.g.\ from the cache or memory that is large enough to contain all data).
The bandwidth depends on the access pattern, so we need to measure it with a micro-benchmark that reflects the access pattern of the algorithm, see \autoref{tab:peak_bandwidth_and_flops}.
The algorithm is characterized by its \emph{compute intensity} $I_c$, which specifies the number of floating point operations per transferred byte.
Of course, the Roofline model is a simplification: in particular, it assumes that data transfers and calculations overlap, which is not realistic if the compute intensity is close to $P_\text{max}/b_s$.
However, the model provides insight into the behavior of the algorithm, and it allows us to assess if a specific implementation achieves a reasonable fraction of the possible performance.

\begin{table}
\centering
\begin{tabular}{lc}
benchmark & measurement \\
\midrule
memory bandwidth (pure load) & 93 GByte/s \\
memory bandwidth (copy) & 70 GByte/s \\
memory bandwidth (stream \cite{McCalpin1995}) & 73 GByte/s \\
memory bandwidth (pure store) & 45 GByte/s \\
double precision performance (AVX512 FMA) & 1009 GFlop/s
\end{tabular}
\caption{Hardware characteristics of a 14-core Intel Xeon Scalable Processor \emph{Skylake} Gold 6132. The data was measured using likwid-bench~\cite{Treibig2010} (version~5.0.1) on a single socket of a 4-socket node. All memory benchmarks use non-temporal stores and AVX512 and an array size of 1~GByte. For this system, the load bandwidth is approximately twice the store bandwidth.}\label{tab:peak_bandwidth_and_flops}
\end{table}

To analyze an algorithm in this paper, we usually first estimate the compute intensity $I_c$ and decide whether the algorithm is compute-bound or memory-bound (limited by main memory bandwidth) on the given hardware.
\begin{align}
I_c \approx \frac{n_\text{flops}}{V_\text{read+write}}.
\end{align}
Then, we calculate the ideal run-time $t_\text{min}$ from the number of floating point operations $n_\text{flops}$, respectively from the main memory data transfer volume $V_\text{read+write}$:
\begin{equation}
t_\text{min} = \left \{
\begin{array}{llll}
\frac{n_\text{flops}}{P_\text{max}}         & \text{if} & I_c > \frac{P_\text{max}}{b_s} & \text{(compute-bound)}\\
\frac{V_\text{read+write}}{b_s}  & \text{if} & I_c < \frac{P_\text{max}}{b_s} & \text{(memory-bound)}
\end{array} \right .
\end{equation}
The quotient $\frac{P_\text{max}}{b_s}$ is called \emph{machine intensity}.

Many supercomputers nowadays also feature accelerator hardware such as GPUs.
We decided not to exploit GPUs in this paper because the TT-SVD accesses the complete data, which typically does not fit into the high bandwidth memory of the device.
The slowest data path then is the PCI/e bus, which would make even the most optimized GPU implementation slower than our CPU code.

\section{Numerical algorithms and required building blocks}\label{sec:algorithm}
In this section we discuss different variants of the TT-SVD algorithm from~\cite{Oseledets2011}.
We focus on algorithmic choices required for an efficient implementation on current hardware that retain numerical accuracy and robustness.
As an important building block, we present a Q-less rank-preserving QR implementation for tall-skinny matrices (Q-less TSQR) based on~\cite{Demmel2015}.

\subsection{TT-SVD}\label{sec:ttsvd}
Based on the original TT format~\cite{Oseledets2011}, several other formats have been suggested, such as the QTT format (see e.g.\ \cite{Khoromskij2011} and the references therein), and the QTT-Tucker format~\cite{Dolgov2013}.
These formats have interesting numerical properties, however, the required operations for calculating a high-order SVD from dense data in these formats are similar.
For simplicity, we describe the algorithm for the TT format,
although it is important that the individual dimensions are as small as possible (e.g.\ two as in the QTT format) to obtain high performance.
For other hierarchical formats such as the $\mathcal H$-Tucker format (see e.g.\ \cite{Grasedyck2011}), the rank is defined differently,
so the complexity analysis of the algorithm is specific to the (Q)TT format.
The algorithmic principles and required building blocks still remain similar for high-order decomposition algorithms for other tree tensor network formats.

\subsubsection{Original TT-SVD algorithm}
We first show how the original TT-SVD algorithm from~\cite{Oseledets2011} can be implemented, see \autoref{alg:TT_SVD}
For ease of notation, we start with dimension $n_d$ (right-most core tensor in the TT-format).
\begin{algorithm}
\caption{TT-SVD}\label{alg:TT_SVD}
\begin{algorithmic}[1]
\REQUIRE $X\in\mathbf R^{n_1 \times \dots \times n_d}$, max. TT-rank $r_\text{max}\ge1$, tolerance $\epsilon$
\ENSURE TT decomposition $\sum_{j_1,\dots,j_{d-1}} T^{(1)}_{1,i_1,j_1} \, T^{(2)}_{j_1,i_2,j_2} \, \cdots \, T^{(d)}_{j_{d-1},i_d,1}= \tilde X_{i_1,\dots,i_d}$ with $\|X-\tilde X\|_F \le \epsilon \|X\|_F$ if $r_\text{max}\ge r_\delta^{(i)}$
\STATE $\delta \gets \frac{\epsilon}{\sqrt{d-1}}\|X\|_F$\hfill\COMMENT{truncation parameter}
\STATE $W \gets X$\hfill\COMMENT{temporary tensor}
\STATE $\bar n \gets \prod_{i=1}^d n_i$\hfill\COMMENT{total size of $W$}
\STATE $r_d \gets 1$
\FOR{$i = d, \dots, 2$}
  \STATE $W \gets \operatorname{reshape}\left(W, \begin{pmatrix} \frac{\bar n}{n_i r_i} & n_i r_i)\end{pmatrix}\right)$
  \STATE Calculate SVD: $U \Sigma V^T = W$ with $\Sigma = \operatorname{diag}(\sigma_1, \dots, \sigma_{n_i r_i})$
  \STATE Choose rank $r_{i-1} = \min(r_\text{max}, r_\delta^{(i)})$,\; $r_\delta^{(i)} = \min\left(j: \sigma_{j+1}^2 + \sigma_{j+2}^2 + \dots \le \delta^2\right)$
  \STATE $T^{(i)} \gets \operatorname{reshape}\left((V_{:,1:r_{i-1}})^T, \begin{pmatrix} r_{i-1} & n_i & r_i\end{pmatrix}\right)$
  \STATE $\bar n \gets \frac{\bar n r_{i-1}}{n_i r_i}$\hfill\COMMENT{new total size of $W$}
  \STATE $W \gets U_{:,1:r_{i-1}} \operatorname{diag}(\sigma_1, \dots, \sigma_{r_{i-1}})$
\ENDFOR
\STATE $T^{(1)} \gets \operatorname{reshape}\left(W, \begin{pmatrix} 1 & n_1 & r_1 \end{pmatrix}\right)$
\end{algorithmic}
\end{algorithm}
The idea of the algorithm is as follows:
Each core tensor is built subsequently from the singular vectors of a truncated SVD of a matricization (first $(i-1)$ times last $(d-i+1)$ dimensions) of the input/intermediate tensor.
In addition, the truncated directions are also removed from the input tensor for subsequent steps.
If $r_\text{max}$ is big enough, the decomposition approximates the input tensor up to the desired accuracy $\epsilon$.
Otherwise, it is less accurate (an a posteriori error bound can be calculated from the truncated singular values of each iteration).
The costly operations in this algorithm are computing the SVD in line~7, and evaluating the reduced matrix $W$ for the next iteration in line~10.
And, depending on the implementation, the reshape operation in line~6 might require copying or reordering the data in memory.
In this algorithm, the total size $\bar n$ of the work matrix $W$ is reduced in each step by a factor $\frac{r_{i-1}}{n_i r_i} \le 1$.
And $W$ is reshaped to very tall and skinny matrices in line~6 except for the last iterations, where $W$ is much smaller due to the reduction in size $\bar n$ in each step.
Therefore, it is advisable to apply the QR trick for calculating the SVD:
\begin{align}\label{eq:QR_trick}
W &= U\Sigma V^T &\Leftrightarrow&& W&=QR, \; R = \bar U \Sigma V^T \text{ with } U = Q \bar U.
\end{align}
This idea has been discussed in the literature in similar settings, see e.g.\ \cite{Constantine2014}, but we can exploit some specific details here.

One can also start the iteration in the middle of the tensor-train by reshaping $W$ into an (almost) square matrix of size approximately $\sqrt{\bar n}\times \sqrt{\bar n}$ and splitting it with an SVD into two independent iterations for a left and a right part.
This approach is not advisable because it requires $O({\bar n}^\frac{3}{2})$ floating-point operations in contrast to $O(\bar n^{1+\frac{1}{d}})$ operations for algorithms that start at the boundaries of the tensor-train (see \autoref{sec:TT-SVD_perf_analysis}).

\subsubsection{Optimized TT-SVD algorithm using TSQR}\label{sec:TSQR_TT-SVD_algorithm}
\autoref{alg:TSQR_TT_SVD} is based on the original TT-SVD (\autoref{alg:TT_SVD}) but avoids some unnecessary computations and data transfers.
It has the same numerical properties as \autoref{alg:TT_SVD} if all required matrix operations are performed accurately: QR and SVD decompositions and multiplications with orthogonal matrices.
\begin{algorithm}
\caption{Optimized TSQR TT-SVD}\label{alg:TSQR_TT_SVD}
\begin{algorithmic}[1]
\REQUIRE $X\in\mathbf R^{n_1 \times \dots \times n_d}$ stored in suitable memory layout, max. TT-rank $r_\text{max}\ge1$, tolerance $\epsilon$
\ENSURE TT decomposition $\sum_{j_1,\dots,j_{d-1}} T^{(1)}_{1,i_1,j_1} \, T^{(2)}_{j_1,i_2,j_2} \, \cdots \, T^{(d)}_{j_{d-1},i_d,1}= \tilde X_{i_1,\dots,i_d}$
\STATE $\bar n \gets \prod_{i=1}^d n_i$
\STATE $r_d \gets 1$
\STATE $W^{(d)} \gets \operatorname{reshape}\left(X, \begin{pmatrix} \frac{\bar n}{n_d} & n_d \end{pmatrix}\right)$ \hfill \COMMENT{only creates a view of $X$}
\FOR{$i = d, \dots, 2$}
  \STATE Calculate $R$ from the QR decomposition: $QR = W^{(i)}$
  \STATE Calculate small SVD: $\bar U \Sigma V^T = R$ with $\Sigma = \operatorname{diag}(\sigma_1, \dots, \sigma_{n_i r_i})$
  \STATE In the first iteration: $\delta \gets \frac{\epsilon}{\sqrt{d-1}}\|\Sigma\|_F$
  \STATE Choose rank $r_{i-1} = \min(r_\text{max}, r_\delta^{(i)})$,\; $r_\delta^{(i)} = \min\left(j: \sigma_{j+1}^2 + \sigma_{j+2}^2 + \dots \le \delta^2\right)$
  \STATE $T^{(i)} \gets \operatorname{reshape}\left((V_{:,1:r_{i-1}})^T, \begin{pmatrix} r_{i-1} & n_i & r_i\end{pmatrix}\right)$
  \STATE $\bar n \gets \frac{\bar n r_{i-1}}{n_i r_i}$
  \STATE $W^{(i-1)} \gets \operatorname{reshape}\left(W^{(i)} V_{:,1:r_{i-1}}, \begin{pmatrix} \frac{\bar n}{n_{i-1} r_{i-1}} & n_{i-1}r_{i-1}\end{pmatrix}\right)$
\ENDFOR
\STATE $T^{(1)} \gets \operatorname{reshape}\left(W^{(1)}, \begin{pmatrix} 1 & n_1 & r_1 \end{pmatrix}\right)$
\end{algorithmic}
\end{algorithm}
In the following, we discuss the three main differences between \autoref{alg:TT_SVD} and \autoref{alg:TSQR_TT_SVD}:
First, an obvious optimization is to calculate the truncation parameter $\delta$ from $\Sigma$ in the first iteration (line~7 in \autoref{alg:TSQR_TT_SVD}).
This avoids calculating the norm of the input in (line~1 of \autoref{alg:TT_SVD}).
Second, using the QR trick \eqref{eq:QR_trick}, we replace the large SVD by a QR decomposition followed by a smaller SVD of the triangular factor $R$ (line~5-6).
In addition, we can use the matrix of right singular vectors $V$ to calculate the work matrix $W^{(i-1)}$ for the next iteration (line~11 in both algorithms).
This has the benefit that we never need the orthogonal factor $Q$ of the QR decomposition which can be exploited in the implementation (see \autoref{sec:TSQR_algorithm}).
Third, we minimize data transfers that would be required for reshaping with appropriate padding of the data to avoid cache thrashing.
So in line~11, we directly store $W^{(i-1)}$ in the desired memory layout for the next iteration. This replaces the additional $\operatorname{reshape}$ operation in line~6 of \autoref{alg:TT_SVD}.
We assume further that the input tensor $X$ already has a suitable memory layout such that we do not need to copy the data for the first iteration (line~3).
The costly operations in this algorithm are the tall-skinny Q-less QR decomposition (line~5) and the tall-skinny matrix-matrix product fused with a $\operatorname{reshape}$ operation (line~11).

\paragraph{Memory layout}
The chosen memory layout has a significant effect on the performance.
A particular problem is cache thrashing (see e.g.\ \cite{Hager2010}).
An example for this is shown in \autoref{sec:numerical_experiments} in \autoref{fig:TSMM_vs_MKL}~(b).
This effect occurs for data accesses with strides that are multiples of $2^k, k\in\mathbf N$ with e.g.\ $k > 6$.
This easily happens in the TT-SVD algorithm if the individual dimensions $n_i$ are multiples of~$2$, for example when storing $W^{(i)}$ in a column-major layout.
To avoid this problem, one can use padding: that means filling-in a few zero entries such that the stride is not close to a multiple of $2^k$
(in our implementation padding is performed for all matrices such as $W^{(i)}$ to obtain strides of the form $2^6(2l+1), l\in\mathbf N$).
In addition, the required matrix operations in the TT-SVD algorithm are memory-bound in many cases (see \autoref{sec:performance_analysis} for a detailed discussion).
That means that data locality in these operations plays a crucial role.
On older multi-core CPUs, a row-major memory layout in operations with tall-skinny matrices operations is favorable, see e.g.\ the comparison in \cite{RoehrigZoellner2015}.
On newer CPUs (Intel Skylake and newer), there is no such performance penalty for using a column-major memory layout (observation of the authors).
Therefore, we employ a column-major memory layout for all matrices in \autoref{alg:TSQR_TT_SVD}.
And the leading dimensions of the input tensor $X$ are stored contiguously (Fortran ordering).
As indicated above, we thus assume that the stride of the last dimension includes appropriate padding.

\subsubsection{Algorithmic variants}\label{sec:TSQR_TT-SVD_variants}
In the following, we discuss some interesting algorithmic variations of the TSQR TT-SVD algorithm.
\paragraph{Thick-bounds variant}
If the dimensions $n_i$ in the first iterations of \autoref{alg:TSQR_TT_SVD} are small, the required tall-skinny matrix operations become strongly memory-bound.
We can increase the compute intensity by combining the right-most dimensions of the input tensor as shown in \autoref{alg:Thick-Bounds_TT_SVD}.
\begin{algorithm}
\caption{Thick-bounds TT-SVD}\label{alg:Thick-Bounds_TT_SVD}
\begin{algorithmic}[1]
\REQUIRE $X\in\mathbf R^{n_1 \times \dots \times n_d}$, min.\ dimension $m_\text{min}$, min.\ reduction factor $f_1^\text{min}$,\\
         estimated TT ranks $\tilde r_i$ \hfill \COMMENT{or simply use $\tilde r_i = r_\text{max}$}
\ENSURE TT decomposition $T^{(1)}, \dots, T^{(d)}$
\STATE Choose \#dimensions $k$ to combine:\\
       minimal $k\in\{1,\dots,d\}$ with $m\ge \max(m_\text{min}, f_1^\text{min} \tilde r_{k-1})$, $m:=\prod_{i={d-k+1}}^d n_i$
\STATE $W \gets \operatorname{reshape}\left(X, \begin{pmatrix} n_1 & \cdots & n_{d-k} & m \end{pmatrix}\right)$
\STATE $T^{(1)},\dots,T^{(d-k)},\bar T^{(d-k+1)} \gets \operatorname{TT-SVD}(W)$
\STATE Recover $T^{(d-k+1)}, \dots, T^{(d)}$ from the TT-SVD of $\bar T^{(d-k+1)}$
\end{algorithmic}
\end{algorithm}
This approach allows a more efficient of the compute resources:
We suggest a heuristic (line 1) based on estimated TT ranks, on a minimal combined boundary dimension ($m_\text{min}$) and on a minimal estimated reduction of the work array size in the first TT-SVD iteration ($f_1^\text{min}$).
One can choose $m_\text{min}$ such that the compute intensity of the first TSQR step is close to the machine intensity.
The reduction factor is discussed in \autoref{sec:TT-SVD_perf_analysis}. The TT cores corresponding to the combined dimensions can be cheaply calculated afterwards (line 4).

\paragraph{Two-sided variant}
The matrix operations in \autoref{alg:TSQR_TT_SVD} become more costly for increasing TT-ranks.
And usually, the ranks are smaller near the left and right boundaries of the tensor-train.
So we can alternatingly calculate TT cores on the left and on the right as depicted in \autoref{alg:Two-Sided_TSQR_TT_SVD}.
The core idea here is to reduce the size of the work array in each iteration with lower computational costs.
\begin{algorithm}
\caption{Two-sided TSQR TT-SVD}\label{alg:Two-Sided_TSQR_TT_SVD}
\begin{algorithmic}[1]
\REQUIRE $X\in\mathbf R^{n_1 \times \dots \times n_d}$
\ENSURE TT decomposition $T^{(1)}, \dots, T^{(d)}$
\STATE $W^{(d)} \gets \operatorname{reshape}\left(X, \begin{pmatrix} \frac{\bar n}{n_d} & n_d \end{pmatrix}\right)$ \hfill \COMMENT{with the total size $\bar n$}
\FOR{$i = d, 1, d-1, 2, d-2, 3, \dots$}
  \STATE Calculate $R$ from the QR decomposition: $QR = W^{(i)}$
  \STATE Calculate small SVD: $\bar U \Sigma V^T = R$
  \IF{ $i < d/2$ }
    \STATE Get the new rank $r_i$ from truncating the SVD. \hfill \COMMENT{left case}
    \STATE $T^{(i)} \gets \operatorname{reshape}\left((V_{:,1:r_i}), \begin{pmatrix} r_{i-1} & n_i & r_i\end{pmatrix}\right)$
    \STATE $\bar  W^{(i)} \gets W^{(i)} V_{:,1:r_i}$
    \STATE Reshape and transpose $\bar W^{(i)}$ to get $W^{(d-i)}$ for the next iteration.
  \ELSE
    \STATE Get the new rank $r_{i-1}$ from truncating the SVD. \hfill \COMMENT{right case}
    \STATE $T^{(i)} \gets \operatorname{reshape}\left((V_{:,1:r_{i-1}})^T, \begin{pmatrix} r_{i-1} & n_i & r_i\end{pmatrix}\right)$
    \STATE $\bar W^{(i)} \gets W^{(i)} V_{:,1:r_{i-1}}$
    \STATE Transpose and reshape $\bar W^{(i)}$ to get $W^{(d-i+1)}$ for the next iteration.
  \ENDIF
\ENDFOR
\IF{ $d$ is even }
  \STATE $T^{(d/2)} \gets \operatorname{reshape}\left(W^{(d/2)}, \begin{pmatrix} r_{d/2-1} & n_{d/2} & r_{d/2} \end{pmatrix}\right)$
\ELSE
  \STATE $T^{((d+1)/2)} \gets \operatorname{reshape}\left((W^{((d+1)/2)})^T, \begin{pmatrix} r_{(d-1)/2} & n_{(d+1)/2} & r_{(d+1)/2} \end{pmatrix}\right)$
\ENDIF
\end{algorithmic}
\end{algorithm}
The algorithm includes additional memory operations (line~9 and~14) to reorder data.
These can be avoided by directly calculating the QR decomposition of the transposed work matrix ($(W^{(i)})^T$).
However, our TSQR implementation requires a specific memory layout which makes the additional reordering necessary.
It is also difficult to fuse the reordering with the preceding tall-skinny matrix multiplication efficiently due to complex index calculations.
So \autoref{alg:Two-Sided_TSQR_TT_SVD} illustrates our slightly sub-optimal implementation.
As the ideas of the thick-bounds variant and the two-sided variant are independent from each other, we can combine them.
In our numerical experiments, we thus directly show timings for a two-sided algorithm with thick-bounds.

\paragraph{Distributed TSQR TT-SVD}
We can extend \autoref{alg:TSQR_TT_SVD} to the case where the input tensor is distributed onto multiple compute nodes.
For simplicity, we assume that the tensor is distributed along the first $k$ dimensions and that the number of processes matches the total size of those dimensions.
This is sketched in \autoref{alg:Distributed_TSQR_TT_SVD}.
\begin{algorithm}
\caption{Distributed TSQR TT-SVD}\label{alg:Distributed_TSQR_TT_SVD}
\begin{algorithmic}[1]
\REQUIRE $X\in\mathbf R^{n_1 \times \dots \times n_d}$ distributed along the first k dimensions $n_1\times\dots\times n_k$, $k\ll d$ onto m processes $j=1,\dots,m$ with $m=\prod_{i=1}^k n_i$
\ENSURE TT decomposition $T^{(1)}, \dots, T^{(d)}$ (duplicated on all processes)
\STATE Read local part: $W^{(j)} \gets X_{i_1^{(j)},\dots,i_k^{(j)},:,\dots,:}$  on process $j$
\STATE $ V^{(j)},T^{(k+1)}, \dots, T^{(d)} \gets \operatorname{TSQR-TT-SVD}(W)$ \hfill \COMMENT{\autoref{alg:TSQR_TT_SVD} with global QR}
\STATE Gather $V \gets \operatorname{reshape}\left(\begin{pmatrix}V^{(1)} & \cdots & V^{(m)}\end{pmatrix},\begin{pmatrix}n_1 & \dots & n_{k-1} & n_kr_k\end{pmatrix}\right)$
\STATE Recover $T^{(1)},\dots,T^{(k)}$ from the TT-SVD of $V$
\end{algorithmic}
\end{algorithm}
The only change required for the distributed case is that the tall-skinny QR decomposition in line~5 of \autoref{alg:TSQR_TT_SVD} needs to perform an additional global reduction of the (local) triangular factors (see discussion in \autoref{sec:TSQR_algorithm}).
All other costly operations of \autoref{alg:TSQR_TT_SVD} are completely independent on all processes.
The work for the small SVDs is duplicated on each process as well as the work for recovering the first few dimensions (line 4 of \autoref{alg:Distributed_TSQR_TT_SVD}).
Of course, the assumption that the data is distributed along the first dimensions is quite restrictive.
For other cases, we could first calculate the TT decomposition with reordered dimensions using this algorithm and in a post-processing step reorder the dimensions in the tensor-train by swapping dimensions through combining and splitting neighboring TT cores (still efficient if the TT representation is exponentially smaller than the input tensor).
We can locally use the thick-bounds variant in the distributed TSQR TT-SVD.
However, we cannot efficiently implement the two-sided variant in a distributed setting as the transpose operations would redistribute the data globally.

\subsection{Rank-preserving Q-less TSQR algorithm}\label{sec:TSQR_algorithm}
In this section, we present our highly efficient rank-preserving tall-skinny QR decomposition based on the communication-avoiding QR factorization in~\cite{Demmel2015}.
The QR decomposition is rank-preserving in the following sense:
It does not break down if the input matrix is rank-deficient, and the resulting triangular matrix $R$ has the same (numerical) rank as the input matrix.
For numerical robustness, we choose an implementation based on Householder reflections~\cite{Householder1958}.
As we do not need the matrix $Q$ in any form, its data is not even stored implicitly as in common LAPACK~\cite{LAPACK} routines to reduce the memory traffic.
The core building block transforms a rectangular matrix with zero lower left triangle to upper triangular form by an orthogonal transformation $Q$:
\begin{align}\label{eq:Qless_HouseholderQR}
&&\begin{pmatrix} M \\ R \end{pmatrix} &= Q \bar R, &&\text{with } M\in\mathbf R^{n_b\times m}, R, \bar R\in R^{m\times m}.
\end{align}
Here, $n_b$ denotes a block-size that is chosen as a multiple of the SIMD width such that the data fits into the CPU caches (e.g.\ $M$ and $R$ fit into L2, multiple Householder reflectors fit into L1 for internal blocking over columns).
This building block is similar to the LAPACK routine \verb|dtpqrt2| (for the special case that $M$ is rectangular).
Our implementation differs in the following three points:
First, \verb|dtpqrt2| overwrites the input matrix $M$ with the Householder reflection vectors.
We do not modify $M$ and store reflection vectors as long as they are needed in an internal buffer.
Second, we assume a special memory layout and alignment of $M$ and $R$; $R$ is overwritten by $\bar R$.
In contrast, LAPACK routines cope with inputs of arbitrary strides and alignment.
Third, our implementation is branch-less and uses fewer flops than the LAPACK reference implementation as discussed below.
Based on this building block, we implement a hybrid-parallel (MPI+OpenMP) TSQR scheme.
The TSQR algorithm is explained in detail in~\cite{Demmel2015}.
The main idea is that with the building block above, one can calculate triangular factors for blocks of the input matrix and combine them, e.g.\ for a flat tree reduction:
\begin{align*}
\begin{pmatrix} M_1 \\ M_2 \\ M_3 \end{pmatrix}
= \begin{pmatrix} Q_1 R_1 \\ M_2 \\ M_3 \end{pmatrix}
= \begin{pmatrix} \begin{pmatrix} & Q_1 \\I & \end{pmatrix} \begin{pmatrix}M_2\\R_1\end{pmatrix} \\ M_3 \end{pmatrix}
= \begin{pmatrix} Q_{12} R_{12} \\ M_3 \end{pmatrix}
= \cdots
= Q_{123} R_{123}.
\end{align*}
Each OpenMP thread performs a flat tree reduction (minimizing data transfers).
The resulting triangular $m\times m$ matrices are combined on the master thread (negligible overhead if the number of rows of the input matrix on each thread is large).
The results on multiple MPI processes are combined using an \verb|MPI_Allreduce| operation with a commutative MPI user reduction.
So the MPI library implementation decides about the actual reduction graph.

Some details of our main TSQR building block are illustrated in \autoref{alg:Qless_HouseholderQR}.
\begin{algorithm}[h!]
\caption{Householder~QR of a rectangular and a triangular matrix}\label{alg:Qless_HouseholderQR}
\begin{algorithmic}[1]
\REQUIRE $M\in\mathbf R^{n_b\times m}$, triangular $R\in\mathbf R^{m\times m}$, $\epsilon_{FP} > 0$\\
\COMMENT{$\epsilon_{FP}$ is the smallest positive normalized floating point number}
\ENSURE triangular $\bar R\in\mathbf R^{m\times m}$ that satisfies \eqref{eq:Qless_HouseholderQR}
\STATE $W_{1:n_b,:} \gets (M ; R)$
\FOR{$j=1, \dots, m$}
  \STATE $u \gets W_{j:n_b+j,j}$ \hfill\COMMENT{$w:=W_{j:n_b+j,j}$}
  \STATE $t \gets \|u\|_2^2 + \epsilon_{FP}, \quad \alpha \gets \sqrt{t + \epsilon_{FP}}$ \hfill\COMMENT{$\Rightarrow\alpha^2=\|w\|_2^2 + 2\epsilon_{FP}$}
  \STATE $\alpha \gets (-1)\cdot \alpha \textbf{ if } u_1 > 0 \textbf{ else } 1\cdot\alpha$ \hfill\COMMENT{implemented without branches}
  \STATE $t \gets t - \alpha u_1, \quad u_1 \gets u_1 - \alpha, \quad \beta \gets 1 / \sqrt{t}$\hfill\COMMENT{$\Rightarrow t=\|w\|_2^2+\epsilon_{FP}-w_1\alpha$}
  \STATE $v \gets \beta u$ \hfill\COMMENT{$\Rightarrow v=(w-\alpha e_1)/\sqrt{t}$}
  \STATE $W_{n_b+1:n_b+j,j} \gets (W_{1:j-1,j};\alpha)$
  \FOR{ $k=j+1, \dots, m$ }
    \STATE $\gamma \gets v^T W_{j:j+n_b,k}$
    \STATE $W_{j:j+n_b,k} \gets W_{j:j+n_b,k} - \gamma v$
  \ENDFOR
\ENDFOR
\STATE $\bar R \gets W_{n_b+1:n_b+m,:}$
\end{algorithmic}
\end{algorithm} 
There are two numerical differences with respect to the LAPACK reference implementation:
First, we calculate scaled Householder reflection vectors $v$ with $\|v\|_2=\sqrt{2}$ to avoid some additional multiplications.
Second, we add the term $\epsilon_{FP}$ in line~4 to prevent a break-down (division by zero in line~6).
In contrast, the reference implementation (\verb|dlarfg|) checks if $\|u\|_2$ is equal to zero or almost zero and performs different (expensive) steps depending on that.
So our implementation avoids a conditional branch at the cost of some numerical robustness.
We emphasize that through adding $\epsilon_{FP}$ twice as in line~4, we obtain in exact arithmetic:
\begin{align}
\|v\|_2^2 &= \frac{\|w - \alpha e_1\|_2^2}{t}
= \frac{\|w\|_2^2 - 2w_1\alpha + \alpha^2}{\|w\|_2^2 + \epsilon_{FP} - w_1\alpha}
=\frac{2\|w\|_2^2 - 2w_1\alpha + 2\epsilon_{FP}}{\|w\|_2^2 + \epsilon_{FP} - w_1\alpha} 
= 2
\end{align}
In inexact arithmetic, this also holds approximately as long as $2\|u\|_2^2+\epsilon_{FP}$ is in the range where the floating point arithmetic is accurate (no denormal numbers, e.g.\ $2\|u\|_2^2\lesssim10^{308}$ and $\epsilon_{FP}\approx10^{-308}$ for double precision).
So $I-vv^T$ is a valid Householder reflection even for $\|u\|_2\approx 0$.

The actual implementation looks more complicated as it uses a recursive blocking of columns:
On each recursion level, it splits the matrix into a left block and a right block and first processes the left block, then applies reflections to the right block and proceeds with the right block.
This is numerically equivalent to the algorithm shown here as it only reorders the loop iterations.
In addition, we avoid the copy in line~1 by just pointing to the actual data.
The conditional sign flip in line~5 is compiled to floating point instructions (masked blending).
Moreover, the vector operations in all iterations use vectors of the same length ($n_b+1$) which facilitates the SIMD parallelization.

\section{Performance analysis}\label{sec:performance_analysis}
In this section we first analyze the performance of the building blocks and then model the run-time of the complete TT-SVD algorithm.
We assume that the dense input tensor is stored in main memory.
If we read the input data from disc, the same principles apply but the gap between the bandwidth and the floating point performance is even larger.

\subsection{Building blocks}
The main building blocks in \autoref{alg:TSQR_TT_SVD} are tall-skinny QR decompositions and matrix-matrix multiplications that we discuss in the following.

\subsubsection{Q-less TSQR algorithm}\label{sec:TSQR_perf_analysis}
For $X\in\mathbf R^{n\times m}$ with $n\gg m$, the TSQR algorithm described in \autoref{sec:TSQR_algorithm} calculates the triangular matrix $R\in\mathbf R^{m\times m}$ of the $QR$ decomposition of $X$.
A cache-friendly implementation only reads $X$ once from main memory ($V_\text{read}=8nm$ bytes for double precision).
Thus, a pure load benchmark shows the upper bound for the possible bandwidth $b_s = b_\text{load}$.
We estimate the required floating point operations of the Householder QR reduction by considering lines 4, 7, 10 and 11 in \autoref{alg:Qless_HouseholderQR}.
We can simplify this to $\sum_{k=1}^m (m-k+1) = \frac{m(m+1)}{2}$ dot products and scaled vector additions (axpy) of length $n_b+1$.
This results in $m(m+1)(n_b+1)$ FMA (fused multiply-add) instructions, respectively $2m(m+1)(n_b+1)$ floating point operations.
We need to perform $n/n_b$ such reduction steps assuming a flat TSQR reduction scheme.
In practice, we perform some additional reduction steps with a different block size $n_b$ depending on the number of OpenMP threads and MPI processes, but these are negligible for large $n$.
Overall, we obtain:
\begin{align}
&&n_\text{flops} &\approx \frac{n}{n_b} \left(2m(m+1)(n_b+1)\right) \approx \left(1+\frac{1}{n_b}\right) 2nm^2,&&\\
&\Rightarrow& I_c &=\frac{n_\text{flops}}{V_\text{read}}\approx \left(1+\frac{1}{n_b}\right)\frac{m}{4}. 
\end{align}
The compute intensity shows that the algorithm is memory-bound for $m$ up to $\sim 50$ (assuming $n_b \gg 1$) on the considered hardware (see \autoref{tab:peak_bandwidth_and_flops}).

\subsubsection{Tall-skinny matrix-matrix multiplication (TSMM)}\label{sec:TSMM_perf_analysis}
For matrices $X\in\mathbf R^{n \times m}$, $M\in\mathbf R^{m \times k}$ and $Y\in\mathbf R^{\hat n \times \hat m}$ with $n\gg m$ and $\hat n\hat m = nk$, the fused kernel for a tall-skinny matrix-matrix multiplication and a $\operatorname{reshape}$ operation calculates:
\begin{align*}
  Y \gets \operatorname{reshape}\left(X M, \begin{pmatrix} \hat n & \hat m\end{pmatrix}\right).
\end{align*}
The $\operatorname{reshape}$ operation just modifies the memory layout of the result and has no influence on the performance.
The matrix-matrix multiplication requires $2nmk$ floating point operations and can exploit modern FMA (fused multiply-add) instructions.
The operation reads $X$ ($8nm$ bytes for double precision) and writes $Y$ ($8nk$ bytes) using non-temporal stores.
The ratio of read to write volume is defined by $m/k$.
In our experiments, we usually have $m/k \approx 2$, so we approximate the limiting bandwidth with a STREAM benchmark: $b_s = b_\text{STREAM}$.
The resulting double precision compute intensity is $I_c = \frac{mk}{4(m+k)} \approx \frac{m}{12}$ for $m/k \approx 2$.
So on the considered hardware, this operation is memory-bound for $m$ up to $\sim 150$ (see \autoref{tab:peak_bandwidth_and_flops}).

\subsection{Complete TT-SVD algorithm}\label{sec:TT-SVD_perf_analysis}
We only analyze the optimized TSQR TT-SVD algorithm depicted in \autoref{alg:TSQR_TT_SVD}.
The analysis includes the idea of the thick-bounds variant in order to adjust algorithmic parameters.

We first consider the case that the number of columns $m$ in the required building blocks is small enough such that they operate in the memory-bound regime (small $r_\text{max}$ and small $n_i$).
For this case, we can estimate a lower bound for the run-time by considering the data transfers in the main building blocks:
One TSQR TT-SVD iteration first reads the work matrix (TSQR) and then reads it again and writes a reduced work matrix (TSMM).
So for each iteration $j=1,\dots,d-1$, we obtain the data volume: $V_\text{read+write} = 2\bar n + f_j\bar n$.
Here, $\bar n$ denotes the total size of the input data of that iteration and $f_j \in (0,1]$ a reduction factor ($f_j = \frac{r_{i-1}}{n_i r_i}$ with $i=d-j+1$ in \autoref{alg:TSQR_TT_SVD}).
This is the lowest data transfer volume possible for one step in the TT-SVD algorithm if we assume that we need to consider all input data before we can compress it (\emph{global} truncated SVD or QR decomposition).
\emph{Local} transformations are possible by e.g.\ calculating truncated SVDs of blocks of the input matrix that fit into the cache and combining them later.
Such a \emph{local} approach could at best improve the performance by roughly a factor of two as it would only read the data once instead of twice.
However, this reduces the accuracy of the approximation (combining multiple \emph{local} approximations instead of one \emph{global} approximation for each step).
For the complete TSQR TT-SVD algorithm, we sum up the data transfers of all iterations:
\begin{align}\label{eq:TT-SVD_data_volume}
\bar V_\text{read+write} &= 2\bar n(1 + f_1 + f_1 f_2 + \dots) + \bar n (f_1 + f_1 f_2 + \dots) \lesssim \frac{2 \bar n}{1-\bar f} + \frac{\bar f \bar n}{1-\bar f},
\end{align}
with $1>\bar f \ge f_j$ and the total size of the input tensor $\bar n$.
To optimize data transfers, we thus need a significant reduction $f_1 \ll 1$ of the size of the work matrix in the first step.
This is exactly the idea of the thick-bounds variant discussed in \autoref{sec:TSQR_TT-SVD_variants}.
Overall, this indicates that small reduction factors $f_j$ would be beneficial.
However, by combining dimensions to reduce $f_j$ in the steps of the algorithm, the compute intensity increases and at some point the building blocks become compute-bound.
For a rank-1 approximation, we can choose a small reduction factor (e.g.\ $\bar f = 1/16$ in our implementation) and for larger maximal rank, we use the choice $\bar f=1/2$. This results in overall transfer volumes of:
\begin{align}
\bar V_\text{read+write} &\lesssim
\begin{cases}
  2.2\bar n  \quad\text{for } \bar f = \frac{1}{16}, \\
  5.0\bar n  \quad\text{for } \bar f = \frac{1}{2}.
\end{cases}.
\end{align}
So for strongly memory-bound cases (small $r_\text{max}$ and $n_i$), we expect a run-time in the order of the time required for copying the input tensor $X$ (in memory).

In contrast, for larger ranks, the problem becomes compute-bound.
The building blocks need approximately $2nm^2$ (TSQR), respectively $2nmk$ (TSMM) floating point operations for an input matrix of size $n\times m$, respectively the multiplication of an $n\times m$ with an $m\times k$ matrix.
In iteration $j$, we have dimensions $nm=\bar n \prod_{l=1}^{j-1} f_j$, $m=k/f_j$ and $k = r_{i-1}$.
So for the complete algorithm, we obtain with $f_j \approx \bar f$ and $f_j \le r_\text{max}$:
\begin{align}\label{eq:TT-SVD_operations}
n_\text{flops} \approx 2\bar n\left(r_{d-1}\frac{1+f_1}{f_1} + f_1r_{d-2}\frac{1+f_2}{f_2} + \dots \right) \lesssim 2\bar n r_\text{max}\left(\frac{1}{\bar f} + \frac{2}{1-\bar f}\right).
\end{align}
This shows that combining more dimensions to reduce $f_1$ in the first step increases the work.
The optimal reduction factor to minimize the number of operations is roughly $\bar f \approx 0.4$.
With the choices for $\bar f$ from above, we obtain:
\begin{align}
n_\text{flops} &\lesssim
\begin{cases}
  36 \bar n r_\text{max} \quad\text{for } \bar f = \frac{1}{16},\\
  12 \bar n r_\text{max} \quad\text{for } \bar f = \frac{1}{2}.
\end{cases}
\end{align}
This approximation neglects the operations of the small SVD calculations of the triangular factors.
So it is only valid for higher dimensions, e.g.\ for $\bar n := \prod n_i \gg (\max n_i)^3$.
For the compute-bound case, we expect roughly a linear increase in run-time for increasing values of $r_\text{max}$ given fixed dimensions and a fixed reduction factor $\bar f$ (this requires combining more dimensions).
For large dimensions $n_i$ the reduction factors become very small ($f_j \sim 1/n_i$ without splitting dimensions) and thus the computational complexity increases.
In our implementation (see \autoref{alg:Thick-Bounds_TT_SVD}), we only combine dimensions at the boundary, so we can only influence the first reduction factor $f_1$.

\section{Numerical experiments}\label{sec:numerical_experiments}
In this section, we first discuss the performance of the building blocks and then consider different variants and implementations of the complete TT-SVD algorithm.
We perform all measurements on a small CPU cluster, see \autoref{tab:peak_bandwidth_and_flops} for information on the hardware.
For most of the experiments, we only use a single CPU socket to avoid NUMA effects (accessing memory from another CPU socket).
We implemented all required algorithms in a templated C++ library~\cite{RoehrigZoellner2021pitts} based on MPI, OpenMP and CPU SIMD intrinsics.
The library includes scripts for all experiments.
Comparisons of building blocks with the Intel\textsuperscript{\textregistered} Math Kernel Library (MKL) are written in python using numpy.
We set up comparisons with other software very carefully:
In particular, we ran benchmarks multiple times and ignored the first runs to avoid measuring initialization overhead.
Furthermore, we checked that a high fraction of the computing time is spent in appropriate building blocks (like MKL functions) and not in some (python) layer above (using the linux tool \verb|perf|).
All calculations use double precision.
The input data in all experiments is uniformly random and we prescribe the dimensions, respectively tensor-train ranks.

\subsection{Building blocks}
The important building blocks are the Q-less TSQR algorithm and the tall-skinny matrix-matrix multiplication (fused with a $\operatorname{reshape}$ of the result).
Depending on the desired TT-rank in the TT-SVD algorithm, the number of columns $m$ changes for the tall-skinny matrices in the building blocks.
Therefore, we need to consider the performance for varying numbers of columns.

\subsubsection{Q-less TSQR algorithm}
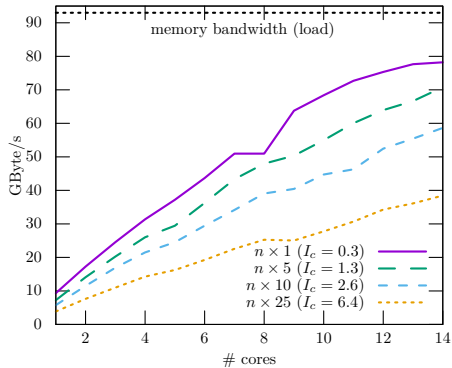
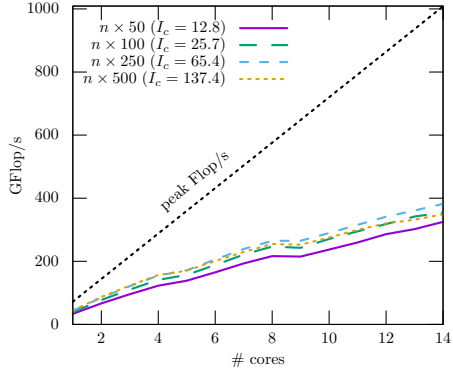
\begin{figure}
\begin{subfigure}[t]{0.47\textwidth}
\begin{tikzpicture}[gnuplot]
\tikzset{every node/.append style={scale=0.60}}
\path (0.000,0.000) rectangle (6.100,5.000);
\gpcolor{color=gp lt color border}
\gpsetlinetype{gp lt border}
\gpsetdashtype{gp dt solid}
\gpsetlinewidth{1.00}
\draw[gp path] (0.680,0.592)--(0.860,0.592);
\draw[gp path] (5.769,0.592)--(5.589,0.592);
\node[gp node right] at (0.570,0.592) {$0$};
\draw[gp path] (0.680,1.036)--(0.860,1.036);
\draw[gp path] (5.769,1.036)--(5.589,1.036);
\node[gp node right] at (0.570,1.036) {$10$};
\draw[gp path] (0.680,1.481)--(0.860,1.481);
\draw[gp path] (5.769,1.481)--(5.589,1.481);
\node[gp node right] at (0.570,1.481) {$20$};
\draw[gp path] (0.680,1.925)--(0.860,1.925);
\draw[gp path] (5.769,1.925)--(5.589,1.925);
\node[gp node right] at (0.570,1.925) {$30$};
\draw[gp path] (0.680,2.370)--(0.860,2.370);
\draw[gp path] (5.769,2.370)--(5.589,2.370);
\node[gp node right] at (0.570,2.370) {$40$};
\draw[gp path] (0.680,2.814)--(0.860,2.814);
\draw[gp path] (5.769,2.814)--(5.589,2.814);
\node[gp node right] at (0.570,2.814) {$50$};
\draw[gp path] (0.680,3.259)--(0.860,3.259);
\draw[gp path] (5.769,3.259)--(5.589,3.259);
\node[gp node right] at (0.570,3.259) {$60$};
\draw[gp path] (0.680,3.703)--(0.860,3.703);
\draw[gp path] (5.769,3.703)--(5.589,3.703);
\node[gp node right] at (0.570,3.703) {$70$};
\draw[gp path] (0.680,4.147)--(0.860,4.147);
\draw[gp path] (5.769,4.147)--(5.589,4.147);
\node[gp node right] at (0.570,4.147) {$80$};
\draw[gp path] (0.680,4.592)--(0.860,4.592);
\draw[gp path] (5.769,4.592)--(5.589,4.592);
\node[gp node right] at (0.570,4.592) {$90$};
\draw[gp path] (1.071,0.592)--(1.071,0.772);
\draw[gp path] (1.071,4.814)--(1.071,4.634);
\node[gp node center] at (1.071,0.407) {$2$};
\draw[gp path] (1.854,0.592)--(1.854,0.772);
\draw[gp path] (1.854,4.814)--(1.854,4.634);
\node[gp node center] at (1.854,0.407) {$4$};
\draw[gp path] (2.637,0.592)--(2.637,0.772);
\draw[gp path] (2.637,4.814)--(2.637,4.634);
\node[gp node center] at (2.637,0.407) {$6$};
\draw[gp path] (3.420,0.592)--(3.420,0.772);
\draw[gp path] (3.420,4.814)--(3.420,4.634);
\node[gp node center] at (3.420,0.407) {$8$};
\draw[gp path] (4.203,0.592)--(4.203,0.772);
\draw[gp path] (4.203,4.814)--(4.203,4.634);
\node[gp node center] at (4.203,0.407) {$10$};
\draw[gp path] (4.986,0.592)--(4.986,0.772);
\draw[gp path] (4.986,4.814)--(4.986,4.634);
\node[gp node center] at (4.986,0.407) {$12$};
\draw[gp path] (5.769,0.592)--(5.769,0.772);
\draw[gp path] (5.769,4.814)--(5.769,4.634);
\node[gp node center] at (5.769,0.407) {$14$};
\draw[gp path] (0.680,4.814)--(0.680,0.592)--(5.769,0.592)--(5.769,4.814)--cycle;
\node[gp node left] at (1.854,4.503) {memory bandwidth (load)};
\node[gp node center,rotate=-270] at (0.175,2.703) {GByte/s};
\node[gp node center] at (3.224,0.130) {\# cores};
\node[gp node right] at (4.819,1.559) {$n\times1$ ($I_c = 0.3$)};
\gpcolor{rgb color={0.580,0.000,0.827}}
\gpsetdashtype{gp dt 1}
\gpsetlinewidth{2.00}
\draw[gp path] (4.929,1.559)--(5.549,1.559);
\draw[gp path] (0.680,1.007)--(1.071,1.359)--(1.463,1.683)--(1.854,1.986)--(2.246,2.245)%
  --(2.637,2.533)--(3.029,2.857)--(3.420,2.857)--(3.812,3.426)--(4.203,3.631)--(4.595,3.823)%
  --(4.986,3.940)--(5.378,4.043)--(5.769,4.068);
\gpcolor{color=gp lt color border}
\node[gp node right] at (4.819,1.334) {$n\times5$ ($I_c = 1.3$)};
\gpcolor{rgb color={0.000,0.620,0.451}}
\gpsetdashtype{gp dt 2}
\draw[gp path] (4.929,1.334)--(5.549,1.334);
\draw[gp path] (0.680,0.914)--(1.071,1.218)--(1.463,1.492)--(1.854,1.748)--(2.246,1.901)%
  --(2.637,2.203)--(3.029,2.518)--(3.420,2.724)--(3.812,2.830)--(4.203,3.035)--(4.595,3.261)%
  --(4.986,3.434)--(5.378,3.554)--(5.769,3.738);
\gpcolor{color=gp lt color border}
\node[gp node right] at (4.819,1.109) {$n\times10$ ($I_c = 2.6$)};
\gpcolor{rgb color={0.337,0.706,0.914}}
\gpsetdashtype{gp dt 3}
\draw[gp path] (4.929,1.109)--(5.549,1.109);
\draw[gp path] (0.680,0.848)--(1.071,1.103)--(1.463,1.344)--(1.854,1.547)--(2.246,1.682)%
  --(2.637,1.901)--(3.029,2.110)--(3.420,2.328)--(3.812,2.389)--(4.203,2.580)--(4.595,2.649)%
  --(4.986,2.923)--(5.378,3.056)--(5.769,3.198);
\gpcolor{color=gp lt color border}
\node[gp node right] at (4.819,0.884) {$n\times25$ ($I_c = 6.4$)};
\gpcolor{rgb color={0.902,0.624,0.000}}
\gpsetdashtype{gp dt 4}
\draw[gp path] (4.929,0.884)--(5.549,0.884);
\draw[gp path] (0.680,0.764)--(1.071,0.931)--(1.463,1.080)--(1.854,1.226)--(2.246,1.314)%
  --(2.637,1.447)--(3.029,1.595)--(3.420,1.715)--(3.812,1.705)--(4.203,1.826)--(4.595,1.956)%
  --(4.986,2.115)--(5.378,2.197)--(5.769,2.301);
\gpcolor{color=gp lt color border}
\gpsetlinetype{gp lt axes}
\gpsetdashtype{gp dt axes}
\draw[gp path] (0.680,4.725)--(0.731,4.725)--(0.783,4.725)--(0.834,4.725)--(0.886,4.725)%
  --(0.937,4.725)--(0.988,4.725)--(1.040,4.725)--(1.091,4.725)--(1.143,4.725)--(1.194,4.725)%
  --(1.245,4.725)--(1.297,4.725)--(1.348,4.725)--(1.400,4.725)--(1.451,4.725)--(1.502,4.725)%
  --(1.554,4.725)--(1.605,4.725)--(1.657,4.725)--(1.708,4.725)--(1.759,4.725)--(1.811,4.725)%
  --(1.862,4.725)--(1.914,4.725)--(1.965,4.725)--(2.017,4.725)--(2.068,4.725)--(2.119,4.725)%
  --(2.171,4.725)--(2.222,4.725)--(2.274,4.725)--(2.325,4.725)--(2.376,4.725)--(2.428,4.725)%
  --(2.479,4.725)--(2.531,4.725)--(2.582,4.725)--(2.633,4.725)--(2.685,4.725)--(2.736,4.725)%
  --(2.788,4.725)--(2.839,4.725)--(2.890,4.725)--(2.942,4.725)--(2.993,4.725)--(3.045,4.725)%
  --(3.096,4.725)--(3.147,4.725)--(3.199,4.725)--(3.250,4.725)--(3.302,4.725)--(3.353,4.725)%
  --(3.404,4.725)--(3.456,4.725)--(3.507,4.725)--(3.559,4.725)--(3.610,4.725)--(3.661,4.725)%
  --(3.713,4.725)--(3.764,4.725)--(3.816,4.725)--(3.867,4.725)--(3.918,4.725)--(3.970,4.725)%
  --(4.021,4.725)--(4.073,4.725)--(4.124,4.725)--(4.175,4.725)--(4.227,4.725)--(4.278,4.725)%
  --(4.330,4.725)--(4.381,4.725)--(4.432,4.725)--(4.484,4.725)--(4.535,4.725)--(4.587,4.725)%
  --(4.638,4.725)--(4.690,4.725)--(4.741,4.725)--(4.792,4.725)--(4.844,4.725)--(4.895,4.725)%
  --(4.947,4.725)--(4.998,4.725)--(5.049,4.725)--(5.101,4.725)--(5.152,4.725)--(5.204,4.725)%
  --(5.255,4.725)--(5.306,4.725)--(5.358,4.725)--(5.409,4.725)--(5.461,4.725)--(5.512,4.725)%
  --(5.563,4.725)--(5.615,4.725)--(5.666,4.725)--(5.718,4.725)--(5.769,4.725);
\gpsetlinetype{gp lt border}
\gpsetdashtype{gp dt solid}
\gpsetlinewidth{1.00}
\draw[gp path] (0.680,4.814)--(0.680,0.592)--(5.769,0.592)--(5.769,4.814)--cycle;
\gpdefrectangularnode{gp plot 1}{\pgfpoint{0.680cm}{0.592cm}}{\pgfpoint{5.769cm}{4.814cm}}
\end{tikzpicture}
\subcaption{Memory-bound case ($I_c \lesssim 11$).}\label{fig:TSQR_bandwidth}
\end{subfigure}
\hfill
\begin{subfigure}[t]{0.47\textwidth}
\begin{tikzpicture}[gnuplot]
\tikzset{every node/.append style={scale=0.60}}
\path (0.000,0.000) rectangle (6.100,5.000);
\gpcolor{color=gp lt color border}
\gpsetlinetype{gp lt border}
\gpsetdashtype{gp dt solid}
\gpsetlinewidth{1.00}
\draw[gp path] (0.900,0.592)--(1.080,0.592);
\draw[gp path] (5.769,0.592)--(5.589,0.592);
\node[gp node right] at (0.790,0.592) {$0$};
\draw[gp path] (0.900,1.429)--(1.080,1.429);
\draw[gp path] (5.769,1.429)--(5.589,1.429);
\node[gp node right] at (0.790,1.429) {$200$};
\draw[gp path] (0.900,2.266)--(1.080,2.266);
\draw[gp path] (5.769,2.266)--(5.589,2.266);
\node[gp node right] at (0.790,2.266) {$400$};
\draw[gp path] (0.900,3.103)--(1.080,3.103);
\draw[gp path] (5.769,3.103)--(5.589,3.103);
\node[gp node right] at (0.790,3.103) {$600$};
\draw[gp path] (0.900,3.939)--(1.080,3.939);
\draw[gp path] (5.769,3.939)--(5.589,3.939);
\node[gp node right] at (0.790,3.939) {$800$};
\draw[gp path] (0.900,4.776)--(1.080,4.776);
\draw[gp path] (5.769,4.776)--(5.589,4.776);
\node[gp node right] at (0.790,4.776) {$1000$};
\draw[gp path] (1.275,0.592)--(1.275,0.772);
\draw[gp path] (1.275,4.814)--(1.275,4.634);
\node[gp node center] at (1.275,0.407) {$2$};
\draw[gp path] (2.024,0.592)--(2.024,0.772);
\draw[gp path] (2.024,4.814)--(2.024,4.634);
\node[gp node center] at (2.024,0.407) {$4$};
\draw[gp path] (2.773,0.592)--(2.773,0.772);
\draw[gp path] (2.773,4.814)--(2.773,4.634);
\node[gp node center] at (2.773,0.407) {$6$};
\draw[gp path] (3.522,0.592)--(3.522,0.772);
\draw[gp path] (3.522,4.814)--(3.522,4.634);
\node[gp node center] at (3.522,0.407) {$8$};
\draw[gp path] (4.271,0.592)--(4.271,0.772);
\draw[gp path] (4.271,4.814)--(4.271,4.634);
\node[gp node center] at (4.271,0.407) {$10$};
\draw[gp path] (5.020,0.592)--(5.020,0.772);
\draw[gp path] (5.020,4.814)--(5.020,4.634);
\node[gp node center] at (5.020,0.407) {$12$};
\draw[gp path] (5.769,0.592)--(5.769,0.772);
\draw[gp path] (5.769,4.814)--(5.769,4.634);
\node[gp node center] at (5.769,0.407) {$14$};
\draw[gp path] (0.900,4.814)--(0.900,0.592)--(5.769,0.592)--(5.769,4.814)--cycle;
\node[gp node left,rotate=40] at (2.024,2.098) {peak Flop/s};
\node[gp node center,rotate=-270] at (0.175,2.703) {GFlop/s};
\node[gp node center] at (3.334,0.130) {\# cores};
\node[gp node right] at (2.990,4.521) {$n\times50$ ($I_c = 12.8$)};
\gpcolor{rgb color={0.580,0.000,0.827}}
\gpsetdashtype{gp dt 1}
\gpsetlinewidth{2.00}
\draw[gp path] (3.100,4.521)--(3.720,4.521);
\draw[gp path] (0.900,0.733)--(1.275,0.871)--(1.649,0.993)--(2.024,1.106)--(2.398,1.172)%
  --(2.773,1.283)--(3.147,1.401)--(3.522,1.498)--(3.896,1.492)--(4.271,1.585)--(4.645,1.679)%
  --(5.020,1.789)--(5.394,1.856)--(5.769,1.952);
\gpcolor{color=gp lt color border}
\node[gp node right] at (2.990,4.296) {$n\times100$ ($I_c = 25.7$)};
\gpcolor{rgb color={0.000,0.620,0.451}}
\gpsetdashtype{gp dt 2}
\draw[gp path] (3.100,4.296)--(3.720,4.296);
\draw[gp path] (0.900,0.754)--(1.275,0.916)--(1.649,1.054)--(2.024,1.183)--(2.398,1.249)%
  --(2.773,1.380)--(3.147,1.518)--(3.522,1.625)--(3.896,1.608)--(4.271,1.725)--(4.645,1.826)%
  --(5.020,1.923)--(5.394,2.023)--(5.769,2.076);
\gpcolor{color=gp lt color border}
\node[gp node right] at (2.990,4.071) {$n\times250$ ($I_c = 65.4$)};
\gpcolor{rgb color={0.337,0.706,0.914}}
\gpsetdashtype{gp dt 3}
\draw[gp path] (3.100,4.071)--(3.720,4.071);
\draw[gp path] (0.900,0.772)--(1.275,0.947)--(1.649,1.100)--(2.024,1.240)--(2.398,1.309)%
  --(2.773,1.453)--(3.147,1.592)--(3.522,1.700)--(3.896,1.702)--(4.271,1.804)--(4.645,1.913)%
  --(5.020,2.020)--(5.394,2.103)--(5.769,2.191);
\gpcolor{color=gp lt color border}
\node[gp node right] at (2.990,3.846) {$n\times500$ ($I_c = 137.4$)};
\gpcolor{rgb color={0.902,0.624,0.000}}
\gpsetdashtype{gp dt 4}
\draw[gp path] (3.100,3.846)--(3.720,3.846);
\draw[gp path] (0.900,0.778)--(1.275,0.956)--(1.649,1.106)--(2.024,1.248)--(2.398,1.305)%
  --(2.773,1.433)--(3.147,1.553)--(3.522,1.659)--(3.896,1.649)--(4.271,1.745)--(4.645,1.847)%
  --(5.020,1.932)--(5.394,1.980)--(5.769,2.054);
\gpcolor{color=gp lt color border}
\gpsetlinetype{gp lt axes}
\gpsetdashtype{gp dt axes}
\draw[gp path] (0.900,0.893)--(0.949,0.933)--(0.998,0.972)--(1.048,1.012)--(1.097,1.052)%
  --(1.146,1.091)--(1.195,1.131)--(1.244,1.170)--(1.293,1.210)--(1.343,1.249)--(1.392,1.289)%
  --(1.441,1.328)--(1.490,1.368)--(1.539,1.408)--(1.589,1.447)--(1.638,1.487)--(1.687,1.526)%
  --(1.736,1.566)--(1.785,1.605)--(1.834,1.645)--(1.884,1.684)--(1.933,1.724)--(1.982,1.764)%
  --(2.031,1.803)--(2.080,1.843)--(2.130,1.882)--(2.179,1.922)--(2.228,1.961)--(2.277,2.001)%
  --(2.326,2.041)--(2.375,2.080)--(2.425,2.120)--(2.474,2.159)--(2.523,2.199)--(2.572,2.238)%
  --(2.621,2.278)--(2.671,2.317)--(2.720,2.357)--(2.769,2.397)--(2.818,2.436)--(2.867,2.476)%
  --(2.916,2.515)--(2.966,2.555)--(3.015,2.594)--(3.064,2.634)--(3.113,2.674)--(3.162,2.713)%
  --(3.212,2.753)--(3.261,2.792)--(3.310,2.832)--(3.359,2.871)--(3.408,2.911)--(3.457,2.950)%
  --(3.507,2.990)--(3.556,3.030)--(3.605,3.069)--(3.654,3.109)--(3.703,3.148)--(3.753,3.188)%
  --(3.802,3.227)--(3.851,3.267)--(3.900,3.306)--(3.949,3.346)--(3.998,3.386)--(4.048,3.425)%
  --(4.097,3.465)--(4.146,3.504)--(4.195,3.544)--(4.244,3.583)--(4.294,3.623)--(4.343,3.663)%
  --(4.392,3.702)--(4.441,3.742)--(4.490,3.781)--(4.539,3.821)--(4.589,3.860)--(4.638,3.900)%
  --(4.687,3.939)--(4.736,3.979)--(4.785,4.019)--(4.835,4.058)--(4.884,4.098)--(4.933,4.137)%
  --(4.982,4.177)--(5.031,4.216)--(5.080,4.256)--(5.130,4.296)--(5.179,4.335)--(5.228,4.375)%
  --(5.277,4.414)--(5.326,4.454)--(5.376,4.493)--(5.425,4.533)--(5.474,4.572)--(5.523,4.612)%
  --(5.572,4.652)--(5.621,4.691)--(5.671,4.731)--(5.720,4.770)--(5.769,4.810);
\gpsetlinetype{gp lt border}
\gpsetdashtype{gp dt solid}
\gpsetlinewidth{1.00}
\draw[gp path] (0.900,4.814)--(0.900,0.592)--(5.769,0.592)--(5.769,4.814)--cycle;
\gpdefrectangularnode{gp plot 1}{\pgfpoint{0.900cm}{0.592cm}}{\pgfpoint{5.769cm}{4.814cm}}
\end{tikzpicture}
\subcaption{Compute-bound case ($I_c \gtrsim 11$).}\label{fig:TSQR_flops}
\end{subfigure}
\caption{Single socket Q-less TSQR compared with the peak bandwidth respectively the peak Flop/s.
Based on \autoref{tab:peak_bandwidth_and_flops}, the machine intensity for this operation (pure load) is $1009/93\approx 11$~[Flops/Byte].
The input dimensions are chosen such that the matrix has a total size of $\sim3/14$~GByte per core.
The TSQR block size is $n_b=592$ for $m\lesssim160$ columns and then reduced linearly with $m$ (e.g.\ $n_b=192$ for $m=500$).}
\end{figure}

As analyzed in \autoref{sec:TSQR_perf_analysis}, the Q-less TSQR algorithm is memory-bound for $m$ up to $\sim50$ columns on the hardware used.
As we do not store the $Q$ matrix of the tall-skinny QR decomposition, its run-time is limited by the load memory bandwidth.
We expect a saturating behavior of the measured bandwidth up to the peak load bandwidth on 1-14 cores.
However, in \autoref{fig:TSQR_bandwidth} we see that the bandwidth is not fully saturated on 14~cores except for the case $n\times1$.
So our implementation partly seems to be limited by the in-core performance even for the memory-bound cases.
This effect increases with the number of columns, respectively with the compute intensity.
This indicates that our implementation is still sub-optimal.
In addition, the simple Roofline model based on the number of floating point operations is too optimistic for this case because the TSQR algorithm includes data dependencies as well as costly $\operatorname{sqrt}$ operations.
Overall we obtain more than 50\% of the peak bandwidth for small numbers of columns.
For the compute-bound case ($m\ge50$ on this hardware), we observe the expected linear scaling with the number of cores (see \autoref{fig:TSQR_flops}).
Our implementation achieves roughly 35\% of the peak performance here, independent of the number of columns.

\begin{figure}
\begin{subfigure}[t]{0.47\textwidth}
\begin{tikzpicture}[gnuplot]
\tikzset{every node/.append style={scale=0.60}}
\path (0.000,0.000) rectangle (6.100,5.000);
\gpcolor{color=gp lt color border}
\gpsetlinetype{gp lt border}
\gpsetdashtype{gp dt solid}
\gpsetlinewidth{1.00}
\draw[gp path] (0.790,0.592)--(0.970,0.592);
\draw[gp path] (5.769,0.592)--(5.589,0.592);
\node[gp node right] at (0.680,0.592) {$0$};
\draw[gp path] (0.790,1.436)--(0.970,1.436);
\draw[gp path] (5.769,1.436)--(5.589,1.436);
\node[gp node right] at (0.680,1.436) {$20$};
\draw[gp path] (0.790,2.281)--(0.970,2.281);
\draw[gp path] (5.769,2.281)--(5.589,2.281);
\node[gp node right] at (0.680,2.281) {$40$};
\draw[gp path] (0.790,3.125)--(0.970,3.125);
\draw[gp path] (5.769,3.125)--(5.589,3.125);
\node[gp node right] at (0.680,3.125) {$60$};
\draw[gp path] (0.790,3.970)--(0.970,3.970);
\draw[gp path] (5.769,3.970)--(5.589,3.970);
\node[gp node right] at (0.680,3.970) {$80$};
\draw[gp path] (0.790,4.814)--(0.970,4.814);
\draw[gp path] (5.769,4.814)--(5.589,4.814);
\node[gp node right] at (0.680,4.814) {$100$};
\draw[gp path] (0.790,0.592)--(0.790,0.772);
\draw[gp path] (0.790,4.814)--(0.790,4.634);
\node[gp node center] at (0.790,0.407) {$1$};
\draw[gp path] (1.345,0.592)--(1.345,0.682);
\draw[gp path] (1.345,4.814)--(1.345,4.724);
\draw[gp path] (1.670,0.592)--(1.670,0.682);
\draw[gp path] (1.670,4.814)--(1.670,4.724);
\draw[gp path] (1.901,0.592)--(1.901,0.682);
\draw[gp path] (1.901,4.814)--(1.901,4.724);
\draw[gp path] (2.079,0.592)--(2.079,0.682);
\draw[gp path] (2.079,4.814)--(2.079,4.724);
\draw[gp path] (2.226,0.592)--(2.226,0.682);
\draw[gp path] (2.226,4.814)--(2.226,4.724);
\draw[gp path] (2.349,0.592)--(2.349,0.682);
\draw[gp path] (2.349,4.814)--(2.349,4.724);
\draw[gp path] (2.456,0.592)--(2.456,0.682);
\draw[gp path] (2.456,4.814)--(2.456,4.724);
\draw[gp path] (2.550,0.592)--(2.550,0.682);
\draw[gp path] (2.550,4.814)--(2.550,4.724);
\draw[gp path] (2.635,0.592)--(2.635,0.772);
\draw[gp path] (2.635,4.814)--(2.635,4.634);
\node[gp node center] at (2.635,0.407) {$10$};
\draw[gp path] (3.190,0.592)--(3.190,0.682);
\draw[gp path] (3.190,4.814)--(3.190,4.724);
\draw[gp path] (3.515,0.592)--(3.515,0.682);
\draw[gp path] (3.515,4.814)--(3.515,4.724);
\draw[gp path] (3.745,0.592)--(3.745,0.682);
\draw[gp path] (3.745,4.814)--(3.745,4.724);
\draw[gp path] (3.924,0.592)--(3.924,0.682);
\draw[gp path] (3.924,4.814)--(3.924,4.724);
\draw[gp path] (4.070,0.592)--(4.070,0.682);
\draw[gp path] (4.070,4.814)--(4.070,4.724);
\draw[gp path] (4.194,0.592)--(4.194,0.682);
\draw[gp path] (4.194,4.814)--(4.194,4.724);
\draw[gp path] (4.301,0.592)--(4.301,0.682);
\draw[gp path] (4.301,4.814)--(4.301,4.724);
\draw[gp path] (4.395,0.592)--(4.395,0.682);
\draw[gp path] (4.395,4.814)--(4.395,4.724);
\draw[gp path] (4.480,0.592)--(4.480,0.772);
\draw[gp path] (4.480,4.814)--(4.480,4.634);
\node[gp node center] at (4.480,0.407) {$100$};
\draw[gp path] (5.035,0.592)--(5.035,0.682);
\draw[gp path] (5.035,4.814)--(5.035,4.724);
\draw[gp path] (5.360,0.592)--(5.360,0.682);
\draw[gp path] (5.360,4.814)--(5.360,4.724);
\draw[gp path] (5.590,0.592)--(5.590,0.682);
\draw[gp path] (5.590,4.814)--(5.590,4.724);
\draw[gp path] (5.769,0.592)--(5.769,0.682);
\draw[gp path] (5.769,4.814)--(5.769,4.724);
\draw[gp path] (0.790,4.814)--(0.790,0.592)--(5.769,0.592)--(5.769,4.814)--cycle;
\node[gp node left] at (1.670,4.350) {roofline model limit};
\node[gp node center,rotate=-270] at (0.175,2.703) {GByte/s};
\node[gp node center] at (3.279,0.130) {columns};
\gpcolor{rgb color={0.580,0.000,0.827}}
\gpsetdashtype{gp dt 1}
\gpsetlinewidth{2.00}
\draw[gp path] (0.790,4.211)--(1.345,3.925)--(1.670,3.745)--(1.901,3.657)--(2.079,3.559)%
  --(2.226,3.353)--(2.349,3.328)--(2.456,3.204)--(2.550,3.144)--(2.635,3.074)--(2.711,2.861)%
  --(2.781,2.867)--(2.845,2.744)--(2.904,2.688)--(2.960,2.601)--(3.011,2.537)--(3.060,2.417)%
  --(3.106,2.429)--(3.149,2.371)--(3.190,2.325)--(3.229,2.263)--(3.266,2.221)--(3.302,2.180)%
  --(3.336,2.283)--(3.369,2.217)--(3.400,2.187)--(3.431,2.157)--(3.460,2.106)--(3.488,2.075)%
  --(3.515,2.044)--(3.541,2.002)--(3.567,1.975)--(3.591,1.936)--(3.615,1.910)--(3.638,1.860)%
  --(3.661,1.867)--(3.683,1.869)--(3.704,1.810)--(3.725,1.851)--(3.745,1.815)--(3.765,1.750)%
  --(3.785,1.730)--(3.803,1.723)--(3.822,1.734)--(3.840,1.715)--(3.857,1.653)--(3.875,1.624)%
  --(3.892,1.649)--(3.908,1.650)--(3.924,1.629)--(3.940,1.621)--(3.956,1.616)--(3.971,1.596)%
  --(3.986,1.610)--(4.001,1.586)--(4.015,1.577)--(4.029,1.535)--(4.043,1.523)--(4.057,1.533)%
  --(4.070,1.513)--(4.084,1.506)--(4.097,1.515)--(4.109,1.486)--(4.122,1.470)--(4.134,1.455)%
  --(4.147,1.449)--(4.159,1.444)--(4.171,1.431)--(4.182,1.428)--(4.194,1.407)--(4.205,1.389)%
  --(4.216,1.411)--(4.227,1.390)--(4.238,1.378)--(4.249,1.373)--(4.260,1.367)--(4.270,1.361)%
  --(4.280,1.348)--(4.291,1.333)--(4.301,1.336)--(4.311,1.322)--(4.321,1.312)--(4.330,1.303)%
  --(4.340,1.314)--(4.349,1.291)--(4.359,1.289)--(4.368,1.292)--(4.377,1.270)--(4.386,1.273)%
  --(4.395,1.261)--(4.404,1.252)--(4.413,1.255)--(4.421,1.251)--(4.430,1.238)--(4.438,1.226)%
  --(4.447,1.231)--(4.455,1.222)--(4.463,1.213)--(4.472,1.209)--(4.480,1.207);
\draw[gp path] (4.480,1.198)--(4.804,1.005)--(5.035,0.914)--(5.214,0.844)--(5.360,0.801)%
  --(5.483,0.767)--(5.590,0.738)--(5.685,0.719)--(5.769,0.702);
\gpcolor{color=gp lt color border}
\gpsetlinetype{gp lt axes}
\gpsetdashtype{gp dt axes}
\draw[gp path] (0.790,4.518)--(0.840,4.518)--(0.891,4.518)--(0.941,4.518)--(0.991,4.518)%
  --(1.041,4.518)--(1.092,4.518)--(1.142,4.518)--(1.192,4.518)--(1.243,4.518)--(1.293,4.518)%
  --(1.343,4.518)--(1.394,4.518)--(1.444,4.518)--(1.494,4.518)--(1.544,4.518)--(1.595,4.518)%
  --(1.645,4.518)--(1.695,4.518)--(1.746,4.518)--(1.796,4.518)--(1.846,4.518)--(1.896,4.518)%
  --(1.947,4.518)--(1.997,4.518)--(2.047,4.518)--(2.098,4.518)--(2.148,4.518)--(2.198,4.518)%
  --(2.248,4.518)--(2.299,4.518)--(2.349,4.518)--(2.399,4.518)--(2.450,4.518)--(2.500,4.518)%
  --(2.550,4.518)--(2.601,4.518)--(2.651,4.518)--(2.701,4.518)--(2.751,4.518)--(2.802,4.518)%
  --(2.852,4.518)--(2.902,4.518)--(2.953,4.518)--(3.003,4.518)--(3.053,4.518)--(3.103,4.518)%
  --(3.154,4.518)--(3.204,4.518)--(3.254,4.518)--(3.305,4.518)--(3.355,4.518)--(3.405,4.518)%
  --(3.456,4.518)--(3.506,4.518)--(3.556,4.518)--(3.606,4.518)--(3.657,4.518)--(3.707,4.518)%
  --(3.757,4.518)--(3.808,4.518)--(3.858,4.294)--(3.908,4.069)--(3.958,3.857)--(4.009,3.659)%
  --(4.059,3.472)--(4.109,3.297)--(4.160,3.132)--(4.210,2.978)--(4.260,2.833)--(4.311,2.696)%
  --(4.361,2.568)--(4.411,2.448)--(4.461,2.335)--(4.512,2.229)--(4.562,2.129)--(4.612,2.036)%
  --(4.663,1.948)--(4.713,1.866)--(4.763,1.788)--(4.813,1.715)--(4.864,1.647)--(4.914,1.583)%
  --(4.964,1.522)--(5.015,1.466)--(5.065,1.413)--(5.115,1.363)--(5.165,1.316)--(5.216,1.272)%
  --(5.266,1.230)--(5.316,1.192)--(5.367,1.155)--(5.417,1.121)--(5.467,1.089)--(5.518,1.058)%
  --(5.568,1.030)--(5.618,1.003)--(5.668,0.978)--(5.719,0.955)--(5.769,0.933);
\gpsetlinetype{gp lt border}
\gpsetdashtype{gp dt solid}
\gpsetlinewidth{1.00}
\draw[gp path] (0.790,4.814)--(0.790,0.592)--(5.769,0.592)--(5.769,4.814)--cycle;
\gpdefrectangularnode{gp plot 1}{\pgfpoint{0.790cm}{0.592cm}}{\pgfpoint{5.769cm}{4.814cm}}
\end{tikzpicture}
\subcaption{Obtained memory bandwidth of our Q-less TSQR implementation compared to the roofline limit: $b_\text{roofline} =\min(P_\text{peak}/I_c, b_s)$.}\label{fig:TSQR_Roofline_bandwidth}
\end{subfigure}
\hfill
\begin{subfigure}[t]{0.47\textwidth}
\input{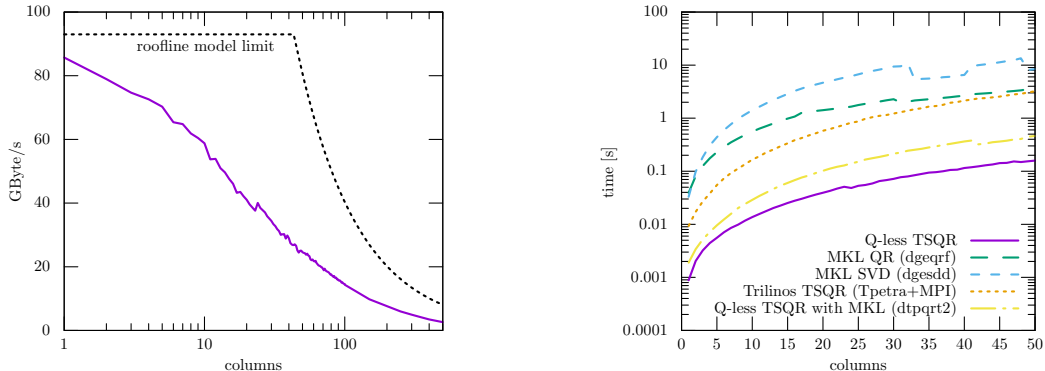}
\subcaption{Run-time for the QR decomposition (respectively a SVD) of a double precision $10^7 \times m$ matrix for $m=1, \dots, 50$ with our Q-less TSQR implementation, Intel MKL 2020.3, and Trilinos 13.0.0.}\label{fig:TSQR_vs_MKL}
\end{subfigure}
\caption{Single socket performance of tall-skinny matrix decompositions for varying numbers of columns.}
\end{figure}

\autoref{fig:TSQR_Roofline_bandwidth} shows the obtained bandwidth on a full socket and the Roofline limit depending on the number of columns $m$.
The kink in the Roofline limit denotes the point where the operation (theoretically) becomes compute-bound.
We see that the obtained bandwidth of our implementation decreases with the number of columns even in the memory-bound regime.
However, our specialized TSQR implementation is still significantly faster than just calling some standard QR algorithm that is not optimized for tall-skinny matrices.
This is illustrated by \autoref{fig:TSQR_vs_MKL}.
The comparison with MKL~QR is fair concerning the hardware setting (single socket with 14 cores, no NUMA effects).
However, it is unfair from the algorithmic point of view because we can directly discard $Q$ and exploit the known memory layout whereas the MKL~QR algorithm must work for all matrix shapes and any given memory layout and strides.
We also show the run-time of the MKL~SVD calculation for the same matrix dimensions.
Calculating the singular values and the right singular vectors from the resulting $R$ of the TSQR algorithm requires no significant additional time (SVD of $m\times m$ matrix with small $m$).
In addition, we measured the run-time of the Trilinos~\cite{Trilinos} TSQR algorithm with the Trilinos Tpetra library on one MPI process per core.
The Trilinos TSQR algorithm explicitly calculates the matrix $Q$ and it does not seem to exploit SIMD parallelism:
We only obtained scalar fused-multiply add instructions (FMA) instead of AVX512 (GCC 10.2 compiler with appropriate flags).
Due to these two reasons, the Trilinos TSQR is still slower than our almost optimal Q-less TSQR implementation by more than a factor of 10.
Finally, we replaced our implementation for reducing a triangular and a rectangular factor to triangular form in our Q-less TSQR implementation by the according low-level MKL routine \verb|dtpqrt2|.
In this case, we need to copy a block of the input matrix to a small buffer to avoid overwriting the input matrix (overhead of less than $\sim10$\% of the total time).
This variant achieves about $1/3$ of the performance of our specialized branch-less Householder~QR implementation.
Overall, the QR trick with our Q-less TSQR implementation reduces the run-time of the SVD calculation by roughly a factor of~50 compared to just calling standard LAPACK (MKL) routines.

\subsubsection{Tall-skinny matrix-matrix multiplication (TSMM)}
As analyzed in \autoref{sec:TSMM_perf_analysis}, the fused tall-skinny matrix-matrix multiplication and reshape is also memory-bound for $m$ up to $\sim150$ columns on the given hardware.
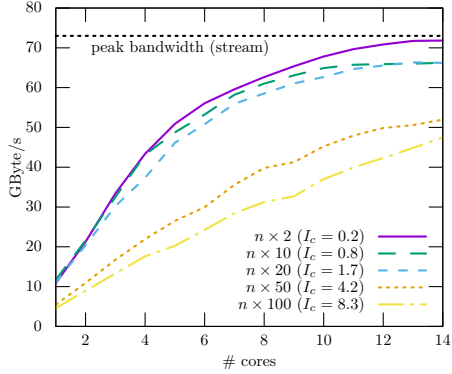
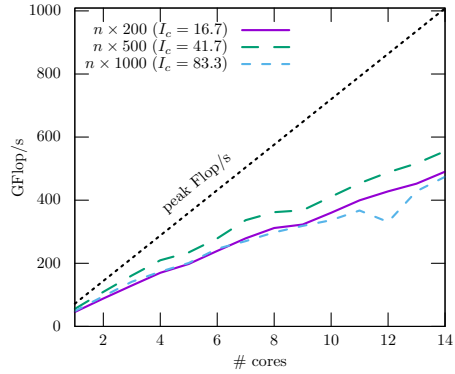
\begin{figure}
\begin{subfigure}[t]{0.47\textwidth}
\begin{tikzpicture}[gnuplot]
\tikzset{every node/.append style={scale=0.60}}
\path (0.000,0.000) rectangle (6.100,5.000);
\gpcolor{color=gp lt color border}
\gpsetlinetype{gp lt border}
\gpsetdashtype{gp dt solid}
\gpsetlinewidth{1.00}
\draw[gp path] (0.680,0.592)--(0.860,0.592);
\draw[gp path] (5.769,0.592)--(5.589,0.592);
\node[gp node right] at (0.570,0.592) {$0$};
\draw[gp path] (0.680,1.120)--(0.860,1.120);
\draw[gp path] (5.769,1.120)--(5.589,1.120);
\node[gp node right] at (0.570,1.120) {$10$};
\draw[gp path] (0.680,1.648)--(0.860,1.648);
\draw[gp path] (5.769,1.648)--(5.589,1.648);
\node[gp node right] at (0.570,1.648) {$20$};
\draw[gp path] (0.680,2.175)--(0.860,2.175);
\draw[gp path] (5.769,2.175)--(5.589,2.175);
\node[gp node right] at (0.570,2.175) {$30$};
\draw[gp path] (0.680,2.703)--(0.860,2.703);
\draw[gp path] (5.769,2.703)--(5.589,2.703);
\node[gp node right] at (0.570,2.703) {$40$};
\draw[gp path] (0.680,3.231)--(0.860,3.231);
\draw[gp path] (5.769,3.231)--(5.589,3.231);
\node[gp node right] at (0.570,3.231) {$50$};
\draw[gp path] (0.680,3.759)--(0.860,3.759);
\draw[gp path] (5.769,3.759)--(5.589,3.759);
\node[gp node right] at (0.570,3.759) {$60$};
\draw[gp path] (0.680,4.286)--(0.860,4.286);
\draw[gp path] (5.769,4.286)--(5.589,4.286);
\node[gp node right] at (0.570,4.286) {$70$};
\draw[gp path] (0.680,4.814)--(0.860,4.814);
\draw[gp path] (5.769,4.814)--(5.589,4.814);
\node[gp node right] at (0.570,4.814) {$80$};
\draw[gp path] (1.071,0.592)--(1.071,0.772);
\draw[gp path] (1.071,4.814)--(1.071,4.634);
\node[gp node center] at (1.071,0.407) {$2$};
\draw[gp path] (1.854,0.592)--(1.854,0.772);
\draw[gp path] (1.854,4.814)--(1.854,4.634);
\node[gp node center] at (1.854,0.407) {$4$};
\draw[gp path] (2.637,0.592)--(2.637,0.772);
\draw[gp path] (2.637,4.814)--(2.637,4.634);
\node[gp node center] at (2.637,0.407) {$6$};
\draw[gp path] (3.420,0.592)--(3.420,0.772);
\draw[gp path] (3.420,4.814)--(3.420,4.634);
\node[gp node center] at (3.420,0.407) {$8$};
\draw[gp path] (4.203,0.592)--(4.203,0.772);
\draw[gp path] (4.203,4.814)--(4.203,4.634);
\node[gp node center] at (4.203,0.407) {$10$};
\draw[gp path] (4.986,0.592)--(4.986,0.772);
\draw[gp path] (4.986,4.814)--(4.986,4.634);
\node[gp node center] at (4.986,0.407) {$12$};
\draw[gp path] (5.769,0.592)--(5.769,0.772);
\draw[gp path] (5.769,4.814)--(5.769,4.634);
\node[gp node center] at (5.769,0.407) {$14$};
\draw[gp path] (0.680,4.814)--(0.680,0.592)--(5.769,0.592)--(5.769,4.814)--cycle;
\node[gp node left] at (1.071,4.286) {peak bandwidth (stream)};
\node[gp node center,rotate=-270] at (0.175,2.703) {GByte/s};
\node[gp node center] at (3.224,0.130) {\# cores};
\node[gp node right] at (4.819,1.784) {$n\times2$ ($I_c = 0.2$)};
\gpcolor{rgb color={0.580,0.000,0.827}}
\gpsetdashtype{gp dt 1}
\gpsetlinewidth{2.00}
\draw[gp path] (4.929,1.784)--(5.549,1.784);
\draw[gp path] (0.680,1.167)--(1.071,1.714)--(1.463,2.355)--(1.854,2.881)--(2.246,3.278)%
  --(2.637,3.551)--(3.029,3.736)--(3.420,3.898)--(3.812,4.042)--(4.203,4.171)--(4.595,4.269)%
  --(4.986,4.331)--(5.378,4.378)--(5.769,4.383);
\gpcolor{color=gp lt color border}
\node[gp node right] at (4.819,1.559) {$n\times10$ ($I_c = 0.8$)};
\gpcolor{rgb color={0.000,0.620,0.451}}
\gpsetdashtype{gp dt 2}
\draw[gp path] (4.929,1.559)--(5.549,1.559);
\draw[gp path] (0.680,1.210)--(1.071,1.732)--(1.463,2.308)--(1.854,2.870)--(2.246,3.165)%
  --(2.637,3.402)--(3.029,3.664)--(3.420,3.812)--(3.812,3.921)--(4.203,4.018)--(4.595,4.062)%
  --(4.986,4.070)--(5.378,4.077)--(5.769,4.087);
\gpcolor{color=gp lt color border}
\node[gp node right] at (4.819,1.334) {$n\times20$ ($I_c = 1.7$)};
\gpcolor{rgb color={0.337,0.706,0.914}}
\gpsetdashtype{gp dt 3}
\draw[gp path] (4.929,1.334)--(5.549,1.334);
\draw[gp path] (0.680,1.178)--(1.071,1.667)--(1.463,2.187)--(1.854,2.564)--(2.246,3.027)%
  --(2.637,3.274)--(3.029,3.544)--(3.420,3.681)--(3.812,3.814)--(4.203,3.900)--(4.595,4.005)%
  --(4.986,4.052)--(5.378,4.093)--(5.769,4.087);
\gpcolor{color=gp lt color border}
\node[gp node right] at (4.819,1.109) {$n\times50$ ($I_c = 4.2$)};
\gpcolor{rgb color={0.902,0.624,0.000}}
\gpsetdashtype{gp dt 4}
\draw[gp path] (4.929,1.109)--(5.549,1.109);
\draw[gp path] (0.680,0.881)--(1.071,1.169)--(1.463,1.475)--(1.854,1.749)--(2.246,1.991)%
  --(2.637,2.173)--(3.029,2.471)--(3.420,2.693)--(3.812,2.770)--(4.203,2.979)--(4.595,3.123)%
  --(4.986,3.223)--(5.378,3.262)--(5.769,3.334);
\gpcolor{color=gp lt color border}
\node[gp node right] at (4.819,0.884) {$n\times100$ ($I_c = 8.3$)};
\gpcolor{rgb color={0.941,0.894,0.259}}
\gpsetdashtype{gp dt 5}
\draw[gp path] (4.929,0.884)--(5.549,0.884);
\draw[gp path] (0.680,0.834)--(1.071,1.066)--(1.463,1.294)--(1.854,1.521)--(2.246,1.661)%
  --(2.637,1.873)--(3.029,2.092)--(3.420,2.241)--(3.812,2.318)--(4.203,2.545)--(4.595,2.699)%
  --(4.986,2.826)--(5.378,2.958)--(5.769,3.099);
\gpcolor{color=gp lt color border}
\gpsetlinetype{gp lt axes}
\gpsetdashtype{gp dt axes}
\draw[gp path] (0.680,4.445)--(0.731,4.445)--(0.783,4.445)--(0.834,4.445)--(0.886,4.445)%
  --(0.937,4.445)--(0.988,4.445)--(1.040,4.445)--(1.091,4.445)--(1.143,4.445)--(1.194,4.445)%
  --(1.245,4.445)--(1.297,4.445)--(1.348,4.445)--(1.400,4.445)--(1.451,4.445)--(1.502,4.445)%
  --(1.554,4.445)--(1.605,4.445)--(1.657,4.445)--(1.708,4.445)--(1.759,4.445)--(1.811,4.445)%
  --(1.862,4.445)--(1.914,4.445)--(1.965,4.445)--(2.017,4.445)--(2.068,4.445)--(2.119,4.445)%
  --(2.171,4.445)--(2.222,4.445)--(2.274,4.445)--(2.325,4.445)--(2.376,4.445)--(2.428,4.445)%
  --(2.479,4.445)--(2.531,4.445)--(2.582,4.445)--(2.633,4.445)--(2.685,4.445)--(2.736,4.445)%
  --(2.788,4.445)--(2.839,4.445)--(2.890,4.445)--(2.942,4.445)--(2.993,4.445)--(3.045,4.445)%
  --(3.096,4.445)--(3.147,4.445)--(3.199,4.445)--(3.250,4.445)--(3.302,4.445)--(3.353,4.445)%
  --(3.404,4.445)--(3.456,4.445)--(3.507,4.445)--(3.559,4.445)--(3.610,4.445)--(3.661,4.445)%
  --(3.713,4.445)--(3.764,4.445)--(3.816,4.445)--(3.867,4.445)--(3.918,4.445)--(3.970,4.445)%
  --(4.021,4.445)--(4.073,4.445)--(4.124,4.445)--(4.175,4.445)--(4.227,4.445)--(4.278,4.445)%
  --(4.330,4.445)--(4.381,4.445)--(4.432,4.445)--(4.484,4.445)--(4.535,4.445)--(4.587,4.445)%
  --(4.638,4.445)--(4.690,4.445)--(4.741,4.445)--(4.792,4.445)--(4.844,4.445)--(4.895,4.445)%
  --(4.947,4.445)--(4.998,4.445)--(5.049,4.445)--(5.101,4.445)--(5.152,4.445)--(5.204,4.445)%
  --(5.255,4.445)--(5.306,4.445)--(5.358,4.445)--(5.409,4.445)--(5.461,4.445)--(5.512,4.445)%
  --(5.563,4.445)--(5.615,4.445)--(5.666,4.445)--(5.718,4.445)--(5.769,4.445);
\gpsetlinetype{gp lt border}
\gpsetdashtype{gp dt solid}
\gpsetlinewidth{1.00}
\draw[gp path] (0.680,4.814)--(0.680,0.592)--(5.769,0.592)--(5.769,4.814)--cycle;
\gpdefrectangularnode{gp plot 1}{\pgfpoint{0.680cm}{0.592cm}}{\pgfpoint{5.769cm}{4.814cm}}
\end{tikzpicture}
\subcaption{Memory-bound case ($I_c \lesssim 14$) measured with $n=10^7$.}\label{fig:TSMM_bandwidth}
\end{subfigure}
\hfill
\begin{subfigure}[t]{0.47\textwidth}
\begin{tikzpicture}[gnuplot]
\tikzset{every node/.append style={scale=0.60}}
\path (0.000,0.000) rectangle (6.100,5.000);
\gpcolor{color=gp lt color border}
\gpsetlinetype{gp lt border}
\gpsetdashtype{gp dt solid}
\gpsetlinewidth{1.00}
\draw[gp path] (0.900,0.592)--(1.080,0.592);
\draw[gp path] (5.769,0.592)--(5.589,0.592);
\node[gp node right] at (0.790,0.592) {$0$};
\draw[gp path] (0.900,1.429)--(1.080,1.429);
\draw[gp path] (5.769,1.429)--(5.589,1.429);
\node[gp node right] at (0.790,1.429) {$200$};
\draw[gp path] (0.900,2.266)--(1.080,2.266);
\draw[gp path] (5.769,2.266)--(5.589,2.266);
\node[gp node right] at (0.790,2.266) {$400$};
\draw[gp path] (0.900,3.103)--(1.080,3.103);
\draw[gp path] (5.769,3.103)--(5.589,3.103);
\node[gp node right] at (0.790,3.103) {$600$};
\draw[gp path] (0.900,3.939)--(1.080,3.939);
\draw[gp path] (5.769,3.939)--(5.589,3.939);
\node[gp node right] at (0.790,3.939) {$800$};
\draw[gp path] (0.900,4.776)--(1.080,4.776);
\draw[gp path] (5.769,4.776)--(5.589,4.776);
\node[gp node right] at (0.790,4.776) {$1000$};
\draw[gp path] (1.275,0.592)--(1.275,0.772);
\draw[gp path] (1.275,4.814)--(1.275,4.634);
\node[gp node center] at (1.275,0.407) {$2$};
\draw[gp path] (2.024,0.592)--(2.024,0.772);
\draw[gp path] (2.024,4.814)--(2.024,4.634);
\node[gp node center] at (2.024,0.407) {$4$};
\draw[gp path] (2.773,0.592)--(2.773,0.772);
\draw[gp path] (2.773,4.814)--(2.773,4.634);
\node[gp node center] at (2.773,0.407) {$6$};
\draw[gp path] (3.522,0.592)--(3.522,0.772);
\draw[gp path] (3.522,4.814)--(3.522,4.634);
\node[gp node center] at (3.522,0.407) {$8$};
\draw[gp path] (4.271,0.592)--(4.271,0.772);
\draw[gp path] (4.271,4.814)--(4.271,4.634);
\node[gp node center] at (4.271,0.407) {$10$};
\draw[gp path] (5.020,0.592)--(5.020,0.772);
\draw[gp path] (5.020,4.814)--(5.020,4.634);
\node[gp node center] at (5.020,0.407) {$12$};
\draw[gp path] (5.769,0.592)--(5.769,0.772);
\draw[gp path] (5.769,4.814)--(5.769,4.634);
\node[gp node center] at (5.769,0.407) {$14$};
\draw[gp path] (0.900,4.814)--(0.900,0.592)--(5.769,0.592)--(5.769,4.814)--cycle;
\node[gp node left,rotate=40] at (2.024,2.098) {peak Flop/s};
\node[gp node center,rotate=-270] at (0.175,2.703) {GFlop/s};
\node[gp node center] at (3.334,0.130) {\# cores};
\node[gp node right] at (2.990,4.521) {$n\times200$ ($I_c = 16.7$)};
\gpcolor{rgb color={0.580,0.000,0.827}}
\gpsetdashtype{gp dt 1}
\gpsetlinewidth{2.00}
\draw[gp path] (3.100,4.521)--(3.720,4.521);
\draw[gp path] (0.900,0.783)--(1.275,0.962)--(1.649,1.134)--(2.024,1.303)--(2.398,1.420)%
  --(2.773,1.593)--(3.147,1.760)--(3.522,1.897)--(3.896,1.943)--(4.271,2.099)--(4.645,2.264)%
  --(5.020,2.383)--(5.394,2.485)--(5.769,2.644);
\gpcolor{color=gp lt color border}
\node[gp node right] at (2.990,4.296) {$n\times500$ ($I_c = 41.7$)};
\gpcolor{rgb color={0.000,0.620,0.451}}
\gpsetdashtype{gp dt 2}
\draw[gp path] (3.100,4.296)--(3.720,4.296);
\draw[gp path] (0.900,0.824)--(1.275,1.050)--(1.649,1.267)--(2.024,1.469)--(2.398,1.573)%
  --(2.773,1.759)--(3.147,2.001)--(3.522,2.107)--(3.896,2.132)--(4.271,2.312)--(4.645,2.490)%
  --(5.020,2.635)--(5.394,2.753)--(5.769,2.916);
\gpcolor{color=gp lt color border}
\node[gp node right] at (2.990,4.071) {$n\times1000$ ($I_c = 83.3$)};
\gpcolor{rgb color={0.337,0.706,0.914}}
\gpsetdashtype{gp dt 3}
\draw[gp path] (3.100,4.071)--(3.720,4.071);
\draw[gp path] (0.900,0.796)--(1.275,0.992)--(1.649,1.184)--(2.024,1.318)--(2.398,1.433)%
  --(2.773,1.622)--(3.147,1.723)--(3.522,1.839)--(3.896,1.925)--(4.271,2.001)--(4.645,2.129)%
  --(5.020,1.977)--(5.394,2.385)--(5.769,2.576);
\gpcolor{color=gp lt color border}
\gpsetlinetype{gp lt axes}
\gpsetdashtype{gp dt axes}
\draw[gp path] (0.900,0.893)--(0.949,0.933)--(0.998,0.972)--(1.048,1.012)--(1.097,1.052)%
  --(1.146,1.091)--(1.195,1.131)--(1.244,1.170)--(1.293,1.210)--(1.343,1.249)--(1.392,1.289)%
  --(1.441,1.328)--(1.490,1.368)--(1.539,1.408)--(1.589,1.447)--(1.638,1.487)--(1.687,1.526)%
  --(1.736,1.566)--(1.785,1.605)--(1.834,1.645)--(1.884,1.684)--(1.933,1.724)--(1.982,1.764)%
  --(2.031,1.803)--(2.080,1.843)--(2.130,1.882)--(2.179,1.922)--(2.228,1.961)--(2.277,2.001)%
  --(2.326,2.041)--(2.375,2.080)--(2.425,2.120)--(2.474,2.159)--(2.523,2.199)--(2.572,2.238)%
  --(2.621,2.278)--(2.671,2.317)--(2.720,2.357)--(2.769,2.397)--(2.818,2.436)--(2.867,2.476)%
  --(2.916,2.515)--(2.966,2.555)--(3.015,2.594)--(3.064,2.634)--(3.113,2.674)--(3.162,2.713)%
  --(3.212,2.753)--(3.261,2.792)--(3.310,2.832)--(3.359,2.871)--(3.408,2.911)--(3.457,2.950)%
  --(3.507,2.990)--(3.556,3.030)--(3.605,3.069)--(3.654,3.109)--(3.703,3.148)--(3.753,3.188)%
  --(3.802,3.227)--(3.851,3.267)--(3.900,3.306)--(3.949,3.346)--(3.998,3.386)--(4.048,3.425)%
  --(4.097,3.465)--(4.146,3.504)--(4.195,3.544)--(4.244,3.583)--(4.294,3.623)--(4.343,3.663)%
  --(4.392,3.702)--(4.441,3.742)--(4.490,3.781)--(4.539,3.821)--(4.589,3.860)--(4.638,3.900)%
  --(4.687,3.939)--(4.736,3.979)--(4.785,4.019)--(4.835,4.058)--(4.884,4.098)--(4.933,4.137)%
  --(4.982,4.177)--(5.031,4.216)--(5.080,4.256)--(5.130,4.296)--(5.179,4.335)--(5.228,4.375)%
  --(5.277,4.414)--(5.326,4.454)--(5.376,4.493)--(5.425,4.533)--(5.474,4.572)--(5.523,4.612)%
  --(5.572,4.652)--(5.621,4.691)--(5.671,4.731)--(5.720,4.770)--(5.769,4.810);
\gpsetlinetype{gp lt border}
\gpsetdashtype{gp dt solid}
\gpsetlinewidth{1.00}
\draw[gp path] (0.900,4.814)--(0.900,0.592)--(5.769,0.592)--(5.769,4.814)--cycle;
\gpdefrectangularnode{gp plot 1}{\pgfpoint{0.900cm}{0.592cm}}{\pgfpoint{5.769cm}{4.814cm}}
\end{tikzpicture}
\subcaption{Compute-bound case ($I_c \gtrsim 14$) measured with $n=5\cdot10^6$.}\label{fig:TSMM_flops}
\end{subfigure}
\caption{Single socket TSMM+reshape compared with the peak bandwidth respectively peak Flop/s.
The input matrices have dimensions $n\times m$ and $m\times m/2$, the result is reshaped to $n/2 \times m$.
Based on \autoref{tab:peak_bandwidth_and_flops}, the machine intensity for this operation (load/store ratio of 2/1 $\widehat{=}$ stream) is $1009/73\approx 14$~[Flops/Byte].}
\end{figure}
\autoref{fig:TSMM_bandwidth} shows the obtained bandwidth for varying numbers of cores.
We observe a saturation of the memory bandwidth for $m<50$.
For $m=100$, we already see a linear scaling with the number of cores.
For the compute-bound case, our implementation roughly obtains 50\% of the peak performance (see \autoref{fig:TSMM_flops}).
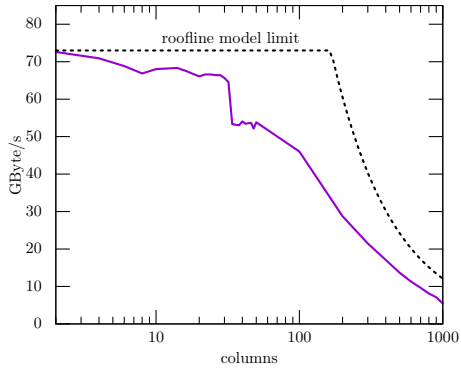
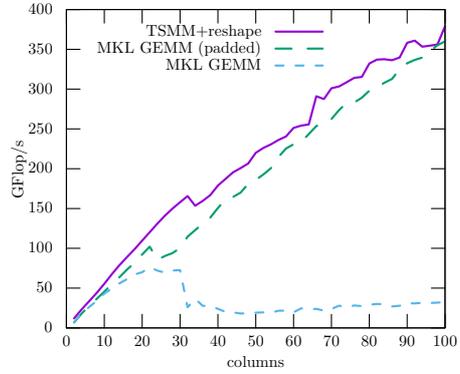
\begin{figure}
\begin{subfigure}[t]{0.47\textwidth}
\begin{tikzpicture}[gnuplot]
\tikzset{every node/.append style={scale=0.60}}
\path (0.000,0.000) rectangle (6.100,5.000);
\gpcolor{color=gp lt color border}
\gpsetlinetype{gp lt border}
\gpsetdashtype{gp dt solid}
\gpsetlinewidth{1.00}
\draw[gp path] (0.680,0.592)--(0.860,0.592);
\draw[gp path] (5.769,0.592)--(5.589,0.592);
\node[gp node right] at (0.570,0.592) {$0$};
\draw[gp path] (0.680,1.089)--(0.860,1.089);
\draw[gp path] (5.769,1.089)--(5.589,1.089);
\node[gp node right] at (0.570,1.089) {$10$};
\draw[gp path] (0.680,1.585)--(0.860,1.585);
\draw[gp path] (5.769,1.585)--(5.589,1.585);
\node[gp node right] at (0.570,1.585) {$20$};
\draw[gp path] (0.680,2.082)--(0.860,2.082);
\draw[gp path] (5.769,2.082)--(5.589,2.082);
\node[gp node right] at (0.570,2.082) {$30$};
\draw[gp path] (0.680,2.579)--(0.860,2.579);
\draw[gp path] (5.769,2.579)--(5.589,2.579);
\node[gp node right] at (0.570,2.579) {$40$};
\draw[gp path] (0.680,3.076)--(0.860,3.076);
\draw[gp path] (5.769,3.076)--(5.589,3.076);
\node[gp node right] at (0.570,3.076) {$50$};
\draw[gp path] (0.680,3.572)--(0.860,3.572);
\draw[gp path] (5.769,3.572)--(5.589,3.572);
\node[gp node right] at (0.570,3.572) {$60$};
\draw[gp path] (0.680,4.069)--(0.860,4.069);
\draw[gp path] (5.769,4.069)--(5.589,4.069);
\node[gp node right] at (0.570,4.069) {$70$};
\draw[gp path] (0.680,4.566)--(0.860,4.566);
\draw[gp path] (5.769,4.566)--(5.589,4.566);
\node[gp node right] at (0.570,4.566) {$80$};
\draw[gp path] (0.680,0.592)--(0.680,0.682);
\draw[gp path] (0.680,4.814)--(0.680,4.724);
\draw[gp path] (1.012,0.592)--(1.012,0.682);
\draw[gp path] (1.012,4.814)--(1.012,4.724);
\draw[gp path] (1.248,0.592)--(1.248,0.682);
\draw[gp path] (1.248,4.814)--(1.248,4.724);
\draw[gp path] (1.430,0.592)--(1.430,0.682);
\draw[gp path] (1.430,4.814)--(1.430,4.724);
\draw[gp path] (1.580,0.592)--(1.580,0.682);
\draw[gp path] (1.580,4.814)--(1.580,4.724);
\draw[gp path] (1.706,0.592)--(1.706,0.682);
\draw[gp path] (1.706,4.814)--(1.706,4.724);
\draw[gp path] (1.815,0.592)--(1.815,0.682);
\draw[gp path] (1.815,4.814)--(1.815,4.724);
\draw[gp path] (1.912,0.592)--(1.912,0.682);
\draw[gp path] (1.912,4.814)--(1.912,4.724);
\draw[gp path] (1.998,0.592)--(1.998,0.772);
\draw[gp path] (1.998,4.814)--(1.998,4.634);
\node[gp node center] at (1.998,0.407) {$10$};
\draw[gp path] (2.566,0.592)--(2.566,0.682);
\draw[gp path] (2.566,4.814)--(2.566,4.724);
\draw[gp path] (2.898,0.592)--(2.898,0.682);
\draw[gp path] (2.898,4.814)--(2.898,4.724);
\draw[gp path] (3.133,0.592)--(3.133,0.682);
\draw[gp path] (3.133,4.814)--(3.133,4.724);
\draw[gp path] (3.316,0.592)--(3.316,0.682);
\draw[gp path] (3.316,4.814)--(3.316,4.724);
\draw[gp path] (3.465,0.592)--(3.465,0.682);
\draw[gp path] (3.465,4.814)--(3.465,4.724);
\draw[gp path] (3.591,0.592)--(3.591,0.682);
\draw[gp path] (3.591,4.814)--(3.591,4.724);
\draw[gp path] (3.701,0.592)--(3.701,0.682);
\draw[gp path] (3.701,4.814)--(3.701,4.724);
\draw[gp path] (3.797,0.592)--(3.797,0.682);
\draw[gp path] (3.797,4.814)--(3.797,4.724);
\draw[gp path] (3.883,0.592)--(3.883,0.772);
\draw[gp path] (3.883,4.814)--(3.883,4.634);
\node[gp node center] at (3.883,0.407) {$100$};
\draw[gp path] (4.451,0.592)--(4.451,0.682);
\draw[gp path] (4.451,4.814)--(4.451,4.724);
\draw[gp path] (4.783,0.592)--(4.783,0.682);
\draw[gp path] (4.783,4.814)--(4.783,4.724);
\draw[gp path] (5.019,0.592)--(5.019,0.682);
\draw[gp path] (5.019,4.814)--(5.019,4.724);
\draw[gp path] (5.201,0.592)--(5.201,0.682);
\draw[gp path] (5.201,4.814)--(5.201,4.724);
\draw[gp path] (5.351,0.592)--(5.351,0.682);
\draw[gp path] (5.351,4.814)--(5.351,4.724);
\draw[gp path] (5.477,0.592)--(5.477,0.682);
\draw[gp path] (5.477,4.814)--(5.477,4.724);
\draw[gp path] (5.586,0.592)--(5.586,0.682);
\draw[gp path] (5.586,4.814)--(5.586,4.724);
\draw[gp path] (5.683,0.592)--(5.683,0.682);
\draw[gp path] (5.683,4.814)--(5.683,4.724);
\draw[gp path] (5.769,0.592)--(5.769,0.772);
\draw[gp path] (5.769,4.814)--(5.769,4.634);
\node[gp node center] at (5.769,0.407) {$1000$};
\draw[gp path] (0.680,4.814)--(0.680,0.592)--(5.769,0.592)--(5.769,4.814)--cycle;
\node[gp node left] at (1.998,4.367) {roofline model limit};
\node[gp node center,rotate=-270] at (0.175,2.703) {GByte/s};
\node[gp node center] at (3.224,0.130) {columns};
\gpcolor{rgb color={0.580,0.000,0.827}}
\gpsetdashtype{gp dt 1}
\gpsetlinewidth{2.00}
\draw[gp path] (0.680,4.198)--(1.248,4.114)--(1.580,4.011)--(1.815,3.913)--(1.998,3.971)%
  --(2.147,3.979)--(2.273,3.987)--(2.383,3.950)--(2.479,3.909)--(2.566,3.875)--(2.644,3.900)%
  --(2.715,3.900)--(2.780,3.891)--(2.841,3.891)--(2.898,3.853)--(2.950,3.798)--(3.000,3.241)%
  --(3.047,3.233)--(3.091,3.227)--(3.133,3.278)--(3.173,3.246)--(3.211,3.253)--(3.248,3.258)%
  --(3.282,3.183)--(3.316,3.266)--(3.883,2.878)--(4.451,2.020)--(4.783,1.659)--(5.019,1.441)%
  --(5.201,1.271)--(5.351,1.154)--(5.477,1.071)--(5.586,0.995)--(5.683,0.946)--(5.769,0.863);
\gpcolor{color=gp lt color border}
\gpsetlinetype{gp lt axes}
\gpsetdashtype{gp dt axes}
\draw[gp path] (0.680,4.218)--(0.731,4.218)--(0.783,4.218)--(0.834,4.218)--(0.886,4.218)%
  --(0.937,4.218)--(0.988,4.218)--(1.040,4.218)--(1.091,4.218)--(1.143,4.218)--(1.194,4.218)%
  --(1.245,4.218)--(1.297,4.218)--(1.348,4.218)--(1.400,4.218)--(1.451,4.218)--(1.502,4.218)%
  --(1.554,4.218)--(1.605,4.218)--(1.657,4.218)--(1.708,4.218)--(1.759,4.218)--(1.811,4.218)%
  --(1.862,4.218)--(1.914,4.218)--(1.965,4.218)--(2.017,4.218)--(2.068,4.218)--(2.119,4.218)%
  --(2.171,4.218)--(2.222,4.218)--(2.274,4.218)--(2.325,4.218)--(2.376,4.218)--(2.428,4.218)%
  --(2.479,4.218)--(2.531,4.218)--(2.582,4.218)--(2.633,4.218)--(2.685,4.218)--(2.736,4.218)%
  --(2.788,4.218)--(2.839,4.218)--(2.890,4.218)--(2.942,4.218)--(2.993,4.218)--(3.045,4.218)%
  --(3.096,4.218)--(3.147,4.218)--(3.199,4.218)--(3.250,4.218)--(3.302,4.218)--(3.353,4.218)%
  --(3.404,4.218)--(3.456,4.218)--(3.507,4.218)--(3.559,4.218)--(3.610,4.218)--(3.661,4.218)%
  --(3.713,4.218)--(3.764,4.218)--(3.816,4.218)--(3.867,4.218)--(3.918,4.218)--(3.970,4.218)%
  --(4.021,4.218)--(4.073,4.218)--(4.124,4.218)--(4.175,4.218)--(4.227,4.218)--(4.278,4.218)%
  --(4.330,4.080)--(4.381,3.867)--(4.432,3.668)--(4.484,3.481)--(4.535,3.305)--(4.587,3.140)%
  --(4.638,2.985)--(4.690,2.839)--(4.741,2.703)--(4.792,2.574)--(4.844,2.454)--(4.895,2.340)%
  --(4.947,2.234)--(4.998,2.134)--(5.049,2.040)--(5.101,1.952)--(5.152,1.869)--(5.204,1.792)%
  --(5.255,1.719)--(5.306,1.650)--(5.358,1.586)--(5.409,1.525)--(5.461,1.468)--(5.512,1.415)%
  --(5.563,1.365)--(5.615,1.318)--(5.666,1.274)--(5.718,1.232)--(5.769,1.193);
\gpsetlinetype{gp lt border}
\gpsetdashtype{gp dt solid}
\gpsetlinewidth{1.00}
\draw[gp path] (0.680,4.814)--(0.680,0.592)--(5.769,0.592)--(5.769,4.814)--cycle;
\gpdefrectangularnode{gp plot 1}{\pgfpoint{0.680cm}{0.592cm}}{\pgfpoint{5.769cm}{4.814cm}}
\end{tikzpicture}
\subcaption{Obtained memory bandwidth of our TSMM+reshape implementation compared to the roofline limit: $b_\text{roofline} =\min(P_\text{peak}/I_c, b_s)$.}\label{fig:TSMM_Roofline_bandwidth}
\end{subfigure}
\hfill
\begin{subfigure}[t]{0.47\textwidth}
\begin{tikzpicture}[gnuplot]
\tikzset{every node/.append style={scale=0.60}}
\path (0.000,0.000) rectangle (6.100,5.000);
\gpcolor{color=gp lt color border}
\gpsetlinetype{gp lt border}
\gpsetdashtype{gp dt solid}
\gpsetlinewidth{1.00}
\draw[gp path] (0.790,0.592)--(0.970,0.592);
\draw[gp path] (5.769,0.592)--(5.589,0.592);
\node[gp node right] at (0.680,0.592) {$0$};
\draw[gp path] (0.790,1.120)--(0.970,1.120);
\draw[gp path] (5.769,1.120)--(5.589,1.120);
\node[gp node right] at (0.680,1.120) {$50$};
\draw[gp path] (0.790,1.648)--(0.970,1.648);
\draw[gp path] (5.769,1.648)--(5.589,1.648);
\node[gp node right] at (0.680,1.648) {$100$};
\draw[gp path] (0.790,2.175)--(0.970,2.175);
\draw[gp path] (5.769,2.175)--(5.589,2.175);
\node[gp node right] at (0.680,2.175) {$150$};
\draw[gp path] (0.790,2.703)--(0.970,2.703);
\draw[gp path] (5.769,2.703)--(5.589,2.703);
\node[gp node right] at (0.680,2.703) {$200$};
\draw[gp path] (0.790,3.231)--(0.970,3.231);
\draw[gp path] (5.769,3.231)--(5.589,3.231);
\node[gp node right] at (0.680,3.231) {$250$};
\draw[gp path] (0.790,3.759)--(0.970,3.759);
\draw[gp path] (5.769,3.759)--(5.589,3.759);
\node[gp node right] at (0.680,3.759) {$300$};
\draw[gp path] (0.790,4.286)--(0.970,4.286);
\draw[gp path] (5.769,4.286)--(5.589,4.286);
\node[gp node right] at (0.680,4.286) {$350$};
\draw[gp path] (0.790,4.814)--(0.970,4.814);
\draw[gp path] (5.769,4.814)--(5.589,4.814);
\node[gp node right] at (0.680,4.814) {$400$};
\draw[gp path] (0.790,0.592)--(0.790,0.772);
\draw[gp path] (0.790,4.814)--(0.790,4.634);
\node[gp node center] at (0.790,0.407) {$0$};
\draw[gp path] (1.288,0.592)--(1.288,0.772);
\draw[gp path] (1.288,4.814)--(1.288,4.634);
\node[gp node center] at (1.288,0.407) {$10$};
\draw[gp path] (1.786,0.592)--(1.786,0.772);
\draw[gp path] (1.786,4.814)--(1.786,4.634);
\node[gp node center] at (1.786,0.407) {$20$};
\draw[gp path] (2.284,0.592)--(2.284,0.772);
\draw[gp path] (2.284,4.814)--(2.284,4.634);
\node[gp node center] at (2.284,0.407) {$30$};
\draw[gp path] (2.782,0.592)--(2.782,0.772);
\draw[gp path] (2.782,4.814)--(2.782,4.634);
\node[gp node center] at (2.782,0.407) {$40$};
\draw[gp path] (3.280,0.592)--(3.280,0.772);
\draw[gp path] (3.280,4.814)--(3.280,4.634);
\node[gp node center] at (3.280,0.407) {$50$};
\draw[gp path] (3.777,0.592)--(3.777,0.772);
\draw[gp path] (3.777,4.814)--(3.777,4.634);
\node[gp node center] at (3.777,0.407) {$60$};
\draw[gp path] (4.275,0.592)--(4.275,0.772);
\draw[gp path] (4.275,4.814)--(4.275,4.634);
\node[gp node center] at (4.275,0.407) {$70$};
\draw[gp path] (4.773,0.592)--(4.773,0.772);
\draw[gp path] (4.773,4.814)--(4.773,4.634);
\node[gp node center] at (4.773,0.407) {$80$};
\draw[gp path] (5.271,0.592)--(5.271,0.772);
\draw[gp path] (5.271,4.814)--(5.271,4.634);
\node[gp node center] at (5.271,0.407) {$90$};
\draw[gp path] (5.769,0.592)--(5.769,0.772);
\draw[gp path] (5.769,4.814)--(5.769,4.634);
\node[gp node center] at (5.769,0.407) {$100$};
\draw[gp path] (0.790,4.814)--(0.790,0.592)--(5.769,0.592)--(5.769,4.814)--cycle;
\node[gp node center,rotate=-270] at (0.175,2.703) {GFlop/s};
\node[gp node center] at (3.279,0.130) {columns};
\node[gp node right] at (3.430,4.521) {TSMM+reshape};
\gpcolor{rgb color={0.580,0.000,0.827}}
\gpsetdashtype{gp dt 1}
\gpsetlinewidth{2.00}
\draw[gp path] (3.540,4.521)--(4.160,4.521);
\draw[gp path] (0.890,0.718)--(0.989,0.837)--(1.089,0.946)--(1.188,1.056)--(1.288,1.177)%
  --(1.387,1.304)--(1.487,1.424)--(1.587,1.532)--(1.686,1.636)--(1.786,1.750)--(1.885,1.861)%
  --(1.985,1.974)--(2.085,2.081)--(2.184,2.174)--(2.284,2.259)--(2.383,2.340)--(2.483,2.212)%
  --(2.582,2.276)--(2.682,2.354)--(2.782,2.480)--(2.881,2.567)--(2.981,2.655)--(3.080,2.711)%
  --(3.180,2.772)--(3.280,2.914)--(3.379,2.978)--(3.479,3.026)--(3.578,3.084)--(3.678,3.131)%
  --(3.777,3.244)--(3.877,3.275)--(3.977,3.290)--(4.076,3.664)--(4.176,3.627)--(4.275,3.770)%
  --(4.375,3.794)--(4.474,3.849)--(4.574,3.908)--(4.674,3.920)--(4.773,4.098)--(4.873,4.150)%
  --(4.972,4.155)--(5.072,4.142)--(5.172,4.178)--(5.271,4.373)--(5.371,4.402)--(5.470,4.321)%
  --(5.570,4.335)--(5.669,4.350)--(5.769,4.591);
\gpcolor{color=gp lt color border}
\node[gp node right] at (3.430,4.296) {      MKL GEMM (padded)};
\gpcolor{rgb color={0.000,0.620,0.451}}
\gpsetdashtype{gp dt 2}
\draw[gp path] (3.540,4.296)--(4.160,4.296);
\draw[gp path] (0.890,0.664)--(0.989,0.782)--(1.089,0.884)--(1.188,0.977)--(1.288,1.073)%
  --(1.387,1.171)--(1.487,1.264)--(1.587,1.366)--(1.686,1.456)--(1.786,1.566)--(1.885,1.670)%
  --(1.985,1.494)--(2.085,1.546)--(2.184,1.584)--(2.284,1.651)--(2.383,1.802)--(2.483,1.889)%
  --(2.582,1.974)--(2.682,2.056)--(2.782,2.182)--(2.881,2.293)--(2.981,2.329)--(3.080,2.386)%
  --(3.180,2.500)--(3.280,2.555)--(3.379,2.626)--(3.479,2.710)--(3.578,2.858)--(3.678,2.971)%
  --(3.777,3.027)--(3.877,3.057)--(3.977,3.163)--(4.076,3.268)--(4.176,3.311)--(4.275,3.362)%
  --(4.375,3.480)--(4.474,3.576)--(4.574,3.581)--(4.674,3.644)--(4.773,3.737)--(4.873,3.785)%
  --(4.972,3.843)--(5.072,3.894)--(5.172,4.024)--(5.271,4.102)--(5.371,4.146)--(5.470,4.177)%
  --(5.570,4.264)--(5.669,4.340)--(5.769,4.392);
\gpcolor{color=gp lt color border}
\node[gp node right] at (3.430,4.071) {MKL GEMM};
\gpcolor{rgb color={0.337,0.706,0.914}}
\gpsetdashtype{gp dt 3}
\draw[gp path] (3.540,4.071)--(4.160,4.071);
\draw[gp path] (0.890,0.666)--(0.989,0.779)--(1.089,0.870)--(1.188,0.950)--(1.288,1.047)%
  --(1.387,1.114)--(1.487,1.183)--(1.587,1.239)--(1.686,1.300)--(1.786,1.331)--(1.885,1.398)%
  --(1.985,1.352)--(2.085,1.325)--(2.184,1.353)--(2.284,1.359)--(2.383,0.866)--(2.483,0.978)%
  --(2.582,0.885)--(2.682,0.882)--(2.782,0.845)--(2.881,0.805)--(2.981,0.798)--(3.080,0.783)%
  --(3.180,0.786)--(3.280,0.794)--(3.379,0.798)--(3.479,0.793)--(3.578,0.820)--(3.678,0.821)%
  --(3.777,0.799)--(3.877,0.852)--(3.977,0.837)--(4.076,0.845)--(4.176,0.821)--(4.275,0.843)%
  --(4.375,0.887)--(4.474,0.866)--(4.574,0.890)--(4.674,0.881)--(4.773,0.900)--(4.873,0.910)%
  --(4.972,0.906)--(5.072,0.877)--(5.172,0.886)--(5.271,0.919)--(5.371,0.917)--(5.470,0.925)%
  --(5.570,0.923)--(5.669,0.927)--(5.769,0.937);
\gpcolor{color=gp lt color border}
\gpsetdashtype{gp dt solid}
\gpsetlinewidth{1.00}
\draw[gp path] (0.790,4.814)--(0.790,0.592)--(5.769,0.592)--(5.769,4.814)--cycle;
\gpdefrectangularnode{gp plot 1}{\pgfpoint{0.790cm}{0.592cm}}{\pgfpoint{5.769cm}{4.814cm}}
\end{tikzpicture}
\subcaption{Obtained performance compared to the Intel MKL 2020.3 for $n=2^{24}$ rows with and without padding to $n^{24}+64$.}
\end{subfigure}
\caption{Single socket performance of TSMM+reshape for varying numbers of columns.
The input matrices have dimensions $n\times m$ and $m\times m/2$, and our implementation directly stores the result reshaped to dimensions $n/2\times m$.
(b) illustrates the effect of cache thrashing (the leading array dimension is a power of two).}\label{fig:TSMM_vs_MKL}
\end{figure}
From \autoref{fig:TSMM_Roofline_bandwidth}, we conclude that our TSMM implementation obtains a high fraction of the maximum possible bandwidth.
Near the kink, the Roofline model is too optimistic because it assumes that data transfers and floating point operations overlap perfectly.
Further insight could be obtained by a more sophisticated performance model such as the ECM (Execution-Cache-Memory) model, see~\cite{Stengel2015}.
For this operation, our implementation and the Intel MKL obtain roughly the same performance, as depicted in \autoref{fig:TSMM_vs_MKL}.
In contrast to the MKL, our implementation exploits a special memory layout, which might explain the small differences in run-time.
So the advantage of our TSMM implementation for the complete TT-SVD algorithm consists mainly in fusing the $\operatorname{reshape}$ operation,
which ensures a suitably padded memory layout for subsequent operations at no additional cost.
Without appropriate padding, the performance can degrade significantly due to cache thrashing (also illustrated in \autoref{fig:TSMM_vs_MKL}),
in particular for operations from tensor algorithms when individual dimensions are multiples of two.

\subsection{TT-SVD}
\begin{figure}
\begin{subfigure}[t]{0.47\textwidth}
\input{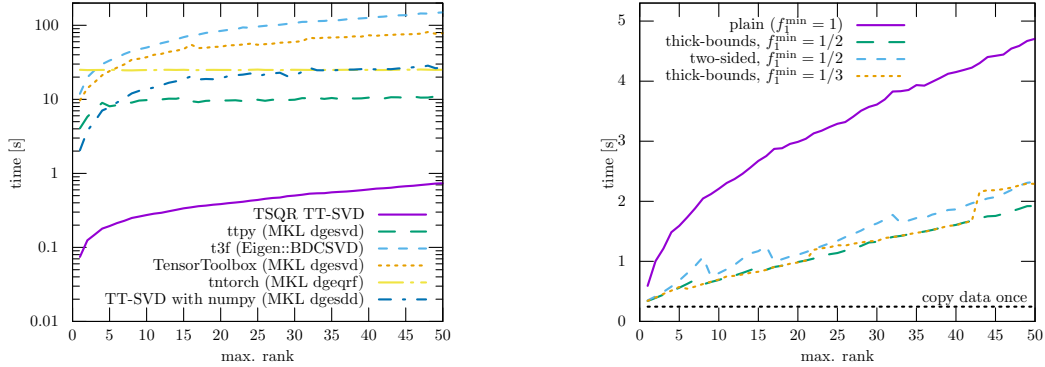}
\subcaption{Different implementations of the classical TT-SVD algorithm for a $2^{27}$ tensor and our TSQR TT-SVD algorithm (without combining dimensions).}\label{fig:TT-SVD_other_libraries}
\end{subfigure}
\hfill
\begin{subfigure}[t]{0.47\textwidth}
\begin{tikzpicture}[gnuplot]
\tikzset{every node/.append style={scale=0.60}}
\path (0.000,0.000) rectangle (6.100,5.000);
\gpcolor{color=gp lt color border}
\gpsetlinetype{gp lt border}
\gpsetdashtype{gp dt solid}
\gpsetlinewidth{1.00}
\draw[gp path] (0.570,0.592)--(0.750,0.592);
\draw[gp path] (5.769,0.592)--(5.589,0.592);
\node[gp node right] at (0.460,0.592) {$0$};
\draw[gp path] (0.570,1.389)--(0.750,1.389);
\draw[gp path] (5.769,1.389)--(5.589,1.389);
\node[gp node right] at (0.460,1.389) {$1$};
\draw[gp path] (0.570,2.185)--(0.750,2.185);
\draw[gp path] (5.769,2.185)--(5.589,2.185);
\node[gp node right] at (0.460,2.185) {$2$};
\draw[gp path] (0.570,2.982)--(0.750,2.982);
\draw[gp path] (5.769,2.982)--(5.589,2.982);
\node[gp node right] at (0.460,2.982) {$3$};
\draw[gp path] (0.570,3.778)--(0.750,3.778);
\draw[gp path] (5.769,3.778)--(5.589,3.778);
\node[gp node right] at (0.460,3.778) {$4$};
\draw[gp path] (0.570,4.575)--(0.750,4.575);
\draw[gp path] (5.769,4.575)--(5.589,4.575);
\node[gp node right] at (0.460,4.575) {$5$};
\draw[gp path] (0.570,0.592)--(0.570,0.772);
\draw[gp path] (0.570,4.814)--(0.570,4.634);
\node[gp node center] at (0.570,0.407) {$0$};
\draw[gp path] (1.090,0.592)--(1.090,0.772);
\draw[gp path] (1.090,4.814)--(1.090,4.634);
\node[gp node center] at (1.090,0.407) {$5$};
\draw[gp path] (1.610,0.592)--(1.610,0.772);
\draw[gp path] (1.610,4.814)--(1.610,4.634);
\node[gp node center] at (1.610,0.407) {$10$};
\draw[gp path] (2.130,0.592)--(2.130,0.772);
\draw[gp path] (2.130,4.814)--(2.130,4.634);
\node[gp node center] at (2.130,0.407) {$15$};
\draw[gp path] (2.650,0.592)--(2.650,0.772);
\draw[gp path] (2.650,4.814)--(2.650,4.634);
\node[gp node center] at (2.650,0.407) {$20$};
\draw[gp path] (3.170,0.592)--(3.170,0.772);
\draw[gp path] (3.170,4.814)--(3.170,4.634);
\node[gp node center] at (3.170,0.407) {$25$};
\draw[gp path] (3.689,0.592)--(3.689,0.772);
\draw[gp path] (3.689,4.814)--(3.689,4.634);
\node[gp node center] at (3.689,0.407) {$30$};
\draw[gp path] (4.209,0.592)--(4.209,0.772);
\draw[gp path] (4.209,4.814)--(4.209,4.634);
\node[gp node center] at (4.209,0.407) {$35$};
\draw[gp path] (4.729,0.592)--(4.729,0.772);
\draw[gp path] (4.729,4.814)--(4.729,4.634);
\node[gp node center] at (4.729,0.407) {$40$};
\draw[gp path] (5.249,0.592)--(5.249,0.772);
\draw[gp path] (5.249,4.814)--(5.249,4.634);
\node[gp node center] at (5.249,0.407) {$45$};
\draw[gp path] (5.769,0.592)--(5.769,0.772);
\draw[gp path] (5.769,4.814)--(5.769,4.634);
\node[gp node center] at (5.769,0.407) {$50$};
\draw[gp path] (0.570,4.814)--(0.570,0.592)--(5.769,0.592)--(5.769,4.814)--cycle;
\node[gp node left] at (4.209,0.911) {copy data once};
\node[gp node center,rotate=-270] at (0.175,2.703) {time [s]};
\node[gp node center] at (3.169,0.130) {max. rank};
\node[gp node right] at (3.320,4.521) {plain ($f_1^\text{min}=1$)};
\gpcolor{rgb color={0.580,0.000,0.827}}
\gpsetdashtype{gp dt 1}
\gpsetlinewidth{2.00}
\draw[gp path] (3.430,4.521)--(4.050,4.521);
\draw[gp path] (0.674,1.064)--(0.778,1.392)--(0.882,1.559)--(0.986,1.775)--(1.090,1.860)%
  --(1.194,1.972)--(1.298,2.091)--(1.402,2.222)--(1.506,2.280)--(1.610,2.351)--(1.714,2.428)%
  --(1.818,2.481)--(1.922,2.550)--(2.026,2.633)--(2.130,2.721)--(2.234,2.785)--(2.338,2.880)%
  --(2.442,2.889)--(2.546,2.946)--(2.650,2.973)--(2.754,3.014)--(2.858,3.086)--(2.962,3.121)%
  --(3.066,3.170)--(3.170,3.213)--(3.273,3.236)--(3.377,3.307)--(3.481,3.386)--(3.585,3.437)%
  --(3.689,3.468)--(3.793,3.536)--(3.897,3.641)--(4.001,3.645)--(4.105,3.661)--(4.209,3.726)%
  --(4.313,3.719)--(4.417,3.770)--(4.521,3.821)--(4.625,3.876)--(4.729,3.899)--(4.833,3.927)%
  --(4.937,3.957)--(5.041,4.018)--(5.145,4.098)--(5.249,4.114)--(5.353,4.131)--(5.457,4.210)%
  --(5.561,4.246)--(5.665,4.311)--(5.769,4.340);
\gpcolor{color=gp lt color border}
\node[gp node right] at (3.320,4.296) {thick-bounds, $f_1^\text{min}=1/2$};
\gpcolor{rgb color={0.000,0.620,0.451}}
\gpsetdashtype{gp dt 2}
\draw[gp path] (3.430,4.296)--(4.050,4.296);
\draw[gp path] (0.674,0.863)--(0.778,0.904)--(0.882,0.947)--(0.986,0.993)--(1.090,1.039)%
  --(1.194,1.090)--(1.298,1.144)--(1.402,1.083)--(1.506,1.115)--(1.610,1.142)--(1.714,1.176)%
  --(1.818,1.201)--(1.922,1.235)--(2.026,1.273)--(2.130,1.311)--(2.234,1.277)--(2.338,1.314)%
  --(2.442,1.331)--(2.546,1.357)--(2.650,1.380)--(2.754,1.400)--(2.858,1.426)--(2.962,1.457)%
  --(3.066,1.483)--(3.170,1.499)--(3.273,1.531)--(3.377,1.558)--(3.481,1.593)--(3.585,1.639)%
  --(3.689,1.648)--(3.793,1.688)--(3.897,1.712)--(4.001,1.726)--(4.105,1.746)--(4.209,1.766)%
  --(4.313,1.784)--(4.417,1.805)--(4.521,1.832)--(4.625,1.855)--(4.729,1.872)--(4.833,1.902)%
  --(4.937,1.926)--(5.041,1.961)--(5.145,1.986)--(5.249,1.991)--(5.353,2.019)--(5.457,2.053)%
  --(5.561,2.085)--(5.665,2.121)--(5.769,2.123);
\gpcolor{color=gp lt color border}
\node[gp node right] at (3.320,4.071) {two-sided, $f_1^\text{min}=1/2$};
\gpcolor{rgb color={0.337,0.706,0.914}}
\gpsetdashtype{gp dt 3}
\draw[gp path] (3.430,4.071)--(4.050,4.071);
\draw[gp path] (0.674,0.871)--(0.778,0.923)--(0.882,0.984)--(0.986,1.054)--(1.090,1.135)%
  --(1.194,1.230)--(1.298,1.328)--(1.402,1.445)--(1.506,1.190)--(1.610,1.232)--(1.714,1.284)%
  --(1.818,1.343)--(1.922,1.424)--(2.026,1.494)--(2.130,1.524)--(2.234,1.590)--(2.338,1.385)%
  --(2.442,1.417)--(2.546,1.452)--(2.650,1.476)--(2.754,1.511)--(2.858,1.549)--(2.962,1.590)%
  --(3.066,1.628)--(3.170,1.661)--(3.273,1.710)--(3.377,1.751)--(3.481,1.797)--(3.585,1.854)%
  --(3.689,1.891)--(3.793,1.942)--(3.897,2.005)--(4.001,1.900)--(4.105,1.930)--(4.209,1.963)%
  --(4.313,1.964)--(4.417,2.009)--(4.521,2.034)--(4.625,2.070)--(4.729,2.078)--(4.833,2.119)%
  --(4.937,2.159)--(5.041,2.190)--(5.145,2.223)--(5.249,2.241)--(5.353,2.284)--(5.457,2.334)%
  --(5.561,2.392)--(5.665,2.430)--(5.769,2.444);
\gpcolor{color=gp lt color border}
\node[gp node right] at (3.320,3.846) {thick-bounds, $f_1^\text{min}=1/3$};
\gpcolor{rgb color={0.902,0.624,0.000}}
\gpsetdashtype{gp dt 4}
\draw[gp path] (3.430,3.846)--(4.050,3.846);
\draw[gp path] (0.674,0.865)--(0.778,0.908)--(0.882,0.957)--(0.986,1.000)--(1.090,1.049)%
  --(1.194,1.026)--(1.298,1.054)--(1.402,1.088)--(1.506,1.118)--(1.610,1.147)--(1.714,1.189)%
  --(1.818,1.197)--(1.922,1.215)--(2.026,1.237)--(2.130,1.251)--(2.234,1.276)--(2.338,1.314)%
  --(2.442,1.331)--(2.546,1.357)--(2.650,1.374)--(2.754,1.399)--(2.858,1.542)--(2.962,1.562)%
  --(3.066,1.575)--(3.170,1.603)--(3.273,1.598)--(3.377,1.623)--(3.481,1.638)--(3.585,1.656)%
  --(3.689,1.657)--(3.793,1.699)--(3.897,1.710)--(4.001,1.727)--(4.105,1.748)--(4.209,1.768)%
  --(4.313,1.780)--(4.417,1.804)--(4.521,1.830)--(4.625,1.859)--(4.729,1.867)--(4.833,1.899)%
  --(4.937,1.922)--(5.041,2.313)--(5.145,2.331)--(5.249,2.332)--(5.353,2.347)--(5.457,2.366)%
  --(5.561,2.395)--(5.665,2.420)--(5.769,2.416);
\gpcolor{color=gp lt color border}
\gpsetlinetype{gp lt axes}
\gpsetdashtype{gp dt axes}
\draw[gp path] (0.674,0.788)--(0.725,0.788)--(0.777,0.788)--(0.828,0.788)--(0.880,0.788)%
  --(0.931,0.788)--(0.983,0.788)--(1.034,0.788)--(1.086,0.788)--(1.137,0.788)--(1.189,0.788)%
  --(1.240,0.788)--(1.292,0.788)--(1.343,0.788)--(1.394,0.788)--(1.446,0.788)--(1.497,0.788)%
  --(1.549,0.788)--(1.600,0.788)--(1.652,0.788)--(1.703,0.788)--(1.755,0.788)--(1.806,0.788)%
  --(1.858,0.788)--(1.909,0.788)--(1.961,0.788)--(2.012,0.788)--(2.064,0.788)--(2.115,0.788)%
  --(2.166,0.788)--(2.218,0.788)--(2.269,0.788)--(2.321,0.788)--(2.372,0.788)--(2.424,0.788)%
  --(2.475,0.788)--(2.527,0.788)--(2.578,0.788)--(2.630,0.788)--(2.681,0.788)--(2.733,0.788)%
  --(2.784,0.788)--(2.836,0.788)--(2.887,0.788)--(2.938,0.788)--(2.990,0.788)--(3.041,0.788)%
  --(3.093,0.788)--(3.144,0.788)--(3.196,0.788)--(3.247,0.788)--(3.299,0.788)--(3.350,0.788)%
  --(3.402,0.788)--(3.453,0.788)--(3.505,0.788)--(3.556,0.788)--(3.607,0.788)--(3.659,0.788)%
  --(3.710,0.788)--(3.762,0.788)--(3.813,0.788)--(3.865,0.788)--(3.916,0.788)--(3.968,0.788)%
  --(4.019,0.788)--(4.071,0.788)--(4.122,0.788)--(4.174,0.788)--(4.225,0.788)--(4.277,0.788)%
  --(4.328,0.788)--(4.379,0.788)--(4.431,0.788)--(4.482,0.788)--(4.534,0.788)--(4.585,0.788)%
  --(4.637,0.788)--(4.688,0.788)--(4.740,0.788)--(4.791,0.788)--(4.843,0.788)--(4.894,0.788)%
  --(4.946,0.788)--(4.997,0.788)--(5.048,0.788)--(5.100,0.788)--(5.151,0.788)--(5.203,0.788)%
  --(5.254,0.788)--(5.306,0.788)--(5.357,0.788)--(5.409,0.788)--(5.460,0.788)--(5.512,0.788)%
  --(5.563,0.788)--(5.615,0.788)--(5.666,0.788)--(5.718,0.788)--(5.769,0.788);
\gpsetlinetype{gp lt border}
\gpsetdashtype{gp dt solid}
\gpsetlinewidth{1.00}
\draw[gp path] (0.570,4.814)--(0.570,0.592)--(5.769,0.592)--(5.769,4.814)--cycle;
\gpdefrectangularnode{gp plot 1}{\pgfpoint{0.570cm}{0.592cm}}{\pgfpoint{5.769cm}{4.814cm}}
\end{tikzpicture}
\subcaption{Algorithmic variants of TSQR TT-SVD for a $2^{30}$ tensor with different prescribed reduction factors $f_1^\text{min}$.}\label{fig:TT-SVD_variants}
\end{subfigure}
\caption{Single socket run-time of different TT-SVD algorithms for varying maximal TT-rank.}
\end{figure}

In the following, we consider the complete TT-SVD algorithm and different variants and implementations of it.
\autoref{fig:TT-SVD_other_libraries} illustrates the run-time of the TT-SVD algorithm in different software libraries.
All cases show the run-time for decomposing a random double precision $2^{27}$ tensor on a single CPU socket with 14 cores with a prescribed maximal TT-rank.
For several of these libraries, we tested different variants and LAPACK back-ends~\cite{LAPACK,IntelMKL}.
Here, we only report the timings for the fastest variant that we could find.
We show results for the following libraries:
\begin{description}
\item[TSQR TT-SVD] The implementation discussed in this paper.
\item[ttpy] A library written in Fortran and Python by the author of~\cite{Oseledets2011}.
\item[t3f] A library based on the TensorFlow framework~\cite{Novikov2020}.
\item[TensorToolbox] A Python library from the author of~\cite{Bigoni2016}.
\item[tntorch] A library based on PyTorch~\cite{tnTorch}.
\item[TT-SVD with numpy] Simple implementation in \texttt{numpy}~\cite{Harris2020} inspired by~\cite{Gelss2019}.
\end{description}
Both \texttt{ttpy} and \texttt{TensorToolbox} use the older (and in many cases slower) \verb|dgesvd| routine for calculating SVD decompositions.
Our classical TT-SVD implementation with \texttt{numpy} uses the newer LAPACK routine \verb|dgesdd|.
The \texttt{ttpy} library is still faster in many cases.
The \texttt{t3f} library is based on TensorFlow, which is optimized for GPUs. It uses the C++ library Eigen~\cite{Eigen} as back-end on CPUs.
However, only some routines in Eigen are parallelized for multi-core CPUs which explains why \texttt{t3f} is slow here.
In contrast to all other variants, the run-time of the TT decomposition in \texttt{tntorch} is almost independent of the maximal TT-rank.
\texttt{tntorch} does not implement the TT-SVD algorithm, but instead first constructs a tensor-train of maximal rank, followed by a left-right orthogonalization step and TT rounding.
The computationally costly part is the left-right orthogonalization step, which is based on QR decompositions whose size only depend on the size of the input tensor and not on the desired rank.

Our TSQR TT-SVD implementation is significantly faster than all other implementations for two reasons.
First, there are multiple combinations of basic linear algebra building blocks that calculate the desired result.
This is an example of the linear algebra mapping problem as discussed in~\cite{Psarras2019}.
Here, we choose a combination of building blocks (Q-less TSQR + multiplication with truncated right singular vectors)  that leads to (almost) minimal data transfers.
Second, common linear algebra software and algorithms are not optimized for avoiding data transfers.
However, for the tall-skinny matrix operations required here, the data transfers determine the performance.
For a detailed overview on communication avoiding linear algebra algorithms, see e.g.\ \cite{Ballard2014} and the references therein.
An interesting discussion that distinguishes between the effects of reading and modifying data can be found in~\cite{Carson2016}.

\autoref{fig:TT-SVD_variants} shows the run-time for the different variants of the TSQR TT-SVD algorithm discussed in \autoref{sec:TSQR_TT-SVD_variants}.
This is the worst case run-time of the algorithm because we approximate a random input matrix and we only prescribe the maximal TT-rank.
For the plain case ($f_1^\text{min}=1$), there is no reduction in the data size in the first steps for $r_\text{max} > 1$.
For the thick-bounds and two-sided variants we set $m_\text{min}=16$ (see \autoref{alg:Two-Sided_TSQR_TT_SVD}).
This reduces the run-time for small TT-ranks (difference between plain and other variants for $r_\text{max}=1$).
See \autoref{tab:tt_svd_variants_dimensions} for some examples on resulting dimensions and TT-ranks.
\begin{table}
\centering
\begin{tabular}{lcccc}
case & $r_\text{max}$ & effective dim. & TT-ranks & reduction factors\\
       &                        &       ($n_i$)    &   ($r_i$)  &      ($f_j$)\\
\midrule
plain                     & 1 & \dots, 2, 2 & \dots, 1, 1 & $\frac{1}{2}$, $\frac{1}{2}$, \dots \\
($f_1^\text{min}=1$)      & 5 & \dots, 2, 2 & \dots, 5, 5, 4, 2 & 1, 1, $\frac{5}{8}$, $\frac{1}{2}$, \dots\\
                           & 16 & \dots, 2, 2 & \dots, 16, 16, 8, 4, 2 & 1, 1, 1, 1, $\frac{1}{2}$, \dots\\
\midrule
thick-bounds             & 1 & \dots, 2, 2, 16 & \dots, 1, 1 & $\frac{1}{16}$, $\frac{1}{2}$, $\frac{1}{2}$, \dots\\
($f_1^\text{min}=1/2$)   & 5 & \dots, 2, 2, 16 & \dots, 5, 5 & $\frac{5}{16}$, $\frac{1}{2}$, $\frac{1}{2}$, \dots\\
                          & 16 & \dots, 2, 2, 32 & \dots, 16, 16 & $\frac{1}{2}$, $\frac{1}{2}$, \dots\\
\end{tabular}
\caption{Examples for the resulting effective dimensions and TT-ranks for the different TT-SVD variants for a $2^d$ tensor.
We consider the right-most dimensions and ranks as our implementation calculates the decomposition from right to left.}\label{tab:tt_svd_variants_dimensions}
\end{table}
As expected, the plain variant is slower as it needs to transfer more data in the first iterations.
For all cases with a prescribed reduction $f_1^\text{min}<1$, we observe roughly a linear scaling with the maximal TT-rank as predicted by the performance analysis for the compute-bound case.
And for small ranks, the run-time is of the order of copying the data in memory.
For our implementation the choice $f_1^\text{min}=1/2$ appears to be optimal even though the theoretical analysis indicates that a smaller $f_1^\text{min}$ could be beneficial.
Decreasing $f_1^\text{min}$ increases the number of columns of the matrices in the first step.
This leads to more work and the obtained bandwidth of our TSQR implementation decreases (see \autoref{fig:TSQR_Roofline_bandwidth}).
The two-sided variant uses thick-bounds as well but it is always slower with our implementation.

\begin{figure}
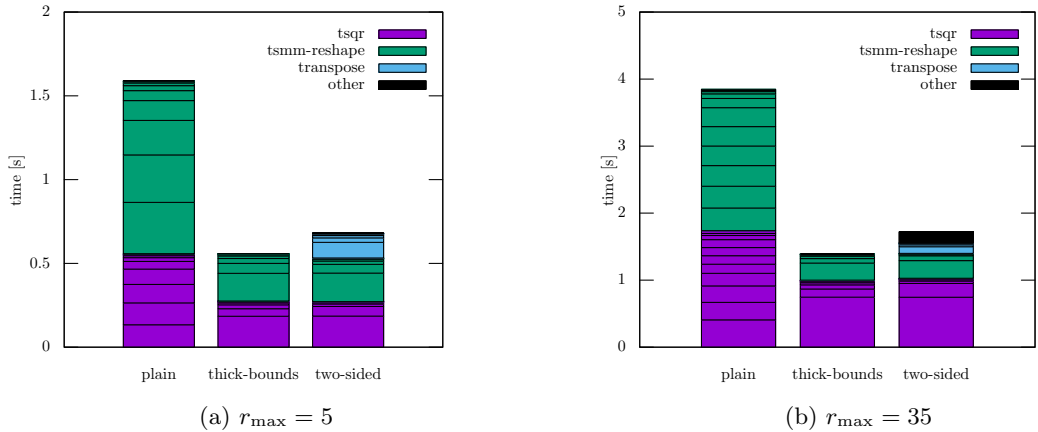

\begin{subfigure}[t]{0.47\textwidth}
\input{fig/tt_from_dense_single_socket_per_function_rank5.tex}
\subcaption{$r_\text{max}=5$}
\end{subfigure}
\hfill
\begin{subfigure}[t]{0.47\textwidth}
\input{fig/tt_from_dense_single_socket_per_function_rank35.tex}
\subcaption{$r_\text{max}=35$}
\end{subfigure}
\caption{Timings for the building blocks in different TT-SVD variants for two cases from \autoref{fig:TT-SVD_variants}.
The \emph{transpose} timings refer to the transpose/reshape operations in lines 9 and 14 of \autoref{alg:Two-Sided_TSQR_TT_SVD}.}\label{fig:TT-SVD_timings_per_function}
\end{figure}
The run-time of the individual steps of the algorithm are illustrated in \autoref{fig:TT-SVD_timings_per_function}.
We clearly see the effect of combining multiple dimensions:
The first TSQR step takes longer but all subsequent steps are faster.
The two-sided variant is only slower due to the additional transpose operation required in our implementation.
For real-world problems, the two-sided variant might be faster depending on the resulting TT-ranks.

\begin{table}
\centering
\begin{tabular}{lccc}
case& $r_\text{max}$ & operations & data transfers\\
& & [GFlop] & [GByte]\\
\midrule
plain (estimate with $\bar f = 1/2$) & 1 & 14 (13) & 43 (43)\\
thick-bounds (estimate with $\bar f = 1/16$) & 1 & 41 (39) & 21 (19)\\
thick-bounds (estimate with $\bar f = 1/2$) & 31 & 417 (399) & 43 (43)\\
\end{tabular}
%
%
%
%
%
\caption{Measured and estimated number of floating point operations and transferred data volume between the CPU and the main memory for the TSQR TT-SVD algorithm with a $2^{30}$ tensor in double precision.
The measurements were performed with likwid-perfctr~\cite{Treibig2010}.
Estimates based on  \eqref{eq:TT-SVD_data_volume} and \eqref{eq:TT-SVD_operations} are shown in parentheses.}\label{tab:TT-SVD_likwid_measurements}
\end{table}
To validate our assumptions in the performance analysis in \autoref{sec:TT-SVD_perf_analysis}, we measured data transfers and flops for several cases using CPU performance counters, see \autoref{tab:TT-SVD_likwid_measurements}.
We compare cases where the simple estimates with the global reduction factor $\bar f$ fit well and we observe a good correlation with the measurements.
Depending on the dimensions and the desired maximal rank, the reduction in the first step can differ from the following steps (see \autoref{tab:tt_svd_variants_dimensions}) which is not captured by \eqref{eq:TT-SVD_data_volume}, \eqref{eq:TT-SVD_operations}.

\begin{figure}
\begin{subfigure}[t]{0.47\textwidth}
\begin{tikzpicture}[gnuplot]
\tikzset{every node/.append style={scale=0.60}}
\path (0.000,0.000) rectangle (6.100,5.000);
\gpcolor{color=gp lt color border}
\gpsetlinetype{gp lt border}
\gpsetdashtype{gp dt solid}
\gpsetlinewidth{1.00}
\draw[gp path] (0.680,0.592)--(0.860,0.592);
\draw[gp path] (5.769,0.592)--(5.589,0.592);
\node[gp node right] at (0.570,0.592) {$0$};
\draw[gp path] (0.680,1.296)--(0.860,1.296);
\draw[gp path] (5.769,1.296)--(5.589,1.296);
\node[gp node right] at (0.570,1.296) {$5$};
\draw[gp path] (0.680,1.999)--(0.860,1.999);
\draw[gp path] (5.769,1.999)--(5.589,1.999);
\node[gp node right] at (0.570,1.999) {$10$};
\draw[gp path] (0.680,2.703)--(0.860,2.703);
\draw[gp path] (5.769,2.703)--(5.589,2.703);
\node[gp node right] at (0.570,2.703) {$15$};
\draw[gp path] (0.680,3.407)--(0.860,3.407);
\draw[gp path] (5.769,3.407)--(5.589,3.407);
\node[gp node right] at (0.570,3.407) {$20$};
\draw[gp path] (0.680,4.110)--(0.860,4.110);
\draw[gp path] (5.769,4.110)--(5.589,4.110);
\node[gp node right] at (0.570,4.110) {$25$};
\draw[gp path] (0.680,4.814)--(0.860,4.814);
\draw[gp path] (5.769,4.814)--(5.589,4.814);
\node[gp node right] at (0.570,4.814) {$30$};
\draw[gp path] (0.680,0.592)--(0.680,0.772);
\draw[gp path] (0.680,4.814)--(0.680,4.634);
\node[gp node center] at (0.680,0.407) {$0$};
\draw[gp path] (1.189,0.592)--(1.189,0.772);
\draw[gp path] (1.189,4.814)--(1.189,4.634);
\node[gp node center] at (1.189,0.407) {$5$};
\draw[gp path] (1.698,0.592)--(1.698,0.772);
\draw[gp path] (1.698,4.814)--(1.698,4.634);
\node[gp node center] at (1.698,0.407) {$10$};
\draw[gp path] (2.207,0.592)--(2.207,0.772);
\draw[gp path] (2.207,4.814)--(2.207,4.634);
\node[gp node center] at (2.207,0.407) {$15$};
\draw[gp path] (2.716,0.592)--(2.716,0.772);
\draw[gp path] (2.716,4.814)--(2.716,4.634);
\node[gp node center] at (2.716,0.407) {$20$};
\draw[gp path] (3.225,0.592)--(3.225,0.772);
\draw[gp path] (3.225,4.814)--(3.225,4.634);
\node[gp node center] at (3.225,0.407) {$25$};
\draw[gp path] (3.733,0.592)--(3.733,0.772);
\draw[gp path] (3.733,4.814)--(3.733,4.634);
\node[gp node center] at (3.733,0.407) {$30$};
\draw[gp path] (4.242,0.592)--(4.242,0.772);
\draw[gp path] (4.242,4.814)--(4.242,4.634);
\node[gp node center] at (4.242,0.407) {$35$};
\draw[gp path] (4.751,0.592)--(4.751,0.772);
\draw[gp path] (4.751,4.814)--(4.751,4.634);
\node[gp node center] at (4.751,0.407) {$40$};
\draw[gp path] (5.260,0.592)--(5.260,0.772);
\draw[gp path] (5.260,4.814)--(5.260,4.634);
\node[gp node center] at (5.260,0.407) {$45$};
\draw[gp path] (5.769,0.592)--(5.769,0.772);
\draw[gp path] (5.769,4.814)--(5.769,4.634);
\node[gp node center] at (5.769,0.407) {$50$};
\draw[gp path] (0.680,4.814)--(0.680,0.592)--(5.769,0.592)--(5.769,4.814)--cycle;
\node[gp node center,rotate=-270] at (0.175,2.703) {time [s]};
\node[gp node center] at (3.224,0.130) {max. rank};
\node[gp node right] at (1.890,4.521) {$2^{30}$ tensor};
\gpcolor{rgb color={0.580,0.000,0.827}}
\gpsetdashtype{gp dt 1}
\gpsetlinewidth{2.00}
\draw[gp path] (2.000,4.521)--(2.620,4.521);
\draw[gp path] (0.782,0.641)--(0.884,0.648)--(0.985,0.656)--(1.087,0.665)--(1.189,0.674)%
  --(1.291,0.682)--(1.392,0.692)--(1.494,0.680)--(1.596,0.686)--(1.698,0.691)--(1.800,0.698)%
  --(1.901,0.700)--(2.003,0.706)--(2.105,0.712)--(2.207,0.720)--(2.308,0.713)--(2.410,0.720)%
  --(2.512,0.723)--(2.614,0.727)--(2.716,0.730)--(2.817,0.734)--(2.919,0.739)--(3.021,0.746)%
  --(3.123,0.749)--(3.225,0.752)--(3.326,0.757)--(3.428,0.762)--(3.530,0.768)--(3.632,0.774)%
  --(3.733,0.776)--(3.835,0.782)--(3.937,0.784)--(4.039,0.789)--(4.141,0.791)--(4.242,0.795)%
  --(4.344,0.795)--(4.446,0.799)--(4.548,0.803)--(4.649,0.808)--(4.751,0.809)--(4.853,0.814)%
  --(4.955,0.818)--(5.057,0.822)--(5.158,0.827)--(5.260,0.828)--(5.362,0.833)--(5.464,0.838)%
  --(5.565,0.844)--(5.667,0.850)--(5.769,0.852);
\gpcolor{color=gp lt color border}
\node[gp node right] at (1.890,4.296) {$4^{15}$ tensor};
\gpcolor{rgb color={0.000,0.620,0.451}}
\gpsetdashtype{gp dt 2}
\draw[gp path] (2.000,4.296)--(2.620,4.296);
\draw[gp path] (0.782,0.639)--(0.884,0.645)--(0.985,0.650)--(1.087,0.658)--(1.189,0.665)%
  --(1.291,0.672)--(1.392,0.680)--(1.494,0.685)--(1.596,0.689)--(1.698,0.690)--(1.800,0.694)%
  --(1.901,0.697)--(2.003,0.701)--(2.105,0.705)--(2.207,0.706)--(2.308,0.711)--(2.410,0.717)%
  --(2.512,0.720)--(2.614,0.724)--(2.716,0.728)--(2.817,0.735)--(2.919,0.737)--(3.021,0.742)%
  --(3.123,0.745)--(3.225,0.755)--(3.326,0.759)--(3.428,0.764)--(3.530,0.770)--(3.632,0.777)%
  --(3.733,0.782)--(3.835,0.789)--(3.937,0.863)--(4.039,0.864)--(4.141,0.870)--(4.242,0.868)%
  --(4.344,0.871)--(4.446,0.874)--(4.548,0.879)--(4.649,0.880)--(4.751,0.882)--(4.853,0.887)%
  --(4.955,0.890)--(5.057,0.895)--(5.158,0.897)--(5.260,0.897)--(5.362,0.899)--(5.464,0.905)%
  --(5.565,0.910)--(5.667,0.916)--(5.769,0.912);
\gpcolor{color=gp lt color border}
\node[gp node right] at (1.890,4.071) {$8^{10}$ tensor};
\gpcolor{rgb color={0.337,0.706,0.914}}
\gpsetdashtype{gp dt 3}
\draw[gp path] (2.000,4.071)--(2.620,4.071);
\draw[gp path] (0.782,0.634)--(0.884,0.646)--(0.985,0.658)--(1.087,0.672)--(1.189,0.678)%
  --(1.291,0.682)--(1.392,0.685)--(1.494,0.688)--(1.596,0.692)--(1.698,0.695)--(1.800,0.701)%
  --(1.901,0.703)--(2.003,0.712)--(2.105,0.714)--(2.207,0.719)--(2.308,0.725)--(2.410,0.735)%
  --(2.512,0.739)--(2.614,0.747)--(2.716,0.749)--(2.817,0.761)--(2.919,0.764)--(3.021,0.773)%
  --(3.123,0.781)--(3.225,0.787)--(3.326,0.797)--(3.428,0.807)--(3.530,0.818)--(3.632,0.827)%
  --(3.733,0.835)--(3.835,0.846)--(3.937,0.920)--(4.039,0.933)--(4.141,0.942)--(4.242,0.954)%
  --(4.344,0.967)--(4.446,0.978)--(4.548,0.994)--(4.649,1.006)--(4.751,1.018)--(4.853,1.033)%
  --(4.955,1.047)--(5.057,1.064)--(5.158,1.081)--(5.260,1.094)--(5.362,1.110)--(5.464,1.130)%
  --(5.565,1.151)--(5.667,1.175)--(5.769,1.187);
\gpcolor{color=gp lt color border}
\node[gp node right] at (1.890,3.846) {$10^9$ tensor};
\gpcolor{rgb color={0.902,0.624,0.000}}
\gpsetdashtype{gp dt 4}
\draw[gp path] (2.000,3.846)--(2.620,3.846);
\draw[gp path] (0.782,0.632)--(0.884,0.640)--(0.985,0.650)--(1.087,0.664)--(1.189,0.676)%
  --(1.291,0.691)--(1.392,0.700)--(1.494,0.702)--(1.596,0.718)--(1.698,0.707)--(1.800,0.712)%
  --(1.901,0.716)--(2.003,0.734)--(2.105,0.730)--(2.207,0.727)--(2.308,0.732)--(2.410,0.749)%
  --(2.512,0.747)--(2.614,0.750)--(2.716,0.751)--(2.817,0.772)--(2.919,0.768)--(3.021,0.773)%
  --(3.123,0.777)--(3.225,0.795)--(3.326,0.796)--(3.428,0.802)--(3.530,0.808)--(3.632,0.819)%
  --(3.733,0.827)--(3.835,0.837)--(3.937,0.846)--(4.039,0.856)--(4.141,0.869)--(4.242,0.875)%
  --(4.344,0.886)--(4.446,0.900)--(4.548,0.912)--(4.649,0.928)--(4.751,0.936)--(4.853,0.950)%
  --(4.955,0.967)--(5.057,0.980)--(5.158,0.998)--(5.260,1.008)--(5.362,1.024)--(5.464,1.042)%
  --(5.565,1.059)--(5.667,1.077)--(5.769,1.156);
\gpcolor{color=gp lt color border}
\node[gp node right] at (1.890,3.621) {$32^6$ tensor};
\gpcolor{rgb color={0.941,0.894,0.259}}
\gpsetdashtype{gp dt 5}
\draw[gp path] (2.000,3.621)--(2.620,3.621);
\draw[gp path] (0.782,0.650)--(0.884,0.657)--(0.985,0.664)--(1.087,0.678)--(1.189,0.685)%
  --(1.291,0.699)--(1.392,0.720)--(1.494,0.741)--(1.596,0.768)--(1.698,0.785)--(1.800,0.812)%
  --(1.901,0.852)--(2.003,0.887)--(2.105,0.930)--(2.207,0.965)--(2.308,1.017)--(2.410,1.079)%
  --(2.512,1.150)--(2.614,1.202)--(2.716,1.262)--(2.817,1.328)--(2.919,1.497)--(3.021,1.594)%
  --(3.123,1.723)--(3.225,1.815)--(3.326,1.924)--(3.428,2.084)--(3.530,2.413)--(3.632,2.348)%
  --(3.733,2.464)--(3.835,2.614)--(3.937,2.767)--(4.039,2.842)--(4.141,2.900)--(4.242,2.971)%
  --(4.344,3.038)--(4.446,3.094)--(4.548,3.152)--(4.649,3.259)--(4.751,3.361)--(4.853,3.439)%
  --(4.955,3.501)--(5.057,3.610)--(5.158,3.721)--(5.260,3.841)--(5.362,3.920)--(5.464,4.031)%
  --(5.565,4.244)--(5.667,4.358)--(5.769,4.505);
\gpcolor{color=gp lt color border}
\gpsetdashtype{gp dt solid}
\gpsetlinewidth{1.00}
\draw[gp path] (0.680,4.814)--(0.680,0.592)--(5.769,0.592)--(5.769,4.814)--cycle;
\gpdefrectangularnode{gp plot 1}{\pgfpoint{0.680cm}{0.592cm}}{\pgfpoint{5.769cm}{4.814cm}}
\end{tikzpicture}
\subcaption{All cases: $2^{30}$ to $32^6$.}
\end{subfigure}
\hfill
\begin{subfigure}[t]{0.47\textwidth}
\begin{tikzpicture}[gnuplot]
\tikzset{every node/.append style={scale=0.60}}
\path (0.000,0.000) rectangle (6.100,5.000);
\gpcolor{color=gp lt color border}
\gpsetlinetype{gp lt border}
\gpsetdashtype{gp dt solid}
\gpsetlinewidth{1.00}
\draw[gp path] (0.790,0.592)--(0.970,0.592);
\draw[gp path] (5.769,0.592)--(5.589,0.592);
\node[gp node right] at (0.680,0.592) {$0$};
\draw[gp path] (0.790,1.120)--(0.970,1.120);
\draw[gp path] (5.769,1.120)--(5.589,1.120);
\node[gp node right] at (0.680,1.120) {$0.5$};
\draw[gp path] (0.790,1.648)--(0.970,1.648);
\draw[gp path] (5.769,1.648)--(5.589,1.648);
\node[gp node right] at (0.680,1.648) {$1$};
\draw[gp path] (0.790,2.175)--(0.970,2.175);
\draw[gp path] (5.769,2.175)--(5.589,2.175);
\node[gp node right] at (0.680,2.175) {$1.5$};
\draw[gp path] (0.790,2.703)--(0.970,2.703);
\draw[gp path] (5.769,2.703)--(5.589,2.703);
\node[gp node right] at (0.680,2.703) {$2$};
\draw[gp path] (0.790,3.231)--(0.970,3.231);
\draw[gp path] (5.769,3.231)--(5.589,3.231);
\node[gp node right] at (0.680,3.231) {$2.5$};
\draw[gp path] (0.790,3.759)--(0.970,3.759);
\draw[gp path] (5.769,3.759)--(5.589,3.759);
\node[gp node right] at (0.680,3.759) {$3$};
\draw[gp path] (0.790,4.286)--(0.970,4.286);
\draw[gp path] (5.769,4.286)--(5.589,4.286);
\node[gp node right] at (0.680,4.286) {$3.5$};
\draw[gp path] (0.790,4.814)--(0.970,4.814);
\draw[gp path] (5.769,4.814)--(5.589,4.814);
\node[gp node right] at (0.680,4.814) {$4$};
\draw[gp path] (0.790,0.592)--(0.790,0.772);
\draw[gp path] (0.790,4.814)--(0.790,4.634);
\node[gp node center] at (0.790,0.407) {$0$};
\draw[gp path] (1.288,0.592)--(1.288,0.772);
\draw[gp path] (1.288,4.814)--(1.288,4.634);
\node[gp node center] at (1.288,0.407) {$5$};
\draw[gp path] (1.786,0.592)--(1.786,0.772);
\draw[gp path] (1.786,4.814)--(1.786,4.634);
\node[gp node center] at (1.786,0.407) {$10$};
\draw[gp path] (2.284,0.592)--(2.284,0.772);
\draw[gp path] (2.284,4.814)--(2.284,4.634);
\node[gp node center] at (2.284,0.407) {$15$};
\draw[gp path] (2.782,0.592)--(2.782,0.772);
\draw[gp path] (2.782,4.814)--(2.782,4.634);
\node[gp node center] at (2.782,0.407) {$20$};
\draw[gp path] (3.280,0.592)--(3.280,0.772);
\draw[gp path] (3.280,4.814)--(3.280,4.634);
\node[gp node center] at (3.280,0.407) {$25$};
\draw[gp path] (3.777,0.592)--(3.777,0.772);
\draw[gp path] (3.777,4.814)--(3.777,4.634);
\node[gp node center] at (3.777,0.407) {$30$};
\draw[gp path] (4.275,0.592)--(4.275,0.772);
\draw[gp path] (4.275,4.814)--(4.275,4.634);
\node[gp node center] at (4.275,0.407) {$35$};
\draw[gp path] (4.773,0.592)--(4.773,0.772);
\draw[gp path] (4.773,4.814)--(4.773,4.634);
\node[gp node center] at (4.773,0.407) {$40$};
\draw[gp path] (5.271,0.592)--(5.271,0.772);
\draw[gp path] (5.271,4.814)--(5.271,4.634);
\node[gp node center] at (5.271,0.407) {$45$};
\draw[gp path] (5.769,0.592)--(5.769,0.772);
\draw[gp path] (5.769,4.814)--(5.769,4.634);
\node[gp node center] at (5.769,0.407) {$50$};
\draw[gp path] (0.790,4.814)--(0.790,0.592)--(5.769,0.592)--(5.769,4.814)--cycle;
\node[gp node left] at (4.275,1.014) {copy data once};
\node[gp node center,rotate=-270] at (0.175,2.703) {time [s]};
\node[gp node center] at (3.279,0.130) {max. rank};
\gpcolor{rgb color={0.580,0.000,0.827}}
\gpsetdashtype{gp dt 1}
\gpsetlinewidth{2.00}
\draw[gp path] (0.890,0.957)--(0.989,1.010)--(1.089,1.071)--(1.188,1.136)--(1.288,1.204)%
  --(1.387,1.270)--(1.487,1.341)--(1.587,1.253)--(1.686,1.294)--(1.786,1.333)--(1.885,1.384)%
  --(1.985,1.404)--(2.085,1.446)--(2.184,1.493)--(2.284,1.553)--(2.383,1.501)--(2.483,1.550)%
  --(2.582,1.576)--(2.682,1.605)--(2.782,1.625)--(2.881,1.658)--(2.981,1.692)--(3.080,1.745)%
  --(3.180,1.772)--(3.280,1.794)--(3.379,1.831)--(3.479,1.869)--(3.578,1.914)--(3.678,1.956)%
  --(3.777,1.971)--(3.877,2.019)--(3.977,2.035)--(4.076,2.067)--(4.176,2.087)--(4.275,2.116)%
  --(4.375,2.115)--(4.474,2.144)--(4.574,2.172)--(4.674,2.213)--(4.773,2.220)--(4.873,2.259)%
  --(4.972,2.286)--(5.072,2.316)--(5.172,2.352)--(5.271,2.359)--(5.371,2.399)--(5.470,2.434)%
  --(5.570,2.480)--(5.669,2.526)--(5.769,2.543);
\gpcolor{rgb color={0.000,0.620,0.451}}
\gpsetdashtype{gp dt 2}
\draw[gp path] (0.890,0.941)--(0.989,0.988)--(1.089,1.030)--(1.188,1.084)--(1.288,1.139)%
  --(1.387,1.190)--(1.487,1.250)--(1.587,1.289)--(1.686,1.321)--(1.786,1.330)--(1.885,1.360)%
  --(1.985,1.381)--(2.085,1.406)--(2.184,1.441)--(2.284,1.451)--(2.383,1.485)--(2.483,1.527)%
  --(2.582,1.554)--(2.682,1.582)--(2.782,1.612)--(2.881,1.664)--(2.981,1.679)--(3.080,1.713)%
  --(3.180,1.743)--(3.280,1.811)--(3.379,1.842)--(3.479,1.883)--(3.578,1.927)--(3.678,1.982)%
  --(3.777,2.016)--(3.877,2.070)--(3.977,2.623)--(4.076,2.632)--(4.176,2.677)--(4.275,2.666)%
  --(4.375,2.683)--(4.474,2.706)--(4.574,2.747)--(4.674,2.753)--(4.773,2.766)--(4.873,2.807)%
  --(4.972,2.830)--(5.072,2.863)--(5.172,2.881)--(5.271,2.882)--(5.371,2.891)--(5.470,2.941)%
  --(5.570,2.980)--(5.669,3.024)--(5.769,2.995);
\gpcolor{rgb color={0.337,0.706,0.914}}
\gpsetdashtype{gp dt 3}
\draw[gp path] (0.890,0.910)--(0.989,0.995)--(1.089,1.085)--(1.188,1.195)--(1.288,1.235)%
  --(1.387,1.263)--(1.487,1.287)--(1.587,1.309)--(1.686,1.344)--(1.786,1.362)--(1.885,1.406)%
  --(1.985,1.427)--(2.085,1.495)--(2.184,1.510)--(2.284,1.548)--(2.383,1.591)--(2.483,1.668)%
  --(2.582,1.695)--(2.682,1.756)--(2.782,1.773)--(2.881,1.857)--(2.981,1.883)--(3.080,1.953)%
  --(3.180,2.009)--(3.280,2.052)--(3.379,2.131)--(3.479,2.204)--(3.578,2.288)--(3.678,2.354)%
  --(3.777,2.413)--(3.877,2.500)--(3.977,3.054)--(4.076,3.151)--(4.176,3.216)--(4.275,3.309)%
  --(4.375,3.403)--(4.474,3.487)--(4.574,3.604)--(4.674,3.700)--(4.773,3.788)--(4.873,3.899)%
  --(4.972,4.003)--(5.072,4.136)--(5.172,4.258)--(5.271,4.358)--(5.371,4.480)--(5.470,4.629)%
  --(5.570,4.787)--(5.585,4.814);
\gpcolor{rgb color={0.902,0.624,0.000}}
\gpsetdashtype{gp dt 4}
\draw[gp path] (0.890,0.889)--(0.989,0.953)--(1.089,1.029)--(1.188,1.131)--(1.288,1.222)%
  --(1.387,1.332)--(1.487,1.403)--(1.587,1.419)--(1.686,1.539)--(1.786,1.456)--(1.885,1.494)%
  --(1.985,1.520)--(2.085,1.653)--(2.184,1.623)--(2.284,1.606)--(2.383,1.642)--(2.483,1.771)%
  --(2.582,1.756)--(2.682,1.775)--(2.782,1.786)--(2.881,1.939)--(2.981,1.915)--(3.080,1.952)%
  --(3.180,1.977)--(3.280,2.118)--(3.379,2.124)--(3.479,2.165)--(3.578,2.215)--(3.678,2.291)%
  --(3.777,2.356)--(3.877,2.427)--(3.977,2.495)--(4.076,2.575)--(4.176,2.671)--(4.275,2.715)%
  --(4.375,2.799)--(4.474,2.900)--(4.574,2.994)--(4.674,3.109)--(4.773,3.172)--(4.873,3.280)%
  --(4.972,3.403)--(5.072,3.499)--(5.172,3.634)--(5.271,3.709)--(5.371,3.834)--(5.470,3.966)%
  --(5.570,4.091)--(5.669,4.229)--(5.767,4.814);
\gpcolor{rgb color={0.941,0.894,0.259}}
\gpsetdashtype{gp dt 5}
\draw[gp path] (0.890,1.026)--(0.989,1.076)--(1.089,1.132)--(1.188,1.235)--(1.288,1.293)%
  --(1.387,1.397)--(1.487,1.551)--(1.587,1.711)--(1.686,1.910)--(1.786,2.037)--(1.885,2.244)%
  --(1.985,2.541)--(2.085,2.808)--(2.184,3.128)--(2.284,3.388)--(2.383,3.779)--(2.483,4.242)%
  --(2.582,4.775)--(2.591,4.814);
\gpcolor{color=gp lt color border}
\gpsetlinetype{gp lt axes}
\gpsetdashtype{gp dt axes}
\draw[gp path] (0.890,0.851)--(0.939,0.851)--(0.988,0.851)--(1.037,0.851)--(1.087,0.851)%
  --(1.136,0.851)--(1.185,0.851)--(1.235,0.851)--(1.284,0.851)--(1.333,0.851)--(1.382,0.851)%
  --(1.432,0.851)--(1.481,0.851)--(1.530,0.851)--(1.580,0.851)--(1.629,0.851)--(1.678,0.851)%
  --(1.727,0.851)--(1.777,0.851)--(1.826,0.851)--(1.875,0.851)--(1.925,0.851)--(1.974,0.851)%
  --(2.023,0.851)--(2.072,0.851)--(2.122,0.851)--(2.171,0.851)--(2.220,0.851)--(2.270,0.851)%
  --(2.319,0.851)--(2.368,0.851)--(2.417,0.851)--(2.467,0.851)--(2.516,0.851)--(2.565,0.851)%
  --(2.615,0.851)--(2.664,0.851)--(2.713,0.851)--(2.762,0.851)--(2.812,0.851)--(2.861,0.851)%
  --(2.910,0.851)--(2.960,0.851)--(3.009,0.851)--(3.058,0.851)--(3.107,0.851)--(3.157,0.851)%
  --(3.206,0.851)--(3.255,0.851)--(3.305,0.851)--(3.354,0.851)--(3.403,0.851)--(3.453,0.851)%
  --(3.502,0.851)--(3.551,0.851)--(3.600,0.851)--(3.650,0.851)--(3.699,0.851)--(3.748,0.851)%
  --(3.798,0.851)--(3.847,0.851)--(3.896,0.851)--(3.945,0.851)--(3.995,0.851)--(4.044,0.851)%
  --(4.093,0.851)--(4.143,0.851)--(4.192,0.851)--(4.241,0.851)--(4.290,0.851)--(4.340,0.851)%
  --(4.389,0.851)--(4.438,0.851)--(4.488,0.851)--(4.537,0.851)--(4.586,0.851)--(4.635,0.851)%
  --(4.685,0.851)--(4.734,0.851)--(4.783,0.851)--(4.833,0.851)--(4.882,0.851)--(4.931,0.851)%
  --(4.980,0.851)--(5.030,0.851)--(5.079,0.851)--(5.128,0.851)--(5.178,0.851)--(5.227,0.851)%
  --(5.276,0.851)--(5.325,0.851)--(5.375,0.851)--(5.424,0.851)--(5.473,0.851)--(5.523,0.851)%
  --(5.572,0.851)--(5.621,0.851)--(5.670,0.851)--(5.720,0.851)--(5.769,0.851);
\gpsetlinetype{gp lt border}
\gpsetdashtype{gp dt solid}
\gpsetlinewidth{1.00}
\draw[gp path] (0.790,4.814)--(0.790,0.592)--(5.769,0.592)--(5.769,4.814)--cycle;
\gpdefrectangularnode{gp plot 1}{\pgfpoint{0.790cm}{0.592cm}}{\pgfpoint{5.769cm}{4.814cm}}
\end{tikzpicture}
\subcaption{Zooming in on faster cases.}
\end{subfigure}
\caption{Timings for TSQR TT-SVD for varying dimensions on a single socket. Uses the thick-bounds variant with $f_1^\text{min}=\frac{1}{2}$ where beneficial.}\label{fig:TT-SVD_varying_n}
\end{figure}
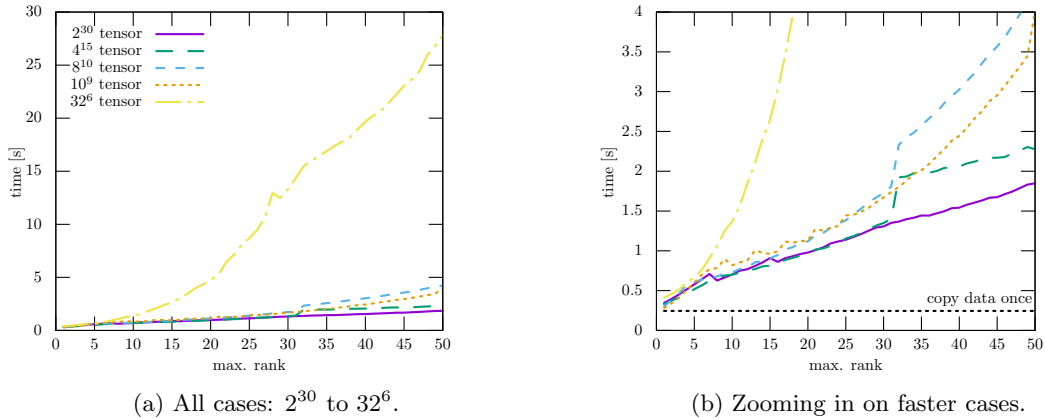
The experiments above use $2^d$ tensors for simplicity (as in the QTT format).
If we increase the size of the individual dimensions, the compute intensity of the TSQR TT-SVD algorithm increases.
\autoref{fig:TT-SVD_varying_n} visualizes the run-time for decomposing tensors of different dimensions with approximately the same total size.
For very small maximal rank ($r_\text{max} < 5$), all cases require similar run-time.
For higher maximal ranks, the cases with a higher individual dimension become more costly.
Near $r_\text{max}= 32 $ there are some interesting deviations in the run-time from the expected linear growth.
We can explain these deviations by the possible choices for combining dimensions in the thick-bounds algorithm: depending on $r_\text{max}$ and $n_i$ there are only few discrete choices for the number of columns $m$ of the first step.
In particular, we obtain $m = 100$ for $r_\text{max}=10, \dots 49$ for the $10^9$ tensor but $m = 512$ for $r_\text{max}=32,\dots,255$ for the $8^{10}$ tensor with a prescribed minimal reduction $f_1^\text{min}=1/2$.
This results in a lower run-time for the $10^9$ tensor as the first step is the most costly part of the algorithm.
As expected, the run-time of the $32^6$ case increases linearly with the maximal rank for $r_\text{max}\ge 16$ and the run-time is significantly higher than for smaller dimensions as the resulting reduction factors are small ($f_j \approx 1/32$).

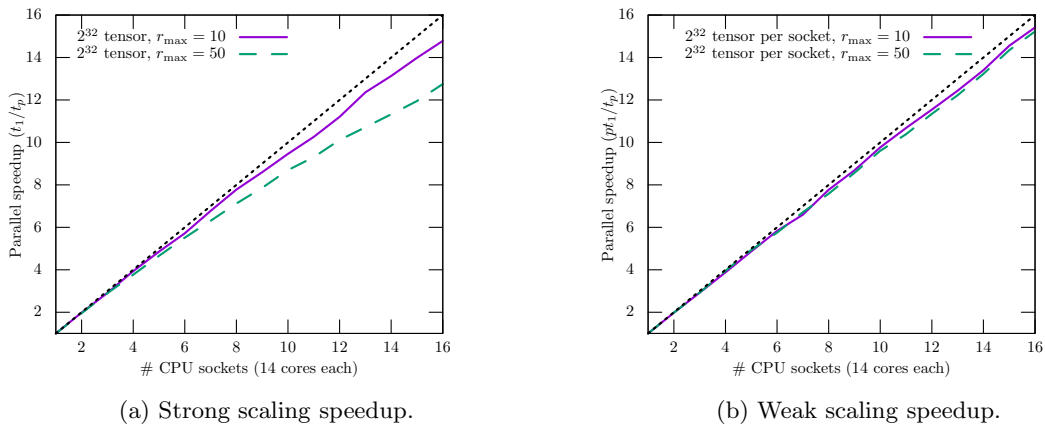
\begin{figure}
\begin{subfigure}[t]{0.47\textwidth}
\begin{tikzpicture}[gnuplot]
\tikzset{every node/.append style={scale=0.60}}
\path (0.000,0.000) rectangle (6.100,5.000);
\gpcolor{color=gp lt color border}
\gpsetlinetype{gp lt border}
\gpsetdashtype{gp dt solid}
\gpsetlinewidth{1.00}
\draw[gp path] (0.680,0.873)--(0.860,0.873);
\draw[gp path] (5.769,0.873)--(5.589,0.873);
\node[gp node right] at (0.570,0.873) {$2$};
\draw[gp path] (0.680,1.436)--(0.860,1.436);
\draw[gp path] (5.769,1.436)--(5.589,1.436);
\node[gp node right] at (0.570,1.436) {$4$};
\draw[gp path] (0.680,1.999)--(0.860,1.999);
\draw[gp path] (5.769,1.999)--(5.589,1.999);
\node[gp node right] at (0.570,1.999) {$6$};
\draw[gp path] (0.680,2.562)--(0.860,2.562);
\draw[gp path] (5.769,2.562)--(5.589,2.562);
\node[gp node right] at (0.570,2.562) {$8$};
\draw[gp path] (0.680,3.125)--(0.860,3.125);
\draw[gp path] (5.769,3.125)--(5.589,3.125);
\node[gp node right] at (0.570,3.125) {$10$};
\draw[gp path] (0.680,3.688)--(0.860,3.688);
\draw[gp path] (5.769,3.688)--(5.589,3.688);
\node[gp node right] at (0.570,3.688) {$12$};
\draw[gp path] (0.680,4.251)--(0.860,4.251);
\draw[gp path] (5.769,4.251)--(5.589,4.251);
\node[gp node right] at (0.570,4.251) {$14$};
\draw[gp path] (0.680,4.814)--(0.860,4.814);
\draw[gp path] (5.769,4.814)--(5.589,4.814);
\node[gp node right] at (0.570,4.814) {$16$};
\draw[gp path] (1.019,0.592)--(1.019,0.772);
\draw[gp path] (1.019,4.814)--(1.019,4.634);
\node[gp node center] at (1.019,0.407) {$2$};
\draw[gp path] (1.698,0.592)--(1.698,0.772);
\draw[gp path] (1.698,4.814)--(1.698,4.634);
\node[gp node center] at (1.698,0.407) {$4$};
\draw[gp path] (2.376,0.592)--(2.376,0.772);
\draw[gp path] (2.376,4.814)--(2.376,4.634);
\node[gp node center] at (2.376,0.407) {$6$};
\draw[gp path] (3.055,0.592)--(3.055,0.772);
\draw[gp path] (3.055,4.814)--(3.055,4.634);
\node[gp node center] at (3.055,0.407) {$8$};
\draw[gp path] (3.733,0.592)--(3.733,0.772);
\draw[gp path] (3.733,4.814)--(3.733,4.634);
\node[gp node center] at (3.733,0.407) {$10$};
\draw[gp path] (4.412,0.592)--(4.412,0.772);
\draw[gp path] (4.412,4.814)--(4.412,4.634);
\node[gp node center] at (4.412,0.407) {$12$};
\draw[gp path] (5.090,0.592)--(5.090,0.772);
\draw[gp path] (5.090,4.814)--(5.090,4.634);
\node[gp node center] at (5.090,0.407) {$14$};
\draw[gp path] (5.769,0.592)--(5.769,0.772);
\draw[gp path] (5.769,4.814)--(5.769,4.634);
\node[gp node center] at (5.769,0.407) {$16$};
\draw[gp path] (0.680,4.814)--(0.680,0.592)--(5.769,0.592)--(5.769,4.814)--cycle;
\node[gp node center,rotate=-270] at (0.175,2.703) {Parallel speedup ($t_1/t_p$)};
\node[gp node center] at (3.224,0.130) {\# CPU sockets (14 cores each)};
\node[gp node right] at (2.990,4.521) {$2^{32}$ tensor, $r_\text{max}=10$};
\gpcolor{rgb color={0.580,0.000,0.827}}
\gpsetdashtype{gp dt 1}
\gpsetlinewidth{2.00}
\draw[gp path] (3.100,4.521)--(3.720,4.521);
\draw[gp path] (0.680,0.592)--(1.019,0.868)--(1.359,1.130)--(1.698,1.418)--(2.037,1.677)%
  --(2.376,1.920)--(2.716,2.217)--(3.055,2.500)--(3.394,2.730)--(3.733,2.973)--(4.073,3.199)%
  --(4.412,3.464)--(4.751,3.789)--(5.090,4.006)--(5.430,4.247)--(5.769,4.471);
\gpcolor{color=gp lt color border}
\node[gp node right] at (2.990,4.296) {$2^{32}$ tensor, $r_\text{max}=50$};
\gpcolor{rgb color={0.000,0.620,0.451}}
\gpsetdashtype{gp dt 2}
\draw[gp path] (3.100,4.296)--(3.720,4.296);
\draw[gp path] (0.680,0.592)--(1.019,0.860)--(1.359,1.126)--(1.698,1.371)--(2.037,1.621)%
  --(2.376,1.862)--(2.716,2.086)--(3.055,2.315)--(3.394,2.523)--(3.733,2.758)--(4.073,2.931)%
  --(4.412,3.157)--(4.751,3.326)--(5.090,3.498)--(5.430,3.674)--(5.769,3.899);
\gpcolor{color=gp lt color border}
\gpsetlinetype{gp lt axes}
\gpsetdashtype{gp dt axes}
\draw[gp path] (0.680,0.592)--(0.731,0.635)--(0.783,0.677)--(0.834,0.720)--(0.886,0.763)%
  --(0.937,0.805)--(0.988,0.848)--(1.040,0.891)--(1.091,0.933)--(1.143,0.976)--(1.194,1.018)%
  --(1.245,1.061)--(1.297,1.104)--(1.348,1.146)--(1.400,1.189)--(1.451,1.232)--(1.502,1.274)%
  --(1.554,1.317)--(1.605,1.360)--(1.657,1.402)--(1.708,1.445)--(1.759,1.488)--(1.811,1.530)%
  --(1.862,1.573)--(1.914,1.616)--(1.965,1.658)--(2.017,1.701)--(2.068,1.743)--(2.119,1.786)%
  --(2.171,1.829)--(2.222,1.871)--(2.274,1.914)--(2.325,1.957)--(2.376,1.999)--(2.428,2.042)%
  --(2.479,2.085)--(2.531,2.127)--(2.582,2.170)--(2.633,2.213)--(2.685,2.255)--(2.736,2.298)%
  --(2.788,2.341)--(2.839,2.383)--(2.890,2.426)--(2.942,2.468)--(2.993,2.511)--(3.045,2.554)%
  --(3.096,2.596)--(3.147,2.639)--(3.199,2.682)--(3.250,2.724)--(3.302,2.767)--(3.353,2.810)%
  --(3.404,2.852)--(3.456,2.895)--(3.507,2.938)--(3.559,2.980)--(3.610,3.023)--(3.661,3.065)%
  --(3.713,3.108)--(3.764,3.151)--(3.816,3.193)--(3.867,3.236)--(3.918,3.279)--(3.970,3.321)%
  --(4.021,3.364)--(4.073,3.407)--(4.124,3.449)--(4.175,3.492)--(4.227,3.535)--(4.278,3.577)%
  --(4.330,3.620)--(4.381,3.663)--(4.432,3.705)--(4.484,3.748)--(4.535,3.790)--(4.587,3.833)%
  --(4.638,3.876)--(4.690,3.918)--(4.741,3.961)--(4.792,4.004)--(4.844,4.046)--(4.895,4.089)%
  --(4.947,4.132)--(4.998,4.174)--(5.049,4.217)--(5.101,4.260)--(5.152,4.302)--(5.204,4.345)%
  --(5.255,4.388)--(5.306,4.430)--(5.358,4.473)--(5.409,4.515)--(5.461,4.558)--(5.512,4.601)%
  --(5.563,4.643)--(5.615,4.686)--(5.666,4.729)--(5.718,4.771)--(5.769,4.814);
\gpsetlinetype{gp lt border}
\gpsetdashtype{gp dt solid}
\gpsetlinewidth{1.00}
\draw[gp path] (0.680,4.814)--(0.680,0.592)--(5.769,0.592)--(5.769,4.814)--cycle;
\gpdefrectangularnode{gp plot 1}{\pgfpoint{0.680cm}{0.592cm}}{\pgfpoint{5.769cm}{4.814cm}}
\end{tikzpicture}
\subcaption{Strong scaling speedup.}
\end{subfigure}
\hfill
\begin{subfigure}[t]{0.47\textwidth}
\begin{tikzpicture}[gnuplot]
\tikzset{every node/.append style={scale=0.60}}
\path (0.000,0.000) rectangle (6.100,5.000);
\gpcolor{color=gp lt color border}
\gpsetlinetype{gp lt border}
\gpsetdashtype{gp dt solid}
\gpsetlinewidth{1.00}
\draw[gp path] (0.680,0.873)--(0.860,0.873);
\draw[gp path] (5.769,0.873)--(5.589,0.873);
\node[gp node right] at (0.570,0.873) {$2$};
\draw[gp path] (0.680,1.436)--(0.860,1.436);
\draw[gp path] (5.769,1.436)--(5.589,1.436);
\node[gp node right] at (0.570,1.436) {$4$};
\draw[gp path] (0.680,1.999)--(0.860,1.999);
\draw[gp path] (5.769,1.999)--(5.589,1.999);
\node[gp node right] at (0.570,1.999) {$6$};
\draw[gp path] (0.680,2.562)--(0.860,2.562);
\draw[gp path] (5.769,2.562)--(5.589,2.562);
\node[gp node right] at (0.570,2.562) {$8$};
\draw[gp path] (0.680,3.125)--(0.860,3.125);
\draw[gp path] (5.769,3.125)--(5.589,3.125);
\node[gp node right] at (0.570,3.125) {$10$};
\draw[gp path] (0.680,3.688)--(0.860,3.688);
\draw[gp path] (5.769,3.688)--(5.589,3.688);
\node[gp node right] at (0.570,3.688) {$12$};
\draw[gp path] (0.680,4.251)--(0.860,4.251);
\draw[gp path] (5.769,4.251)--(5.589,4.251);
\node[gp node right] at (0.570,4.251) {$14$};
\draw[gp path] (0.680,4.814)--(0.860,4.814);
\draw[gp path] (5.769,4.814)--(5.589,4.814);
\node[gp node right] at (0.570,4.814) {$16$};
\draw[gp path] (1.019,0.592)--(1.019,0.772);
\draw[gp path] (1.019,4.814)--(1.019,4.634);
\node[gp node center] at (1.019,0.407) {$2$};
\draw[gp path] (1.698,0.592)--(1.698,0.772);
\draw[gp path] (1.698,4.814)--(1.698,4.634);
\node[gp node center] at (1.698,0.407) {$4$};
\draw[gp path] (2.376,0.592)--(2.376,0.772);
\draw[gp path] (2.376,4.814)--(2.376,4.634);
\node[gp node center] at (2.376,0.407) {$6$};
\draw[gp path] (3.055,0.592)--(3.055,0.772);
\draw[gp path] (3.055,4.814)--(3.055,4.634);
\node[gp node center] at (3.055,0.407) {$8$};
\draw[gp path] (3.733,0.592)--(3.733,0.772);
\draw[gp path] (3.733,4.814)--(3.733,4.634);
\node[gp node center] at (3.733,0.407) {$10$};
\draw[gp path] (4.412,0.592)--(4.412,0.772);
\draw[gp path] (4.412,4.814)--(4.412,4.634);
\node[gp node center] at (4.412,0.407) {$12$};
\draw[gp path] (5.090,0.592)--(5.090,0.772);
\draw[gp path] (5.090,4.814)--(5.090,4.634);
\node[gp node center] at (5.090,0.407) {$14$};
\draw[gp path] (5.769,0.592)--(5.769,0.772);
\draw[gp path] (5.769,4.814)--(5.769,4.634);
\node[gp node center] at (5.769,0.407) {$16$};
\draw[gp path] (0.680,4.814)--(0.680,0.592)--(5.769,0.592)--(5.769,4.814)--cycle;
\node[gp node center,rotate=-270] at (0.175,2.703) {Parallel speedup ($pt_1/t_p$)};
\node[gp node center] at (3.224,0.130) {\# CPU sockets (14 cores each)};
\node[gp node right] at (4.200,4.521) {$2^{32}$ tensor per socket, $r_\text{max}=10$};
\gpcolor{rgb color={0.580,0.000,0.827}}
\gpsetdashtype{gp dt 1}
\gpsetlinewidth{2.00}
\draw[gp path] (4.310,4.521)--(4.930,4.521);
\draw[gp path] (0.680,0.592)--(1.019,0.866)--(1.359,1.133)--(1.698,1.402)--(2.037,1.680)%
  --(2.376,1.951)--(2.716,2.167)--(3.055,2.497)--(3.394,2.755)--(3.733,3.059)--(4.073,3.317)%
  --(4.412,3.560)--(4.751,3.811)--(5.090,4.083)--(5.430,4.407)--(5.769,4.649);
\gpcolor{color=gp lt color border}
\node[gp node right] at (4.200,4.296) {$2^{32}$ tensor per socket, $r_\text{max}=50$};
\gpcolor{rgb color={0.000,0.620,0.451}}
\gpsetdashtype{gp dt 2}
\draw[gp path] (4.310,4.296)--(4.930,4.296);
\draw[gp path] (0.680,0.592)--(1.019,0.868)--(1.359,1.134)--(1.698,1.422)--(2.037,1.699)%
  --(2.376,1.928)--(2.716,2.200)--(3.055,2.445)--(3.394,2.721)--(3.733,3.014)--(4.073,3.234)%
  --(4.412,3.504)--(4.751,3.751)--(5.090,4.034)--(5.430,4.348)--(5.769,4.597);
\gpcolor{color=gp lt color border}
\gpsetlinetype{gp lt axes}
\gpsetdashtype{gp dt axes}
\draw[gp path] (0.680,0.592)--(0.731,0.635)--(0.783,0.677)--(0.834,0.720)--(0.886,0.763)%
  --(0.937,0.805)--(0.988,0.848)--(1.040,0.891)--(1.091,0.933)--(1.143,0.976)--(1.194,1.018)%
  --(1.245,1.061)--(1.297,1.104)--(1.348,1.146)--(1.400,1.189)--(1.451,1.232)--(1.502,1.274)%
  --(1.554,1.317)--(1.605,1.360)--(1.657,1.402)--(1.708,1.445)--(1.759,1.488)--(1.811,1.530)%
  --(1.862,1.573)--(1.914,1.616)--(1.965,1.658)--(2.017,1.701)--(2.068,1.743)--(2.119,1.786)%
  --(2.171,1.829)--(2.222,1.871)--(2.274,1.914)--(2.325,1.957)--(2.376,1.999)--(2.428,2.042)%
  --(2.479,2.085)--(2.531,2.127)--(2.582,2.170)--(2.633,2.213)--(2.685,2.255)--(2.736,2.298)%
  --(2.788,2.341)--(2.839,2.383)--(2.890,2.426)--(2.942,2.468)--(2.993,2.511)--(3.045,2.554)%
  --(3.096,2.596)--(3.147,2.639)--(3.199,2.682)--(3.250,2.724)--(3.302,2.767)--(3.353,2.810)%
  --(3.404,2.852)--(3.456,2.895)--(3.507,2.938)--(3.559,2.980)--(3.610,3.023)--(3.661,3.065)%
  --(3.713,3.108)--(3.764,3.151)--(3.816,3.193)--(3.867,3.236)--(3.918,3.279)--(3.970,3.321)%
  --(4.021,3.364)--(4.073,3.407)--(4.124,3.449)--(4.175,3.492)--(4.227,3.535)--(4.278,3.577)%
  --(4.330,3.620)--(4.381,3.663)--(4.432,3.705)--(4.484,3.748)--(4.535,3.790)--(4.587,3.833)%
  --(4.638,3.876)--(4.690,3.918)--(4.741,3.961)--(4.792,4.004)--(4.844,4.046)--(4.895,4.089)%
  --(4.947,4.132)--(4.998,4.174)--(5.049,4.217)--(5.101,4.260)--(5.152,4.302)--(5.204,4.345)%
  --(5.255,4.388)--(5.306,4.430)--(5.358,4.473)--(5.409,4.515)--(5.461,4.558)--(5.512,4.601)%
  --(5.563,4.643)--(5.615,4.686)--(5.666,4.729)--(5.718,4.771)--(5.769,4.814);
\gpsetlinetype{gp lt border}
\gpsetdashtype{gp dt solid}
\gpsetlinewidth{1.00}
\draw[gp path] (0.680,4.814)--(0.680,0.592)--(5.769,0.592)--(5.769,4.814)--cycle;
\gpdefrectangularnode{gp plot 1}{\pgfpoint{0.680cm}{0.592cm}}{\pgfpoint{5.769cm}{4.814cm}}
\end{tikzpicture}
\subcaption{Weak scaling speedup.}
\end{subfigure}
\caption{Speedup for the TSQR TT-SVD (thick-bounds variant with $f_1^\text{min}=\frac{1}{2}$) on varying number of CPU sockets and nodes with 1 MPI process per socket. Each node has 4~sockets with 14~cores. The reference time is measured on a single socket.}\label{fig:TT-SVD_scaling}
\end{figure}
Finally, we also tested the distributed variant of the TSQR TT-SVD algorithm using MPI.
\autoref{fig:TT-SVD_scaling} shows strong and weak scaling results for input tensors of dimension $2^{32}$ (strong scaling) and $2^{32}$ to $2^{36}$ (weak scaling).
We observe a good weak scaling behavior (parallel efficiency of about $\sim95\%$).
The biggest considered case has an input tensor of size $2^{36}$ ($\sim550$~GByte).
For strong scaling, the problem size per CPU socket gets smaller.
So in particular for bigger TT ranks, the relative overhead due to duplicating the work of the small SVD increases.
The same holds for the relative parallelization overhead in the TSQR algorithm.
The TSQR MPI reduction only amounts to about 3\% of the overall runtime (with 16 CPU sockets and $r_\text{max}=50$, similar for both strong and weak scaling).
Summing up, the distributed variant allows to tackle problems where the dense input tensor is too large for the memory of a single node or where the input tensor is generated by a distributed program on a cluster.
The communication overhead is low.
Only for strong scaling, we observe a significant overhead due to non-parallelized parts of the algorithm.

\section{Conclusion and future work}\label{sec:conclusion}
In this paper we analyzed the node-level performance of the tensor-train SVD algorithm that calculates a low-rank approximation of a high-dimensional tensor.
The results can also be transferred to other non-randomized high-order SVD algorithms.
We considered the case where the input tensor is large and dense, but not too large to be processed completely, i.e., to be read from main memory or disk as a whole.
The theoretical minimal run-time depends on the desired accuracy of the approximation.
For small tensor-train ranks (low accuracy), the algorithm is memory-bound and the ideal run-time on current hardware is approximately twice the time required for reading the data (transferring it from the memory to the CPU).
For larger tensor-train ranks (higher accuracy), the algorithm becomes compute-bound and the ideal run-time increases linearly with the maximal TT-rank.
We presented different variants of the TT-SVD algorithm.
In order to reduce the computational complexity, these variants start with the calculation of the TT-cores at the boundaries of the tensor-train and reduce the data size in each step.
The key ingredient is a Q-less tall-skinny QR decomposition based on Householder reflections that
handles rank-deficient matrices without pivoting by clever use of floating point arithmetic.
We performed numerical experiments with $2^d$ tensors of size up to 550~GByte ($d=36$) on up 224 cores on a small cluster.
Our hybrid-parallel (MPI+OpenMP) TT-SVD implementation achieves almost optimal run-time for small ranks and about 25\% peak performance for larger TT-ranks.
On a single CPU socket, our implementation is about $50\times$ faster compared to TT-SVD algorithms in other libraries.
We provide a lower bound for the run-time: reading the data twice from main memory.
This also indicates that randomized algorithms can cirumvent this lower bound by not considering all data.

For future work, we see three interesting directions:
First, here, we use random input data and prescribe the TT-ranks.
In real applications, usually a certain truncation accuracy is prescribed instead, and the TT-ranks depend on the desired accuracy.
For optimal performance one needs to combine, rearrange or split dimensions based on some heuristic such that the first step leads to a sufficient reduction in data size.
Second, we only analyzed one tensor-train operation for dense input.
Similar performance gains might be possible for other important operations involving large dense data.
Handling sparse input data efficiently is more challenging as the reduction in dimensions in each step does not necessarily lead to a reduction in data size.
And finally, it would be interesting to analyze the performance of randomized decomposition algorithms and to deduce lower bounds for their run-time on current hardware.

\clearpage
\bibliographystyle{siamplain}
\bibliography{almost_optimal_tt_svd}

\begin{thebibliography}{10}

\bibitem{affleck1985}
{\sc I.~Affleck}, {\em Large-$n$ limit of $\mathrm{SU}(n)$ quantum "spin"
  chains}, Physical Review Letters, 54 (1985), pp.~966--969,
  \url{https://doi.org/10.1103/PhysRevLett.54.966}.

\bibitem{affleck1987}
{\sc I.~Affleck, T.~Kennedy, E.~H. Lieb, and H.~Tasaki}, {\em Rigorous results
  on valence-bond ground states in antiferromagnets}, Physical Review Letters,
  59 (1987), pp.~799--802, \url{https://doi.org/10.1103/PhysRevLett.59.799}.

\bibitem{LAPACK}
{\sc E.~Anderson, Z.~Bai, C.~Bischof, L.~S. Blackford, J.~Demmel, J.~Dongarra,
  J.~D. Croz, A.~Greenbaum, S.~Hammarling, A.~McKenney, and D.~Sorensen}, {\em
  {LAPACK} Users' Guide}, Society for Industrial and Applied Mathematics, 1999,
  \url{https://doi.org/10.1137/1.9780898719604}.

\bibitem{Ballard2014}
{\sc G.~Ballard, E.~Carson, J.~Demmel, M.~Hoemmen, N.~Knight, and O.~Schwartz},
  {\em Communication lower bounds and optimal algorithms for numerical linear
  algebra}, Acta Numerica, 23 (2014), pp.~1--155,
  \url{https://doi.org/10.1017/s0962492914000038}.

\bibitem{tnTorch}
{\sc R.~Ballester-Ripoll}, {\em {tntorch - Tensor Network Learning with
  PyTorch}}.
\newblock \url{https://tntorch.readthedocs.io}, 2019.
\newblock Revision 8c81a1cb.

\bibitem{Berry2005}
{\sc M.~W. Berry, S.~A. Pulatova, and G.~W. Stewart}, {\em Algorithm 844},
  {ACM} Transactions on Mathematical Software, 31 (2005), pp.~252--269,
  \url{https://doi.org/10.1145/1067967.1067972}.

\bibitem{Bigoni2016}
{\sc D.~Bigoni, A.~P. Engsig-Karup, and Y.~M. Marzouk}, {\em Spectral
  tensor-train decomposition}, SIAM Journal on Scientific Computing, 38 (2016),
  pp.~A2405--A2439, \url{https://doi.org/10.1137/15m1036919}.

\bibitem{Carson2016}
{\sc E.~Carson, J.~Demmel, L.~Grigori, N.~Knight, P.~Koanantakool, O.~Schwartz,
  and H.~V. Simhadri}, {\em Write-avoiding algorithms}, in 2016 {IEEE}
  International Parallel and Distributed Processing Symposium ({IPDPS}),
  {IEEE}, may 2016, pp.~648--658, \url{https://doi.org/10.1109/ipdps.2016.114}.

\bibitem{Chen2019}
{\sc C.~Chen, K.~Batselier, C.-Y. Ko, and N.~Wong}, {\em A support tensor train
  machine}, in 2019 International Joint Conference on Neural Networks
  ({IJCNN}), {IEEE}, Jul 2019,
  \url{https://doi.org/10.1109/ijcnn.2019.8851985}.

\bibitem{Chen2022}
{\sc C.~Chen, K.~Batselier, W.~Yu, and N.~Wong}, {\em Kernelized support tensor
  train machines}, Pattern Recognition, 122 (2022), p.~108337,
  \url{https://doi.org/10.1016/j.patcog.2021.108337}.

\bibitem{Constantine2014}
{\sc P.~G. Constantine, D.~F. Gleich, Y.~Hou, and J.~Templeton}, {\em Model
  reduction with {MapReduce}-enabled tall and skinny singular value
  decomposition}, SIAM Journal on Scientific Computing, 36 (2014),
  pp.~S166--S191, \url{https://doi.org/10.1137/130925219}.

\bibitem{Daas2022}
{\sc H.~A. Daas, G.~Ballard, and P.~Benner}, {\em Parallel algorithms for
  tensor train arithmetic}, {SIAM} Journal on Scientific Computing, 44 (2022),
  pp.~C25--C53, \url{https://doi.org/10.1137/20m1387158}.

\bibitem{Demmel2012}
{\sc J.~Demmel, L.~Grigori, M.~Hoemmen, and J.~Langou}, {\em
  Communication-optimal parallel and sequential {QR} and {LU} factorizations},
  SIAM Journal on Scientific Computing, 34 (2012), pp.~A206--A239,
  \url{https://doi.org/10.1137/080731992}.

\bibitem{Demmel2015}
{\sc J.~W. Demmel, L.~Grigori, M.~Gu, and H.~Xiang}, {\em Communication
  avoiding rank revealing {QR} factorization with column pivoting}, SIAM
  Journal on Matrix Analysis and Applications, 36 (2015), pp.~55--89,
  \url{https://doi.org/10.1137/13092157x}.

\bibitem{Dolgov2013}
{\sc S.~Dolgov and B.~Khoromskij}, {\em Two-level {QTT}-tucker format for
  optimized tensor calculus}, SIAM Journal on Matrix Analysis and Applications,
  34 (2013), pp.~593--623, \url{https://doi.org/10.1137/120882597}.

\bibitem{Fan2016}
{\sc H.-Y. Fan, L.~Zhang, E.~w.~Chu, and Y.~Wei}, {\em Q-less {QR}
  decomposition in inner product spaces}, Linear Algebra and Its Applications,
  491 (2016), pp.~292--316, \url{https://doi.org/10.1016/j.laa.2015.08.035}.

\bibitem{Gelss2019}
{\sc P.~Gelß, S.~Klus, M.~Scherer, F.~Nüske, and M.~Lücke}, {\em {Scikit-TT}
  tensor-train computations in python}.
\newblock \url{https://github.com/PGelss/scikit_tt}, 2019.
\newblock Revision 1dfd64a.

\bibitem{Grasedyck2011}
{\sc L.~Grasedyck and W.~Hackbusch}, {\em An introduction to hierarchical (h-)
  rank and {TT}-rank of tensors with examples}, Computational Methods in
  Applied Mathematics, 11 (2011), pp.~291--304,
  \url{https://doi.org/10.2478/cmam-2011-0016}.

\bibitem{Grasedyck2013}
{\sc L.~Grasedyck, D.~Kressner, and C.~Tobler}, {\em A literature survey of
  low-rank tensor approximation techniques}, GAMM Mitteilungen, 36 (2013),
  pp.~53--78, \url{https://doi.org/10.1002/gamm.201310004}.

\bibitem{Eigen}
{\sc G.~Guennebaud, B.~Jacob, et~al.}, {\em Eigen v3}.
\newblock \url{http://eigen.tuxfamily.org}, 2010.
\newblock Version 3.3.9.

\bibitem{Hager2010}
{\sc G.~Hager and G.~Wellein}, {\em Introduction to High Performance Computing
  for Scientists and Engineers}, {CRC} Press, jul 2010,
  \url{https://doi.org/10.1201/ebk1439811924}.

\bibitem{Harris2020}
{\sc C.~R. Harris, K.~J. Millman, S.~J. van~der Walt, R.~Gommers, P.~Virtanen,
  D.~Cournapeau, E.~Wieser, J.~Taylor, S.~Berg, N.~J. Smith, R.~Kern, M.~Picus,
  S.~Hoyer, M.~H. van Kerkwijk, M.~Brett, A.~Haldane, J.~F. del R{\'{\i}}o,
  M.~Wiebe, P.~Peterson, P.~G{\'{e}}rard-Marchant, K.~Sheppard, T.~Reddy,
  W.~Weckesser, H.~Abbasi, C.~Gohlke, and T.~E. Oliphant}, {\em Array
  programming with {NumPy}}, Nature, 585 (2020), pp.~357--362,
  \url{https://doi.org/10.1038/s41586-020-2649-2}.

\bibitem{Householder1958}
{\sc A.~S. Householder}, {\em Unitary triangularization of a nonsymmetric
  matrix}, Journal of the ACM, 5 (1958), pp.~339--342,
  \url{https://doi.org/10.1145/320941.320947}.

\bibitem{IntelMKL}
{\sc Intel}, {\em Math kernel library (mkl)}.
\newblock \url{https://software.intel.com/mkl}.
\newblock Version 2020.3.

\bibitem{Khoromskij2011}
{\sc B.~N. Khoromskij}, {\em {$O(d\log N)$}-quantics approximation of n-d
  tensors in high-dimensional numerical modeling}, Constructive Approximation,
  34 (2011), pp.~257--280, \url{https://doi.org/10.1007/s00365-011-9131-1}.

\bibitem{Klus2019}
{\sc S.~Klus and P.~Gel{\ss}}, {\em Tensor-based algorithms for image
  classification}, Algorithms, 12 (2019), p.~240,
  \url{https://doi.org/10.3390/a12110240}.

\bibitem{Kolda2020}
{\sc T.~G. Kolda and D.~Hong}, {\em Stochastic gradients for large-scale tensor
  decomposition}, SIAM Journal on Mathematics of Data Science, 2 (2020),
  pp.~1066--1095, \url{https://doi.org/10.1137/19m1266265}.

\bibitem{Kour2020}
{\sc K.~Kour, S.~Dolgov, M.~Stoll, and P.~Benner}, {\em Efficient
  structure-preserving support tensor train machine},  (2020),
  \url{https://doi.org/10.48550/ARXIV.2002.05079}.

\bibitem{Larsen2020}
{\sc B.~W. Larsen and T.~G. Kolda}, {\em Practical leverage-based sampling for
  low-rank tensor decomposition}, 2020,
  \url{https://doi.org/10.48550/ARXIV.2006.16438}.

\bibitem{McCalpin1995}
{\sc J.~D. McCalpin}, {\em Memory bandwidth and machine balance in current high
  performance computers}, IEEE Computer Society Technical Committee on Computer
  Architecture (TCCA) Newsletter,  (1995), pp.~19--25.

\bibitem{Novikov2020}
{\sc A.~Novikov, P.~Izmailov, V.~Khrulkov, M.~Figurnov, and I.~Oseledets}, {\em
  {Tensor Train Decomposition on TensorFlow (T3F)}}, Journal of Machine
  Learning Research, 21 (2020), pp.~1--7,
  \url{http://jmlr.org/papers/v21/18-008.html}.

\bibitem{oseledets2009a}
{\sc I.~V. Oseledets}, {\em A new tensor decomposition}, Doklady Mathematics,
  80 (2009), pp.~495--496, \url{https://doi.org/10.1134/S1064562409040115}.

\bibitem{Oseledets2011}
{\sc I.~V. Oseledets}, {\em {Tensor-Train Decomposition}}, SIAM Journal on
  Scientific Computing, 33 (2011), pp.~2295--2317,
  \url{https://doi.org/10.1137/090752286}.

\bibitem{oseledets2009b}
{\sc I.~V. Oseledets and E.~E. Tyrtyshnikov}, {\em Breaking the curse of
  dimensionality, or how to use svd in many dimensions}, SIAM Journal on
  Scientific Computing, 31 (2009), pp.~3744--3759,
  \url{https://doi.org/10.1137/090748330}.

\bibitem{Psarras2019}
{\sc C.~Psarras, H.~Barthels, and P.~Bientinesi}, {\em The linear algebra
  mapping problem. current state of linear algebra languages and libraries},
  (2019), \url{https://doi.org/10.48550/ARXIV.1911.09421}.

\bibitem{RoehrigZoellner2021pitts}
{\sc M.~Röhrig-Zöllner}, {\em {PITTS - Parallel Iterative Tensor-Train
  Solvers}}, 9 2021, \url{https://doi.org/10.5281/zenodo.5534544},
  \url{https://github.com/melven/pitts}.

\bibitem{RoehrigZoellner2015}
{\sc M.~Röhrig-Zöllner, J.~Thies, M.~Kreutzer, A.~Alvermann, A.~Pieper,
  A.~Basermann, G.~Hager, G.~Wellein, and H.~Fehske}, {\em Increasing the
  performance of the jacobi--davidson method by blocking}, {SIAM} Journal on
  Scientific Computing, 37 (2015), pp.~C697--C722,
  \url{https://doi.org/10.1137/140976017}.

\bibitem{Stengel2015}
{\sc H.~Stengel, J.~Treibig, G.~Hager, and G.~Wellein}, {\em Quantifying
  performance bottlenecks of stencil computations using the
  execution-cache-memory model}, in Proceedings of the 29th {ACM} on
  International Conference on Supercomputing - {ICS} '15, {ACM} Press, 2015,
  \url{https://doi.org/10.1145/2751205.2751240}.

\bibitem{Treibig2010}
{\sc J.~Treibig, G.~Hager, and G.~Wellein}, {\em {LIKWID}: A lightweight
  performance-oriented tool suite for x86 multicore environments}, in 2010 39th
  International Conference on Parallel Processing Workshops, {IEEE}, Sept.
  2010, \url{https://doi.org/10.1109/icppw.2010.38}.

\bibitem{Trilinos}
{\sc T.~{T}rilinos~{P}roject {T}eam}, {\em The {T}rilinos {P}roject {W}ebsite}.
\newblock \url{https://trilinos.github.io}, 2020.
\newblock Accessed May 22, 2020.

\bibitem{verstraete2006mps}
{\sc F.~Verstraete and J.~I. Cirac}, {\em Matrix product states represent
  ground states faithfully}, Physical Review B, 73 (2006), p.~094423,
  \url{https://doi.org/10.1103/PhysRevB.73.094423}.

\bibitem{white1992dmrg}
{\sc S.~R. White}, {\em Density matrix formulation for quantum renormalization
  groups}, Physical Review Letters, 69 (1992), pp.~2863--2866,
  \url{https://doi.org/10.1103/PhysRevLett.69.2863}.

\bibitem{Williams2009}
{\sc S.~Williams, A.~Waterman, and D.~Patterson}, {\em Roofline: An insightful
  visual performance model for multicore architectures}, Communications of the
  ACM, 52 (2009), pp.~65--76, \url{https://doi.org/10.1145/1498765.1498785}.

\end{thebibliography}

\end{document}